\documentclass{amsart}
\usepackage{amssymb,amsfonts,amscd,graphicx,psfrag}
\usepackage[all]{xy}
\numberwithin{equation}{section}
\let\cal\mathcal

\def\Lscr{{\cal L}}

\def\Oscr{{\cal O}}

\def\Rscr{{\cal R}}
\def\Sscr{{\cal S}}

\let\blb\mathbb

\def \PP{{\blb P}}

\def \ZZ{{\blb Z}}

\def \NN{{\blb N}}

\def\ev{\operatorname{ev}}

\let\st\ast
\let\at\ast
\def\Der{\operatorname{Der}}
\def\Inn{\operatorname{Inn}}

\def\quot{/\!\!/}

\def\Lie{\operatorname{Lie}}

\def\Gl{\operatorname {Gl}}
\def\Rep{\operatorname {Rep}}

\def\Hom{\operatorname {Hom}}

\def\tr{\operatorname {tr}}
\def\Tr{\operatorname {Tr}}

\def\ker{\operatorname {ker}}

\def\r{\rightarrow}

\def\DR{\operatorname{DR}}

\def\inv{\operatorname{inv}}
\def\coinv{\operatorname{coinv}}

\newtheorem{lemma}{Lemma}[section]
\newtheorem{proposition}[lemma]{Proposition}
\newtheorem{theorem}[lemma]{Theorem}

\newtheorem{lemmas}{Lemma}[subsection]
\newtheorem{propositions}[lemmas]{Proposition}
\newtheorem{theorems}[lemmas]{Theorem}
\newtheorem{corollarys}[lemmas]{Corollary}

\theoremstyle{definition}

\newtheorem{example}[lemma]{Example}

\newtheorem{examples}[lemmas]{Example}
\newtheorem{definitions}[lemmas]{Definition}

{

}

\theoremstyle{remark}

\newtheorem{remarks}[lemmas]{Remark}

\newdimen\uboxsep \uboxsep=1ex
\def\uboxn#1{\vtop to 0pt{\hrule height 0pt depth 0pt\vskip\uboxsep
\hbox to 0pt{\hss #1\hss}\vss}}

\def\uboxs#1{\vbox to 0pt{\vss\hbox to 0pt{\hss #1\hss}
\vskip\uboxsep\hrule height 0pt depth 0pt}}

\def\ldb{\mathopen{\{\!\!\{}} \def\rdb{\mathclose{\}\!\!\}}}
 
\def\ldbgg{\mathopen{\Bigl\{\!\!\Bigl\{}} \def\rdbgg{\mathclose{\Bigr\}\!\!\Bigr\}}}

\title{Double Poisson algebras}
\author{Michel Van den Bergh}
 \address{Departement WNI,  Limburgs Universitair Centrum, 
 3590 Diepenbeek, Belgium.}
  \email{michel.vandenbergh@uhasselt.be}
\thanks{The author is a senior researcher at the FWO}
\begin{document}
\begin{abstract} 
  In this paper we develop Poisson geometry for
  non-com\-mu\-ta\-tive algebras. This generalizes the bi-symplectic
  geometry which was recently, and independently, introduced by
  Crawley-Boevey, Etingof and Ginzburg.

Our (quasi-)Poisson
brackets induce classical (quasi-)Poisson brackets on representation spaces. 
As an application we show that the moduli spaces of representations
  associated to the deformed multiplicative preprojective algebras
  recently introduced by Crawley-Boevey and Shaw carry a natural
  Poisson structure. 
 \end{abstract}
\maketitle
\setcounter{tocdepth}{1}
\tableofcontents
\section{Introduction}
In this introduction we assume that $k$ is an algebraically closed
field of characteristic zero. We start with our original motivating
example (taken from \cite{CBShaw}).
Let $Q=(Q,I,h,t)$ be a finite quiver with vertex set
$I=\{1,\ldots,n\}$ and edge set $Q$. The maps $t,h:Q\r I$ associate
with every edge its starting and ending vertex.  We let $\bar{Q}$ be the double of
$Q$. $\bar{Q}$ is obtained from $Q$ by adjoining for every arrow $a$
an opposite arrow $a^\ast$. We define $\epsilon:\bar{Q}\r \{\pm 1\}$
as the function which is $1$ on $Q$ and $-1$ on $\bar{Q}-Q$.

Let $\alpha=(\alpha_1,\ldots,\alpha_n)\in \NN^n$ be a dimension vector and fix
scalars $q=(q_1,\ldots,q_n)\in (k^\ast)^n$. Put
$\Rscr_\alpha=\prod_{a\in \bar{Q}} M_{\alpha_{t(a)}\times
  \alpha_{h(a)}}$.  The group $\Gl_\alpha=\prod_i \Gl_{\alpha_i}$ acts on $\Rscr_\alpha$
 by conjugation.

Let $\Sscr_{\alpha,q}$ be the $\Gl_\alpha$ invariant subscheme of $\Rscr_\alpha$
consisting of matrices $(X_a)_{a\in \bar{Q}}$ such that $1+X_a
X_{a^\ast}$ is invertible for all $a\in \bar{Q}$ and such that the
following equations are satisfied for all $i\in I$.
\[
\prod_{a\in\bar{Q},h(a)=i} (1+X_a X_{a^\ast})^{\epsilon(a)}=q_i
\]
(it is shown in \cite{CBShaw} that $\Sscr_{\alpha,q}$ is independent of the ordering on
these products). We prove the following result in this paper.
\begin{theorem} 
\label{ref-1.1-0}
The GIT quotient $\Sscr_{\alpha,q}\quot \Gl_\alpha$ is in a natural
way a Poisson variety.
\end{theorem}
This result is not unexpected since if $Q$ is a ``star'' then it is
shown in \cite{CBShaw} that the points in $\Sscr_{\alpha,q}\quot
\Gl_\alpha$ correspond to local systems on $\PP^1$ whose monodromy
lies in the closure of specific conjugacy classes. The result then
follows from the work of Atiyah and Bott \cite{AB1}.

Different proofs of the Atiyah-Bott result were given in \cite{AKM, AMM} using 
quasi-Hamilto\-nian reduction and fusion. 
It is possible to give
a proof of Theorem  \ref{ref-1.1-0} in the same spirit. 
First one considers
the small quiver  consisting of two vertices and two arrows $a,a^\ast$
and one shows that in that case $\Rscr_\alpha$ is quasi-Hamiltonian. This is then
extended to general quivers using a process called ``fusion''. Finally we obtain
a Poisson structure on $\Sscr_{\alpha,q}$ by quasi-Hamiltonian reduction.

While working out this proof I noticed that all computations could 
be done directly in the path algebra $k\bar{Q}$ of $\bar{Q}$
(suitably localized).  If computations are organized this way explicit
matrices occur, somewhat as an afterthought, only in the very last step.  Trying to understand why
this is so then became the second motivation for writing this
paper.

\medskip

So we restart this introduction! Throughout $A$ is a $k$-algebra which for simplicity
we assume to be finitely generated. For $N\in \NN$ the 
associated representation space of $A$ is defined as
\[
\Rep(A,N)=\Hom(A,M_N(k))
\]
The group $\Gl_{{N}}$ acts on $\Rep(A,{{N}})$ by conjugation on $M_{{N}}(k)$.

A well-known philosophy in non-commutative algebraic geometry
(probably first formulated by Maxim Kontsevich) is that for a property of
the non-commutative ring $A$ to have geometric meaning it should
induce standard geometric properties  on all $\Rep(A,{{N}})$. The case of
symplectic geometry was worked out in \cite{LebBock,Ginzburg1,Kosymp}.
In this paper we discuss Poisson geometry.  More precisely we work out
what kind of structure we need on $A$ in order that all $\Rep(A,{{N}})$
are Poisson varieties.

To motivate our definitions we have to look in more detail at the
coordinate ring $\Oscr(\Rep(A,{{N}}))$ of $\Rep(A,{{N}})$.  For every $a\in A$
we have a corresponding matrix valued function
$(a_{ij})_{i,j=1,\ldots,N}$ on $\Rep(A,{{N}})$. It is easy to see that the
ring $\Oscr(\Rep(A,{{N}}))$ is generated by the functions $a_{ij}$,
subject to the relations
\[
(ab)_{ij}=a_{il}b_{lj}
\]
(where here and below we sum over repeated indices). Hence to define 
a Poisson bracket $\{-,-\}$ on $\Rep(A,{{N}})$ we have to fix the values
of $\{a_{ij},b_{uv}\}$ for all $a,b\in A$. Now $\{a_{ij},b_{uv}\}$ depends
on four indices so it is natural to assume that it comes from an element
of $A\otimes A$. This leads to the following definition. A \emph{double
bracket} on $A$ is a bilinear map 
\[
\ldb-,-\rdb:A\times A\r A\otimes A
\]
which is a derivation in its second argument (for the outer bimodule
structure on $A$) and which satisfies
\[
\ldb a,b\rdb =-\ldb b,a\rdb^\circ
\]
where $(u\otimes v)^\circ=v\otimes u$. We say that $A$ is a double Poisson
algebra if $\ldb-,-\rdb$ satisfies in addition a natural analog of the Jacobi 
identity (see \S\ref{ref-2.3-9}). A special case of  
one of our results is  the following (see \S\ref{ref-7.5-102}). 
\begin{proposition} If $A,\ldb-,-\rdb$ is a double Poisson algebra then
$\Oscr(\Rep(A,n))$ is a Poisson algebra, with Poisson bracket given by
\begin{equation}
\label{ref-1.1-1}
\{ a_{ij},b_{uv}\}=\ldb a,b\rdb'_{uj} \ldb a,b\rdb''_{iv}
\end{equation}
where by convention we write an element $x$ of $A\otimes A$ as $x'\otimes x''$
(i.e.\ we drop the summation sign).
\end{proposition}
\begin{example} (see Examples \ref{ref-2.3.3-15},\ref{ref-7.5.3-106}) If $A=k[t]/(t^n)$ then $\Rep(A,n)$ consists nilpotent $n\times n$ matrices.
$A$ has a double Poisson
bracket  which is uniquely defined by the property.
\[
\ldb t,t\rdb =t\otimes 1-1\otimes t
\]
This double bracket induces the standard Poisson bracket on nilpotent matrices.
\end{example}
Let $\ldb-,-\rdb$ be a double Poisson bracket on $A$. We define the associated
bracket as 
\[
\{-,-\}:A\times A\r A: (a,b)\mapsto \ldb a,b\rdb'\ldb a,b\rdb''
\]
\begin{proposition}
\label{ref-1.4-2} Assume that $A,\ldb-,-\rdb$ is a double Poisson algebra.
Then the following holds
\begin{enumerate}
\item $\{-,-\}$ is a derivation in its second argument and vanishes
on commutators in it is first argument.
\item $\{-,-\}$ is anti-symmetric modulo commutators. 
\item $\{-,-\}$ makes $A$ into a left Loday algebra \cite{KS1,Loday}. I.e.\ 
$\{-,-\}$ satisfies the following version of the Jacobi identity
\[
\{a,\{b,c\}\}=\{\{a,b\} ,c\} +\{b,\{a,c\}\}
\]
\item $\{-,-\}$ makes $A/[A,A]$ into a Lie algebra.
\end{enumerate}
\end{proposition}

In commutative geometry it is customary to describe a Poisson bracket
on a smooth variety $X$ in terms of a bivector field, i.e. in terms of
a section $P$ of $\bigwedge^2 T_{X}$ satisfying $\{P,P\}=0$ where $\{-,-\}$ is
the so-called \emph{Schouten-Nijenhuis bracket} on $\Gamma(X,\bigwedge T_X)$.
Our next aim is to give non-commutative version of this. 

In the rest of this introduction we assume for simplicity that $A$ is \emph{smooth} by
which we mean that $A$ is finitely generated and
$\Omega_A=\ker(A\otimes A\r A)$ is a 
projective $A$-bimodule.  It is easy to see that this implies that all spaces
$\Rep(A,N)$ are smooth over $k$.

We first have to find the correct non-commutative analogue of a vector field.
There are in fact \emph{two} good answers to this. If we insist
that a vector field on $A$  induces vector fields on all $\Rep(A,{{N}})$ 
then a vector field on $A$ should simply be a derivation $\Delta:A\r A$. 
The induced derivation $\delta$ on $\Oscr(\Rep(A,{{N}}))$ is then given by
\[
\delta(a_{ij})=\Delta(a)_{ij}
\]
A second point of view is that a vector field $\Delta$ on $A$ should induce
\emph{matrix valued vector fields} $(\Delta_{ij})_{i,j=1,\ldots,n}$ on all $\Rep(A,{{N}})$. Since now 
$\Delta_{ij}(a_{uv})$ depends on four indices $\Delta(a)$ should be an element
of $A\otimes A$. 

In this paper we accept the second point of view, i.e.\ 
 vector fields on $A$ will be elements of  $D_{A}\overset{\text{def}}{=}\Der(A,A\otimes A)$ where 
as usual we put the outer bimodule structure on $A\otimes A$.  The 
corresponding matrix valued vector fields on  $\Rep(A,{{N}})$ are then given
by
\[
\Delta_{ij}(a_{uv})=\Delta(a)'_{uj}\Delta(a)''_{iv}
\]
$D_A$ contains a remarkable element $E$ which acts as $E(a)=a\otimes
1-1\otimes a$. We will call this element the \emph{gauge element} since we have
\begin{proposition} The matrix valued vector field $(E_{ji})_{ij}$ on
$\Rep(A,n)$ is the derivative of the action of $\Gl_{{N}}$ by conjugation.
\end{proposition}
The importance
of $\Der(A,A\otimes A)$ was first emphasized in \cite{CB3}.

Starting with $D_A$ we define \emph{the algebra of p{{oly-vector}} fields} $DA$ 
on $A$ as the tensor algebra $T_A D_A$ of $D_A$ where we make $D_A$ into an
$A$-bimodule by using the inner bimodule structure on $A\otimes A$.

Another main result of this paper is  (see \S\ref{ref-3.2-44})
\begin{proposition}  The graded algebra
$DA$ has the structure of a \emph{double Gerstenhaber algebra}
  i.e.\ a (super) double Poisson algebra with a double Poisson bracket
$\ldb-,-\rdb$ of degree $-1$. 
 \end{proposition}
We call $\ldb-,-\rdb$ the \emph{Schouten-Nijenhuis bracket} on $DA$. It is
somewhat hard to construct, but as we will see below, in the case
of quivers it takes a very trivial form. 

The elements
of $DA$ define matrix valued p{{oly-vector}} fields on $\Rep(A,N)$ by the rule
\[
(\delta_1\cdots \delta_m)_{ij}=\delta_{1,il_1}\delta_{1,l_1l_2}\cdots \delta_{m,l_mj}
\]
The compatibility between the matrix valued p{{oly-vector}} fields and the
Schouten brackets on $DA$ and $\Gamma(\Rep(A,N),T_{\Rep(A,N)})$ is
given by a formula which is entirely similar to \eqref{ref-1.1-1}
\[
\{P_{ij},Q_{uv}\}=\ldb P,Q\rdb'_{uj}\ldb P,Q\rdb''_{iv}
\]
Let us write $\tr(P)=P_{ii}$. Then the previous formula yields a morphism 
of graded Lie algebras
\[
\tr:DA/[DA,DA]\r \Gamma(\Rep(A,N),\wedge T_{\Rep(A,N)})^{\Gl_N}
\]

To reconnect with double Poisson structures on $A$ we show that there is a bijection
\[
\bigl(DA/[DA,DA]\bigr)_2\leftrightarrow \{\text{double brackets on $A$}\}
\]
which sends $\delta_1\delta_2$ for  $\delta_1,\delta_2\in D_A$ to the
double bracket
\[
\ldb a,b\rdb =\delta_2(b)'\delta_1(a)''\otimes \delta_1(a)'\delta_2(b)''
-\delta_1(b)'\delta_2(a)''\otimes \delta_2(a)'\delta_1(b)''
\]
An element $P\in (DA)_2$ corresponds to a double Poisson bracket if
and only if 
\[
\{P,P\}=0\qquad\text{modulo commutators}
\]

Having a rudimentary differential geometric formalism in place we can now
define various related notions. For example we say that $\mu\in A$ is a
\emph{moment map} for a double Poisson bracket $P$ if the following identity 
holds in $D_A$:
\[
\{P,\mu\}=-E
\]
The reason is of course that if $\mu\in A$ then the corresponding 
matrix valued function $(\mu_{ij})_{ij}$ defines a moment map
$\Rep(A,N)\r M_{{N}}$ for the action of $\Gl_{{N}}$ on $\Rep(A,{{N}})$ (where
we identify, as is customary, $M_{{N}}$ with its dual through the trace map).

We can also define the corresponding multiplicative notions (see \cite{AKM}).
An element $P\in (DA)_2$ is said to be a \emph{quasi-Poisson bracket} if
the following identity holds
\[
\{P,P\}=\frac{1}{6} E^3\qquad \text{modulo commutators}
\]
and an element $\Phi\in A$ is a multiplicative moment map for $P$ if it is a unit and
\[
\{P,\Phi\}=-\frac{1}{2}(E\Phi+\Phi E)
\]
in $D_A$. 
Again these notions induce the corresponding notions on representation spaces. 
\begin{proposition}
\label{ref-1.7-3}
\begin{enumerate}
\item Proposition \ref{ref-1.4-2} goes through unmodified for double quasi-Poisson algebras.
\item Let $A$ be either a double Poisson algebra with moment map $\mu$
  or a double quasi-Poisson algebra with moment map $\Phi$. Put
  $A^\lambda=A/(\mu-\lambda)$ in the first case with $\lambda\in k$ and
  $A^q=A/(\Phi-q)$ with $q\in k^\ast$ in the second case. Then the
  associated (quasi-)Poisson brackets on $\Oscr(\Rep(A,{{N}}))$ induces
  Poisson brackets on $\Oscr(\Rep(A,{{N}}))^{\Gl_{{N}}}$ and 
$\Oscr(\Rep(A^{\text{$\lambda$ or $q$}},{{N}}))^{\Gl_{{N}}}$.
\end{enumerate}
\end{proposition}
The second part of this theorem is an application of (quasi) Hamiltonian
reduction \cite{AKM}.

\medskip

Now we discuss quivers. Thus we return to the setting in the beginning
of this introduction. In order for things to work nicely we must set things up
in a relative setting.  I.e.\ we let $B$ a fixed commutative semi-simple
algebra of the form $ke_1\oplus\cdots \oplus ke_n$ with $e_i^2=e_i$. A $B$ algebra is a
$k$-algebra $A$ equipped with a morphism of $k$-algebras $B\r A$. For
$B$-algebras we may define relative versions of the notions introduced
above e.g.\ $D_{A/B}=\Der_B(A,A\otimes A)$, $D_B A=T_A D_{A/B}$. Representation
spaces are now indexed by an $n$-tuples $(\alpha_1,\ldots,\alpha_n)\in \NN^n$. By definition
\[
\Rep(A,\alpha)=\Hom_B(A,M_N(k))
\]
where $N=\alpha_1+\cdots+\alpha_n$ and we view $B$ as being diagonally embedded 
$M_N(k)$.

Now let $A=kQ$. In this case the idempotents $e_i$ are the paths of length
zero corresponding to the vertices of $Q$. 
 For $a\in Q$ we define the element $\frac{\partial\ }{\partial a}\in D_BA$ 
which on  $b\in Q$ acts as
\[
\frac{\partial b}{\partial a}
=
\begin{cases}
e_{t(a)}\otimes e_{h(a)}&\text{if $a=b$}\\
0&\text{otherwise}
\end{cases}
\]
It is clear that $D_{A/B}$ is generated by $\left(\frac{\partial\
  }{\partial a}\right)_{a\in Q}$ as an $A$-bimodule. Hence $D_BA$ is the
tensor algebra over $A$ generated by $\left(\frac{\partial\ }{\partial
    a}\right)_a$. The matrix valued vector field corresponding to $\frac{\partial\ }{\partial
    a}$ is given by
\[
\left(\frac{\partial\ }{\partial a}
\right)_{ij}=
\begin{cases}
\displaystyle \frac{\partial\ }{\partial a_{ji}}&\text{\rm if $\phi(i)=h(a)$, $\phi(j)=t(a)$}\\
0&\text{\rm otherwise}
\end{cases}
\]
where for $i\in \{1,\ldots,N\}$ we put $\phi(i)=p\in\{1,\ldots, n\}$ 
if $i$ is in the  $p$'th subinterval
of $[1\,.\,.\,N]$ when we decompose the latter into  intervals of length $(\alpha_p)_p$.

The Schouten bracket on $A$ is as follows. Let  $a,b\in Q$. Then
\begin{align*}
\ldb a,b\rdb &=0 \\
\ldbgg \frac{\partial\ }{\partial a},b\rdbgg&=
\begin{cases}
e_{t(a)}\otimes e_{h(a)}&\text{if $a=b$}\\
0&\text{otherwise}
\end{cases}\\
\ldbgg \frac{\partial\ }{\partial a},\frac{\partial\ }{\partial b}\rdbgg&=0
\end{align*}
We prove 
\begin{theorem} (Theorem \ref{ref-6.3.1-88} in the body of the paper) $A=k\bar{Q}$ has a
double Poisson bracket given by
\[
P=\sum_{a\in Q} \frac{\partial\ }{\partial a}
\frac{\partial\ }{\partial a^\ast}
\]
and a corresponding moment map
\[
\mu=\sum_{a\in Q} [a,a^\ast]
\]
\end{theorem}
This theorem is more or less a reformulation of known results. The
induced Lie algebra structure on $k\bar{Q}/[k\bar{Q},k\bar{Q}]$ is the
so-called \emph{necklace Lie algebra} \cite{LebBock,Ginzburg1,Kosymp}.
However it is noteworthy that this Lie algebra structure is induced
from a Loday algebra structure on $k\bar{Q}$.  In
\S\ref{ref-6.4-90} we work out what it is.

 The algebra
$A^\lambda$ introduced in Proposition \ref{ref-1.7-3}
is the so-called \emph{deformed preprojective algebra} \cite{CH} $\Pi^\lambda$.
The Poisson bracket on $\Rep(\Pi^\lambda,\alpha)/Gl(\alpha)$ is obtained from the standard
Poisson bracket on $\Rscr_\alpha=\Rep(k\bar{Q},\alpha)$ given by
\[
\sum_{a\in Q} \frac{\partial\ }{\partial (X_a)_{ij}}
\frac{\partial\ }{\partial (X_{a^\ast})_{ji}}
\]
in the notations of the first paragraph.

We then prove the main result of this paper.
\begin{theorem} 
(Theorem \ref{ref-6.7.1-96})
Let $A$ be obtained from $k\bar{Q}$  by inverting all elements
$(1+aa^\ast)_{a\in \bar{Q}}$. Fix an arbitrary total ordering
  on $\bar{Q}$.  Then $A$ has a quasi-Poisson bracket given by
\[
P=\frac{1}{2}\left(\sum_{a\in \bar{Q}}
\left(\epsilon(a)
(1
+
a^\ast a)\frac{\partial\ }{\partial a} \frac{\partial\ }{\partial a^\ast}\right)
-
\sum_{a<b\in\bar{Q}}\left(\frac{\partial }{\partial  a^\ast} a^\ast-a\frac{\partial }{\partial  a}\right)\left(
\frac{\partial }{\partial  b^\ast} b^\ast-b\frac{\partial }{\partial  b}
\right)\right)
\]
and a corresponding moment map given by
\[
\Phi=\prod_{a\in \bar{Q}}(1+aa^\ast)^{\epsilon(a)}
\]
In the definition of $\Phi$ the product is taken with respect to the
chosen ordering on~$\bar{Q}$.
\end{theorem}
The algebra $A_q$ introduced in Proposition \ref{ref-1.7-3} is now 
the \emph{deformed multiplicative preprojective algebra} $\Lambda^q$
as introduced in
\cite{CBShaw}. Combining the previous theorem with Proposition 
\ref{ref-1.7-3} proves Theorem \ref{ref-1.1-0} since $\Sscr_{\alpha,q}=
\Rep(\Lambda^q,\alpha)$.

\bigskip

\subsection*{Relation with bi-symplectic geometry}
The first version of this paper was written independently of the paper
\cite{CBEG} which appeared around the same time on the ArXiv and which
discusses a non-commutative analogue of symplectic geometry. In an
appendix we outline the connection between the two papers. In
particular we prove that an algebra with a bi-symplectic form is a
double Poisson algebra. This allows us to strengthen some results of
\cite{CBEG}. For example: if $A$ has a bi-symplectic form then the
associated Lie bracket on $A/[A,A]$ is obtained from a Loday bracket
on~$A$.
\subsection*{Relation with Crawley-Boevey's Poisson structures}
In \cite{CB4} Crawley-Boevey introduces non-commutative Poisson structures
 and shows that they
induce classical Poisson structures on  moduli spaces of
representations.  In the commutative case a Poisson
structure is the same as a classical Poisson structure.

We show below that a double Poisson bracket or a double quasi-Poisson
bracket  induces a Poisson structure (Lemmas \ref{ref-2.6.2-39}
and \ref{ref-5.1.3-73}). By considering commutative algebras one easily sees
that the converse is false.

The concept of a double Poisson structure is more adapted to
algebras which are smooth in a non-commutative sense \cite{CQ}. For example a
semi-simple algebra has no Poisson structures \cite[Rem.\ 1.2]{CB4}
but it has many double Poisson structures \cite{weyer}.

\subsection*{A note on the organization of this paper}
The reader will find that this paper is rather peculiarly organized.
We have seen  that our motivating example is not a double
Poisson algebra but only a double quasi-Poisson algebra. But it seemed
difficult to treat double quasi-Poisson algebras without first
introducing the algebra of p{{oly-vector}} fields and its Schouten
bracket. This Schouten bracket is a graded version of a double Poisson
bracket. But again it seemed unreasonable to start this paper with
graded double Poisson brackets since the many signs would have
obscured the simplicity of the theory.  Jean-Louis Loday pointed out
to me that the sign problems can be mitigated by writing the
definitions in terms of functions instead of elements.  However since
the Schouten bracket has degree $-1$ some signs would still
remain\footnote{Unless one is prepared to raise the level of
  abstraction by writing the formulas in terms of operators which take
  functions as arguments!}.

So we have chosen to treat double Poisson brackets first, and then to
accept the (routine) generalizations of our statements to super
Poisson brackets without further proof or discussion.
\subsection*{Acknowledgment} This paper came out of discussions with
Crawley-Boevey and Alexei Bondal during the year on non-commutative
algebraic geometry at the Mittag-Leffler institute. It was
Crawley-Boevey who suggested that the element $\Phi$ occurring in the
definition of a multiplicative preprojective algebra could perhaps be
interpreted as a multiplicative moment map. I am very grateful for
this. The principle that one can meaningfully study non-commutative
notions through their effect on representation spaces I learned from
Lieven Le Bruyn \cite{Leb4}.

In addition I wish to thank Victor Ginzburg, Jean-Louis Loday and Geert Van
de Weyer for interesting discussions and comments. 

Jean-Louis Loday objected to my use of the term Loday algebra (taken
from  \cite{KS1}) and wished me to use the original terminology of
Leibniz algebras instead. After some consideration I decided that it is never
bad to name a concept after its inventor so I left things as they
were.

Finally I am very grateful to the anonymous referee for his thorough
reading of the manuscript and for pointing out numerous misprints and
inaccuracies.
\section{Double brackets and double Poisson algebras}
\subsection{Generalities}
Throughout we work over a field $k$ of characteristic zero although this is
not an essential condition.  Unadorned tensor products are over $k$.
If $V$, $W$ are $k$-vector spaces then an element $a\in V\otimes W$ is
written as $a'\otimes a''$. This is a short hand for $\sum_i a'_i\otimes a''_i$.
A similar convention is sometimes used for longer tensor products.
We put $a^\circ=a''\otimes a'$, i.e. $a^\circ=\sum_i a''_i\otimes a_i'$.

If $(V_i)_{i=1,\ldots,n}$ are $k$-vector spaces and $s\in S_n$ then
for $a=a_1\otimes\cdots \otimes a_n \in V_1\otimes\cdots \otimes V_n$
we put
\[
\tau_s(a)=a_{s^{-1}(1)}\otimes\cdots \otimes a_{s^{-1}(n)}
\]
so that $\tau_{st}(a)=\tau_s(\tau_t(a))$. 

Below we fix a $k$-algebra $A$. Throughout we denote the multiplication map $A^{\otimes n}\r A$ by $m$.
We will also view $A^{\otimes n}$ as an $A$-bimodule via the
\emph{outer} bimodule structure
\[
b(a_1\otimes\cdots \otimes a_n)c=ba_1\otimes\cdots \otimes a_nc
\]
Of course $A^{\otimes n}$ has many other bimodule structures. For $n=2$  we
will frequently use the \emph{inner} bimodule structure on $A^{\otimes 2}$ given
by
\[
b\ast(a_1\otimes a_2)\ast c=a_1c\otimes ba_2
\]
If $B$ is a (not necessarily commutative) $k$-algebra then a $B$-algebra will be an $k$-algebra equipped
with an (unnamed) $k$-algebra map $B\r A$.

An  element
of $\Der_B(A,A\otimes A)$ will be called a \emph{double $B$-derivation}.
\subsection{Double brackets}
\begin{definitions}
An $n$-bracket is a linear map
\[
\ldb-,\cdots,-\rdb:A^{\otimes n}\r A^{\otimes n}
\]
which is a derivation $A\r A^{\otimes n}$ in its last argument for the
outer bimodule structure on $A^{\otimes n}$  i.e.
\begin{equation}
\label{ref-2.1-4}
\ldb a_1,a_2, \ldots,a_{n-1},a_na'_n\rdb =a_n\ldb a_1,a_2, \ldots,a_{n-1},a'_n\rdb+\ldb a_1,a_2, \ldots,a_{n-1},a_n\rdb a'_n
\end{equation}
and which is cyclically anti-symmetric in the sense
\[
\tau_{(1\cdots n)}^{\phantom{-1}}\circ \ldb-,\cdots,-\rdb \circ \tau^{-1}_{(1\cdots n)}
=(-1)^{n+1}  \ldb-,\cdots,-\rdb 
\]
If $A$ is a $B$-algebra then  an $n$-bracket is \emph{$B$-linear} if it vanishes when its
last argument is in the image of $B$.
\end{definitions}
Clearly a $1$-bracket is just a derivation $A\r A$. We will call a
$2$- and a $3$-bracket respectively a \emph{double} and a
\emph{triple} bracket.  A double bracket satisfies.
\begin{align}
\label{ref-2.2-5}
\ldb a,b\rdb &=-\ldb b,a\rdb ^\circ \\
\ldb a,bc\rdb &=b\ldb a,c\rdb+\ldb a,b\rdb c
\label{ref-2.3-6}
\end{align}
The formulas \eqref{ref-2.2-5} \eqref{ref-2.3-6} imply that $\ldb -,-\rdb $ is a derivation $A\r A\otimes A$  in
its first argument
 for the inner bimodule structure on $A\otimes A$. I.e.
\begin{equation}
\label{ref-2.4-7}
\ldb ab,c\rdb =a\ast  \ldb b,c\rdb +\ldb a,c\rdb \ast b
\end{equation}
where by ``$\ast$'' we mean the inner action.
Combining \eqref{ref-2.3-6}\eqref{ref-2.4-7} we obtain 
\begin{equation}
\label{ref-2.5-8}
\ldb a_1\cdots a_m,b_1\cdots b_n\rdb = \sum_{p,q}
b_1\cdots b_{q-1} \ldb a_p,b_q\rdb 'a_{p+1}\cdots a_m
\otimes
a_1\cdots a_{p-1}\ldb a_p,b_q\rdb ''b_{q+1}\cdots b_{n}
\end{equation}
\subsection{The double Jacobi identity}
\label{ref-2.3-9}
If $a\in A$, $b=b_1\otimes \cdots \otimes b_n\in A^{\otimes n}$ then we define
\begin{align*}
\ldb a,b\rdb _L&=\ldb a,b_1\rdb \otimes b_2\otimes \cdots \otimes b_n\\
\ldb a,b\rdb_R&=b_1\otimes\cdots\otimes b_{n-1}\otimes \ldb a,b_n\rdb
\end{align*}
Associated to a double bracket $\ldb -,-\rdb $ we define a tri-ary operation
$\ldb-,-,-\rdb$ as follows:
\[
\ldb a,b,c\rdb =\ldb a,\ldb b,c\rdb \rdb _L
+\tau_{(123)}\ldb b,\ldb c,a\rdb \rdb _L
+\tau_{(132)}\ldb c,\ldb a,b\rdb \rdb _L
\]
Or in more intrinsic notations
\[
\ldb -,-,-\rdb =\ldb -,\ldb -,-\rdb \rdb _L 
+\tau_{(123)}\ldb -,\ldb -,-\rdb \rdb _L \tau^{-1}_{(123)}
+\tau^2_{(123)}\ldb -,\ldb -,-\rdb \rdb _L \tau^{-2}_{(123)}
\]
So $\ldb -,-,-\rdb $ is cyclically invariant, in the sense that
\begin{equation}
\label{ref-2.6-10}
\ldb -,-,-\rdb =\tau_{(123)}\circ \ldb -,-,-\rdb \circ \tau_{(123)}^{-1}
\end{equation}
\begin{propositions} 
\label{ref-2.3.1-11}
$\ldb -,-,-\rdb $ is a triple bracket.
\end{propositions}
\begin{proof}
  The cyclic invariance property has already been established. We now check
the derivation property.
\begin{equation}
\label{ref-2.7-12}
\begin{split}
\ldb a,\ldb b,cd\rdb \rdb _L&=\ldb a,\ldb b,c\rdb d\rdb _L+\ldb a,c\ldb b,d\rdb \rdb_L\\
&=\ldb a,\ldb b,c\rdb \rdb _Ld+\ldb a,c\rdb \ldb b,d\rdb +c\,\ldb a,\ldb b,d\rdb \rdb _L
\end{split}
\end{equation}
where in the second line we use the convention that $(x\otimes y)(s\otimes t)=x\otimes ys\otimes t$. We
will often use the same convention below.
\begin{align*}
\ldb b,\ldb cd,a\rdb \rdb_L &=\ldb b,c\ast \ldb d,a\rdb \rdb _L+\ldb b,\ldb c,a\rdb \ast d\rdb _L\\
&=\ldb b,\ldb d,a\rdb '\otimes c\ldb d,a\rdb ''\rdb _L+\ldb b,\ldb c,a\rdb 'd\otimes \ldb c,a\rdb ''\rdb _L\\
&=\ldb b,\ldb d,a\rdb '\rdb '\otimes \ldb b,\ldb d,a\rdb '\rdb ''\otimes c\ldb d,a\rdb ''
+\\&\qquad\ldb b,\ldb c,a\rdb '\rdb '\otimes \ldb b,\ldb c,a\rdb '\rdb ''d\otimes \ldb c,a\rdb ''+
\ldb c,a\rdb '\ldb b,d\rdb '\otimes \ldb b,d\rdb ''\otimes \ldb c,a\rdb ''
\end{align*}
Thus we find
\begin{equation}
\label{ref-2.8-13}
\tau_{(123)}\ldb b,\ldb cd,a\rdb \rdb _L=c\tau_{(123)}\ldb b,\ldb d,a\rdb \rdb _L+
\tau_{(123)}\ldb b,\ldb c,a\rdb \rdb _L d-\ldb a,c\rdb \ldb b,d\rdb 
\end{equation}
Finally
\begin{align*}
\ldb cd,\ldb a,b\rdb \rdb _L&=\ldb cd,\ldb a,b\rdb '\rdb \otimes \ldb a,b\rdb ''\\
&=c\ast\ldb d,\ldb a,b\rdb '\rdb \otimes \ldb a,b\rdb ''+\ldb c,\ldb a,b\rdb '\rdb \ast d\otimes \ldb a,b\rdb ''\\
&=\ldb d,\ldb a,b\rdb '\rdb '\otimes c\ldb d,\ldb a,b\rdb '\rdb ''\otimes \ldb a,b\rdb ''
+\ldb c,\ldb a,b\rdb '\rdb 'd\otimes \ldb c,\ldb a,b\rdb '\rdb ''\otimes \ldb a,b\rdb ''
\end{align*}
which yields
\begin{equation}
\label{ref-2.9-14}
\tau_{(132)}\ldb cd,\ldb a,b\rdb \rdb _L=c\tau_{(132)}\ldb d,\ldb a,b\rdb \rdb _L+
\tau_{(132)}\ldb c,\ldb a,b\rdb \rdb _Ld
\end{equation}
Taking the sum of \eqref{ref-2.7-12}\eqref{ref-2.8-13}\eqref{ref-2.9-14} yields
the desired result.
\end{proof}
\begin{definitions} A double bracket $\ldb -,-\rdb $ on $A$ is 
\emph{a double Poisson}
bracket if $\ldb -,-,-\rdb =0$.
An algebra with a double Poisson bracket is \emph{a double Poisson} algebra.
\end{definitions}
We will call the identity $\ldb -,-,-\rdb =0$ \emph{the double Jacobi identity}.
\begin{examples}
  \label{ref-2.3.3-15} Put $A=k[t]$. It is easy to check that up to
  automorphisms of $A$ the only double Poisson brackets on $A$ are
  given by
\begin{equation}
\label{ref-2.10-16}
\ldb t,t\rdb =t\otimes 1-1\otimes t
\end{equation}
and
\begin{equation}
\label{ref-2.11-17}
\ldb t,t\rdb =t^2\otimes t-t\otimes t^2
\end{equation}
These two brackets are related. Extend \eqref{ref-2.10-16} to $k[t,t^{-1}]$
(this is possible by Proposition \ref{ref-2.5.3-30} below).  It turns out that the resulting double 
bracket preserves $k[t^{-1}]$ and the corresponding restriction 
is precisely \eqref{ref-2.11-17} up to changing the sign and replacing $t$ by $t^{-1}$.

Assume that \eqref{ref-2.10-16} holds. An easy computation shows
\[
\ldb u(t),v(t)\rdb=\frac{(u(t_1)-u(t_2))(v(t_1)-v(t_2))}{t_1-t_2}
\]
where $t_1=t\otimes 1$, $t_2=1\otimes t$. From this formula it follows that any quotient
of $k[t]$ has an induced double Poisson bracket. This is for example the case
for $k[t]/(t^n)$.
\end{examples}

\subsection{Brackets associated to double brackets}
If $\ldb -,\cdots,-\rdb $ is an $n$-bracket then we put $\{-,\cdots,-\}=m\circ \ldb -,\cdots,-\rdb$.
If $n=2$ then we call $\{-,-\}$ the \emph{bracket} associated to $\ldb-,-\rdb$. By definition  $\{a,b\}=\ldb a,b\rdb '\cdot \ldb a,b\rdb ''$.   It is clear
that $\{-,-\}$ is a derivation in its second argument. I.e.
\begin{equation}
\label{ref-2.12-18}
\{a,bc\}=\{a,b\}c+b\{a,c\}
\end{equation}
and furthermore by \eqref{ref-2.2-5} 
\begin{equation}
\label{ref-2.13-19}
\{b,a\}\cong -\{a,b\}\quad \text{mod }[A,A]
\end{equation}
Finally an easy computation shows
\begin{equation}
\label{ref-2.14-20}
\{bc,a\}=\{cb,a\}
\end{equation}
\begin{lemmas}
\label{ref-2.4.1-21} 
$\{-,-\}$ induces well defined maps
\begin{equation}
\label{ref-2.15-22}
A/[A,A]\times A\r A
\end{equation}
and
\begin{equation}
\label{ref-2.16-23}
A/[A,A]\times A/[A,A]\r A/[A,A]
\end{equation}
where the latter one is anti-symmetric. 
\end{lemmas}
\begin{proof}
The map \eqref{ref-2.15-22}  is well defined by \eqref{ref-2.14-20}. From \eqref{ref-2.14-20}
together with
\eqref{ref-2.13-19} it follows that $\{a,bc\}$ is symmetric in $b,c$ modulo
commutators. 
Thus  \eqref{ref-2.16-23} is
well defined as well. Its anti-symmetry follows also from \eqref{ref-2.13-19}.
\end{proof}
\begin{propositions} If $\ldb -,-\rdb$ is a double bracket on $A$ then
the following identity holds in $A\otimes A$.
\begin{equation}
\label{ref-2.17-24}
\{a,\ldb b,c\rdb \}-\ldb \{a,b\} ,c\rdb -\ldb b,\{a,c\}\rdb
=(m\otimes 1)\ldb a,b,c\rdb -(1\otimes m) \ldb b, a,c\rdb
\end{equation}
where $m$ is the multiplication map and $\{a,-\}$ acts on tensors by $\{a,u\otimes v\}=\{a,u\}\otimes v
+u\otimes \{a,v\}$.
\end{propositions}
\begin{proof}
First we record a useful identity.
\begin{equation}
\begin{split}
\ldb a,\ldb c,b\rdb \rdb _R&=-\ldb a,\ldb b,c\rdb ''\otimes \ldb b,c\rdb '\rdb _R\\
&=-\ldb b,c\rdb ''\otimes \ldb a,\ldb b,c\rdb '\rdb \\
&=-\tau_{(123)}\biggl(\ldb a,
\ldb b,c\rdb '\rdb \otimes \ldb b,c\rdb ''\biggr)\\
&=-\tau_{(123)}\ldb a,\ldb b,c\rdb \rdb _L
\end{split}
\end{equation}
We now compute
\begin{align*}
\{a,\ldb b,c\rdb \}&=\{a,\ldb b,c,\rdb'\}\otimes \ldb b,c\rdb''
+\ldb b,c\rdb'\otimes \{a,\ldb b,c\rdb'' \}\\
&=(m\otimes 1) \ldb a,\ldb b,c\rdb \rdb_L +(1\otimes m)\ldb a,\ldb b,c \rdb \rdb_R\\
&=(m\otimes 1)\ldb a,\ldb b,c\rdb \rdb_L -(1\otimes m)\tau_{(123)} \ldb a,\ldb c,b\rdb \rdb_L
\end{align*}
\begin{align*}
\ldb \{a,b\},c\rdb&=-\ldb c,\{a,b\}\rdb^\circ\\
&=-\tau_{(12)}(\ldb c,\ldb a,b\rdb '\rdb \ldb a,b\rdb ''
+
\ldb a,b\rdb '\,\ldb c,\ldb a,b\rdb ''\rdb 
)\\
&=-(m\otimes 1)\tau_{(132)}\ldb c,\ldb a,b\rdb \rdb _L
-(1\otimes m)\tau_{(123)}\ldb c,\ldb a,b\rdb \rdb _R\\
&=-(m\otimes 1)\tau_{(132)}\ldb c,\ldb a,b\rdb \rdb _L
+(1\otimes m)\tau_{(132)}\ldb c,\ldb b,a\rdb \rdb _L
\end{align*}
\begin{align*}
\ldb b,\{a,c\}\rdb 
&=\ldb b,\ldb a,c\rdb'\rdb \ldb a,c\rdb''
+\ldb a,c\rdb' \ldb b, \ldb a,c\rdb''\rdb\\
&=(1\otimes m) \ldb b, \ldb a,c\rdb\rdb_L +(m\otimes 1) \ldb b,\ldb a,c\rdb\rdb_R\\
&=(1\otimes m) \ldb b, \ldb a,c\rdb\rdb _L-(m\otimes 1) \tau_{(123)} \ldb b,\ldb c,a\rdb \rdb_L
\end{align*}
Collecting everything we obtain the desired result.
\end{proof}
\begin{definitions} A \emph{left Loday (or Leibniz) algebra}  \cite{KS1,Loday}  is a  vector space equipped with a bilinear operation
$[-,-]$  such that the following
version of the Jacobi
identity is satisfied
\[
[a,[b,c]]=[[a,b] ,c] +[b,[a,c]]
\]
\end{definitions}
\begin{corollarys}
\label{ref-2.4.4-25}  Assume  $\ldb -,-\rdb $ is a double bracket on $A$. 
Then the following identity holds in $A$:
\begin{equation}
\label{ref-2.19-26}
\{a,\{b,c\} \}-\{ \{a,b\} ,c\} -\{ b,\{a,c\}\}
=\{ a,b,c\} -\{ b, a,c\}
\end{equation}
If $\{-,-,-\}=0$ (e.g.\ when $A$ is a double Poisson algebra)
 then $A$ becomes
a left Loday algebra.
\end{corollarys}
\begin{proof} Applying the multiplication map to  \eqref{ref-2.17-24}
  we obtain \eqref{ref-2.19-26} which in case $\{ -,-,-\}=0$ 
yields
\[
\{a,\{b,c\} \}=\{ \{a,b\} ,c\} +\{ b,\{a,c\}\}
\]
i.e.\ the defining equation for a left Loday algebra.
\end{proof}
\begin{remarks} Jean-Louis Loday asks if the relations between $\{-,-\}$ and
$\{-,-,-\}$ can be explained by some kind of Leibniz-brace-algebra structure
on $A$.
\end{remarks}
\begin{corollarys} 
\label{ref-2.4.6-27}
If $\ldb -,-\rdb $ is a a double bracket on $A$ such that
$\{-,-,-\}=0$ then
$A/[A,A]$ equipped with the bracket $\{-,-\}$ is a Lie algebra.
\end{corollarys}
\subsection{Induced brackets}
\label{ref-2.5-28}
In this section we discuss the compatibility of double brackets and 
double Poisson brackets with some natural constructions. 
\begin{propositions} Assume that $A,\ldb-,-\rdb $, $A',\ldb-,-\rdb$ are
double brackets over $B$ then there is a unique double bracket on
$A\ast_B A'$ extending the double brackets on $A$ and $A'$ with the additional
property
\[
\forall a\in A, \forall a'\in A':\ldb a,a'\rdb=0
\]
 If $A$, $A'$ are
double Poisson then so is $A\ast_B A'$.
\end{propositions}
\begin{proof} Easy.
\end{proof}
A special case is when the bracket on $A'$ is trivial. In that
case we obtain a $A'$ linear bracket on $A\ast_B A'$.  Hence we obtain:
\begin{corollarys}
\label{ref-2.5.2-29}
Double (Poisson) brackets are compatible with base change.
\end{corollarys}
Quite similarly we have
\begin{propositions} 
\label{ref-2.5.3-30}
Double (Poisson) brackets are compatible with
universal localization. I.e.\ if $S\subset A$ and $\ldb-,-\rdb$ is
a double bracket on $A$ then there is a unique extended double bracket on $A_S$. If
$A$ is double Poisson then so is $A_S$.
\end{propositions}
\begin{proof}
Left to the reader.
\end{proof}
\begin{propositions}
\label{ref-2.5.4-31}
  Assume that $e\in B$ is an idempotent. Then a $B$-linear double
  bracket $\ldb-,-\rdb$ on $A$ induces a $eBe$-linear double bracket on $eAe$.
If $\ldb-,-\rdb$ is double Poisson then so is the induced bracket on
$eAe$.
\end{propositions}
\begin{proof} This follows from
\[
\ldb eae,ebe\rdb=e\ldb a,b\rdb'e \otimes e\ldb a,b\rdb''e\qed
\]
\def\qed{}\end{proof} 
\begin{examples} In \cite{LBVW} Lieven Le Bruyn and Geert Van de
Weyer define $\sqrt[n]{A}$ as the $B$-algebra which represents the
functor of $B$-algebras to sets given by
\[
\Hom_B(A,M_n(-))
\]
A concrete realization of $\sqrt[n]{A}$ is $e(A\ast_B M_n(B))e$ where
$e$ is the upper left corner idempotent of $M_n(k)$.  So we obtain
that if $A,\ldb-,-\rdb$ is a double Poisson algebra over $B$ then
there is an induced double Poisson structure on
$\sqrt[n]{A}$.

Another realization of $\sqrt[n]{A}$ is the algebra generated by symbols
$a_{ij}$ for $a\in A$ and $i,j=1,\ldots,n$ which are linear in $a$ and
which satisfy in addition the following relations
\begin{align*}
a_{ij}b_{jk}&=(ab)_{ik}\\
b_{ij}&=\delta_{ij} b\qquad \text{if $b\in B$}
\end{align*}
where we sum over repeated indices. In this realization the double bracket
on $\sqrt[n]A$ is given by the formula
\begin{equation}
\label{ref-2.20-32}
\ldb a_{ij},b_{uv}\rdb =\ldb a,b\rdb_{uj}'\otimes \ldb a,b\rdb_{iv}''
\end{equation}
\end{examples}

Now we discuss ``fusion''. This is a
procedure which allows one to collapse two idempotents into one.  In
the case of quivers it amount to gluing vertices.  This is explained
in more detail in \S\ref{ref-6-85}.

Assume that $e_1,e_2\in B$ are 
orthogonal idempotents. Construct $\bar{A}$ from $A$ by
formally adjoining two variables $e_{12}$, $e_{21}$ satisfying the
usual matrix relations $e_{uv}e_{wt}=\delta_{wv} e_{ut}$ (with $e_{ii}=e_i$).
We have
\[
\bar{A}=A\ast_{ke_1\oplus ke_2 \oplus k\mu} (M_2(k)\oplus k\mu)=A\ast_{B}\bar{B}
\]
where $\mu=1-e_1-e_2$. 
The fusion algebra of $A$ along $e_1,e_2$ is defined as 
\[
A^f=\epsilon \bar{A}\epsilon
\]
where $\epsilon =1-e_2$. Clearly $\bar{A}$ is a $\bar{B}$-algebra and 
$A^f$ is a $B^f$-algebra. Combining Corollary \ref{ref-2.5.2-29}  and Proposition \ref{ref-2.5.4-31} we obtain:
\begin{corollarys}
\label{ref-2.5.6-33}
If $A,\ldb-,-\rdb$ is a double Poisson algebra over $B$ 
then there are associated double Poisson algebras $\bar{A}$ and 
$A^f$ over $\bar{B}$ and $B^f$ respectively.
\end{corollarys}

We recall the definition of the trace map. Let $e\in B$ be an
idempotent such that $BeB=B$. Write $1=\sum_i p_i e q_i$. Then we put
\[
\Tr:A\r eAe: a\mapsto \sum_i e q_i a p_i e
\]
The trace map depends on the chosen decomposition $1=\sum_i p_i e q_i$.
However it gives a uniquely defined isomorphism
\[
A/[A,A]\r eAe/[eAe,eAe]
\]
which is an inverse to the obvious map
\[
eAe/[eAe,eAe]\r A/[A,A]
\]
\begin{propositions}
\label{ref-2.5.7-34} We have for $a,b\in A$
\[
\Tr\{  a,b\} =\{ \Tr(a),\Tr(b)\}
\]
\end{propositions}
\begin{proof}
This is a simple computation.
\begin{align*}
\ldb \Tr(a),\Tr(b)\rdb&=\ldb \sum_i eq_i a p_ie,\sum_j eq_j b p_je\rdb\\
&=\sum_{i,j} eq_j (    eq_i* \ldb a,b\rdb       *p_ie            )  p_je\\
&=\sum_{i,j} eq_j\ldb a,b\rdb' p_ie\otimes eq_i\ldb a,b\rdb'' p_je
\end{align*}
and hence
\begin{align*}
\{ \Tr(a),\Tr(b)\}&=\sum_{i,j} eq_j\ldb a,b\rdb' p_i eq_i\ldb a,b\rdb'' p_je\\
&=\sum_{j} eq_j\{a,b\} p_je\\
&=\Tr\{a,b\}  
\end{align*}
\def\qed{}
\end{proof}

\subsection{Poisson structures and moment maps}
\label{ref-2.6-35}
By $\Der_B(A,A)$ we denote the 
$B$-derivations $A\r A$ and by $\Inn_B(A,A)\subset \Der_B(A,A)$ the subvector space of 
inner derivations (i.e.\ those derivations which are of the form $[a,-]$ for
$a$ in the centralizer of $B$). For an arbitrary linear map 
\begin{equation}
\label{ref-2.21-36}
p:A/[A,A]\r \Der_B(A,A)/\Inn_B(A,A)
\end{equation}
and for $\bar{a},\bar{b}\in A/[A,A]$ put 
\begin{equation}
\label{ref-2.22-37}
\{\bar{a},\bar{b}\}_p=
\overline{p(\bar{a})\,\tilde{}\,\,(b)}\in A/[A,A]
\end{equation}
 where 
 $p(\bar{a})\,\tilde{}$ is an arbitrary lift of $p(\bar{a})$. It is easy to show
that this is well-defined. Following Crawley-Boevey \cite[Rem.\ 1.3]{CB4}  we define
\begin{definitions} \label{ref-2.6.1-38} \cite{CB4}
If $p$ is as in \eqref{ref-2.21-36} then we say that $p$ is a \emph{Poisson structure}
on $A$ over $B$ is $\{\bar{a},\bar{b}\}_p$ is a Lie bracket on $A/[A,A]$.
\end{definitions}
We then have the following result
\begin{lemmas}
\label{ref-2.6.2-39}
If $\ldb-,-\rdb$ is a double Poisson bracket on $A$ then the map
\[
p:A/[A,A]\r \Der_B(A)/\Inn_B(A):\bar{a}\mapsto \overline{ \{a,-\}}
\]
defines a Poisson bracket of $A$ over $B$.
\end{lemmas}
\begin{proof} 
This is a combination of \eqref{ref-2.15-22} and Corollary \ref{ref-2.4.6-27}.
\end{proof}
\begin{remarks} Note that a Poisson structure is in fact a map 
\[
HH_0(A)\r HH^1(A)
\]
where ``$HH$'' denotes Hochschild (co)homology.
\end{remarks} 
The following definition will be motivated afterward. We assume that $B=ke_1\oplus \cdots
\oplus ke_n$ is semi-simple.
\begin{definitions}
\label{ref-2.6.4-40}
  Let $A,\ldb-,-\rdb$ be a double Poisson algebra. A
  \emph{moment map} for $A$ is an element $\mu=(\mu_i)_i\in \oplus_i e_i A
  e_i$ such that for all $a\in A$ we have
\[
\ldb \mu_i,a\rdb=ae_i\otimes e_i-e_i\otimes e_i a
\]
A double Poisson algebra equipped with a moment map is said to be
a \emph{Hamiltonian algebra}.
\end{definitions}
One application of a moment map is the following.
\begin{propositions}
\label{ref-2.6.5-41}
 Let $A,\ldb-,-\rdb,\mu$ be a 
Hamiltonian algebra.  Fix $\lambda\in B$ and put $\bar{A}=A/(\mu-\lambda)$.
  Then the associated Poisson structure
\[
p:A/[A,A]\r \Der_B(A,A)/\Inn_B(A,A)
\]
descends to a Poisson structure on $\bar{A}/[\bar{A},\bar{A}]$
\end{propositions}
\begin{proof}
Left to the reader.
\end{proof}
It is easy to verify that the existence of a moment map is compatible
with the induction procedures described in \S\ref{ref-2.5-28}. For further
reference we record the following.
\begin{propositions} \label{ref-2.6.6-42}(Fusion)
Assume that $A,\ldb-,-\rdb$ is a double Poisson
algebra over $B$ with moment map $\mu$. Then $A^f$ considered as a
double Poisson algebra over $B^f=
ke_1\oplus ke_3\oplus\cdots \oplus ke_n$ has a moment map given by 
\[
\mu^f_i=
\begin{cases}
\mu_1+e_{12}\mu_2 e_{21}&\text{if $i=1$}\\
\mu_i & \text{if $i\ge 3$}
\end{cases}
\]
\end{propositions}
\begin{proof} Left to the reader.
\end{proof}
\subsection{Super version}
As usual it is possible to define $\ZZ$-graded super versions of
double Poisson algebras. As usual the signs 
are determined by the Koszul convention.
We write $|a|$ for the degree of a
homogeneous element $a$ of a graded vector space.

If $V_i$, $i=1,\ldots,n$ are graded vector spaces and
$a=a_1\otimes\cdots \otimes a_n$ is a homogeneous element of
$V_1\otimes\cdots \otimes V_n$ and $s\in S_n$ then
\[
\sigma_s(a)=(-1)^{t}a_{s^{-1}(1)}\otimes\cdots \otimes a_{s^{-1}(n)}
\]
where 
\[
t=\sum_{\begin{smallmatrix} i<j \\ s^{-1}(i)>s^{-1}(j)\end{smallmatrix} }
|a_{s^{-1}(i)}||a_{s^{-1}(j)}|
\]

Let $D$ be a graded algebra.  We will call $D$ 
\emph{a double Gerstenhaber algebra}
if it is equipped  with a graded bilinear map
\[
\ldb -,-\rdb :D\otimes D\r D\otimes D
\]
of degree $-1$ such that the following identities hold:
\[
\ldb a,bc\rdb =(-1)^{(|a|-1)|b|}b\,\ldb a,c\rdb +\ldb a,b\rdb c
\]
\[
\ldb a,b\rdb =-\sigma_{(12)}(-1)^{(|a|-1)(|b|-1)}\ldb b,a\rdb 
\]
\[
0=\ldb a,\ldb b,c\rdb \rdb _L
+(-1)^{(|a|-1)(|b|+|c|)}\sigma_{(123)}\ldb b,\ldb c,a\rdb \rdb _L
+(-1)^{(|c|-1)(|a|+|b|)}\sigma_{(132)}\ldb c,\ldb a,b\rdb \rdb _L
\]
We will omit the routine verifications of graded generalizations
of the ungraded statements we have proved.
In particular if
$\ldb a,b\rdb =\ldb a,b\rdb '\ldb a,b\rdb ''$ then as in the ungraded case one proves that 
if $D,\ldb -,-\rdb $ is a double Gerstenhaber algebra then $D/[D,D][1]$ equipped
with $\{ -,-\} $ is a graded Lie algebra.
\section{P{{oly-vector}} fields and the double Schouten-Nijenhuis bracket}
\label{ref-3-43}
\subsection{Generalities}
In this section we assume that $A$ is a finitely generated $B$-algebra. 
Following \cite{CB3} we define
\[
D_{A/B}=\Hom_{A^e}(\Omega_{A/B},A\otimes A)=\Der_B(A,A\otimes A)
\]
The bimodule structure on $A\otimes A$ is the outer structure.
The surviving inner bimodule structure on $A^{\otimes 2}$ makes $D_{A/B}$ 
into an $A$-bimodule. Put $D_{B}A =T_A D_{A/B}$.  
\subsection{The double Schouten-Nijenhuis bracket}
\label{ref-3.2-44}
Our aim is to define the structure of a double 
Gerstenhaber algebra on $D_B A$.

\begin{propositions} Let $\delta,\Delta\in D_{A/B}$. Then
\begin{align*}
  \ldb \delta,\Delta\rdb \,\tilde{}_l&=(\delta\otimes 1)\Delta-(1\otimes \Delta)\delta\\
  \ldb \delta,\Delta\rdb \,\tilde{}_r&=(1\otimes
  \delta)\Delta-(\Delta\otimes 1)\delta=-\ldb \Delta,\delta\rdb\,\tilde{}_l
\end{align*}
define $B$-derivations $A\r A^{\otimes 3}$, 
where the bimodule structure on $A^{\otimes 3}$ is the outer structure.
\end{propositions}
\begin{proof} Left to the reader.
\end{proof}
Since $\Omega_{A/B}$ is finitely generated we obtain
\[
\Der_B(A,A^{\otimes 3})\cong\Hom_{A^e}(\Omega_{A/B},A\otimes A)\otimes A
\]
We will view $\ldb \delta,\Delta\rdb \,\tilde{}_l$ and $\ldb \delta,\Delta\rdb \,\tilde{}_r$ as 
elements of $D_{A/B}\otimes_k A$ and $A\otimes_k D_{A/B}$ respectively.  To this end we
define
\begin{align*}
\ldb \delta,\Delta\rdb _l&=\tau_{(23)}\circ \ldb \delta,\Delta\rdb \,\tilde{}_l\\
\ldb \delta,\Delta\rdb _r&=\tau_{(12)}\circ \ldb \delta,\Delta\rdb \,\tilde{}_r
\end{align*}
and we write
\begin{equation}
\begin{split}
\label{ref-3.1-45}
\ldb \delta,\Delta\rdb _l&=\ldb \delta,\Delta\rdb '_l\otimes \ldb \delta,\Delta\rdb ''_l\\
\ldb \delta,\Delta\rdb _r&=\ldb \delta,\Delta\rdb '_r\otimes \ldb \delta,\Delta\rdb ''_r
\end{split}
\end{equation}
with $\ldb \delta,\Delta\rdb ''_l,\ldb \delta,\Delta\rdb '_r\in A$,
$\ldb \delta,\Delta\rdb '_l,\ldb \delta,\Delta\rdb ''_r$ in $D_{A/B}$.
An easy verification shows that
\begin{equation}
\label{ref-3.2-46}
\ldb \delta,\Delta \rdb_r=-\ldb \Delta,\delta \rdb_l^\circ
\end{equation}

\medskip

For $a,b\in A$, $\delta,\Delta\in D_{A/B}$ we put
\begin{equation}
\label{ref-3.3-47}
\begin{split}
\ldb a,b\rdb &=0\\
\ldb \delta,a\rdb &=\delta(a)\\
\ldb \delta,\Delta\rdb &=\ldb \delta,\Delta\rdb _l+\ldb \delta,\Delta\rdb _r
\end{split}
\end{equation}
Here we consider the righthand sides of \eqref{ref-3.3-47} as elements
of $D_BA$.
\begin{theorems}
The definitions in \eqref{ref-3.3-47} define a unique structure of a 
double Gerstenhaber algebra on $D_B A$.
\end{theorems}

\begin{proof}
  Uniqueness is clear.  Furthermore it is easy to see that the
  derivation property and anti-symmetry of $\ldb-,-\rdb$ have to be
  checked only on generators.  Using (a graded version of) Proposition
  \ref{ref-2.3.1-11} and \eqref{ref-2.6-10} it follows that we have to
  check the double Jacobi identity only on generators also.  Thus we
  need to check the following list of identities.  For $a\in A$,
  $\alpha,\beta,\gamma\in D_{A/B}$ we need

\begin{equation}
\label{ref-3.4-48}
\ldb \alpha,\beta\rdb =-\sigma_{(12)} \ldb \beta,\alpha\rdb 
\end{equation}
\begin{equation}
\label{ref-3.5-49}
\ldb \alpha,a\beta\rdb =a\,\ldb \alpha,\beta\rdb + \ldb \alpha,a\rdb \beta
\end{equation}
\begin{equation}
\label{ref-3.6-50}
\ldb \alpha,\beta a\rdb =\ldb \alpha,\beta\rdb a+\beta \ldb \alpha,a\rdb 
\end{equation}
\begin{equation}
\label{ref-3.7-51}
0=\ldb a,\ldb \alpha,\beta\rdb \rdb _L+\sigma_{(123)}\ldb \alpha,\ldb \beta,a\rdb \rdb _L+
\sigma_{(132)}\ldb \beta,\ldb a,\alpha\rdb \rdb _L
\end{equation}
\begin{equation}
\label{ref-3.8-52}
0=\ldb \alpha,\ldb \beta,\gamma\rdb \rdb _L+\sigma_{(123)}
\ldb \beta,\ldb \gamma,\alpha\rdb \rdb _L+
\sigma_{(132)}\ldb \gamma,\ldb \alpha,\beta\rdb \rdb _L
\end{equation}
Identities \eqref{ref-3.4-48}-\eqref{ref-3.6-50} take place in $D_{A/B}\otimes A\oplus
A\otimes D_{A/B}$. \eqref{ref-3.7-51} takes place in $A^{\otimes 3}$ and 
\eqref{ref-3.8-52} takes place $D_{A/B}\otimes A\otimes A\oplus 
A\otimes D_{A/B}\otimes A\oplus A\otimes A\otimes D_{A/B}$. Taking into account
cyclic symmetry it is sufficient to prove  \eqref{ref-3.4-48} and \eqref{ref-3.8-52}
 after projection
on the first factor. So it is sufficient to prove the following identities.
\begin{equation}
\label{gensid11}\tag{\ref{ref-3.4-48}-1}
\ldb \alpha,\beta\rdb _l=-\sigma_{(12)}\ldb \beta,\alpha\rdb _r
\end{equation}
\begin{equation}
\label{gensid21}\tag{\ref{ref-3.5-49}-1}
\ldb \alpha,a\beta\rdb _l=a\,\ldb \alpha,\beta\rdb _l
\end{equation}
\begin{equation}
\label{gensid22}\tag{\ref{ref-3.5-49}-2}
\ldb \alpha,a\beta\rdb _r=a\,\ldb \alpha,\beta\rdb _r+\ldb \alpha,a\rdb \beta
\end{equation}
\begin{equation}
\label{gensid31}\tag{\ref{ref-3.6-50}-1}
\ldb \alpha,\beta a\rdb _r=\ldb \alpha,\beta\rdb _r a
\end{equation}
\begin{equation}
\label{gensid32}\tag{\ref{ref-3.6-50}-2}
\ldb \alpha,\beta a\rdb _l=\ldb \alpha,\beta\rdb _l a+\beta\,\ldb \alpha,a\rdb 
\end{equation}
\begin{equation}
\label{gensid41}\tag{\ref{ref-3.7-51}-1}
0=\ldb a,\ldb \alpha,\beta\rdb _l\rdb _L+\sigma_{(123)}\ldb \alpha,\ldb \beta,a\rdb \rdb _L+
\sigma_{(132)}\ldb \beta,\ldb a,\alpha\rdb \rdb _L
\end{equation}
\begin{equation}
\label{gensid51}\tag{\ref{ref-3.8-52}-1}
0=\ldb \alpha,\ldb \beta,\gamma\rdb _l\rdb _{l,L}+
\sigma_{(123)}\ldb \beta,\ldb \gamma,\alpha\rdb _r\rdb _L+
\sigma_{(132)}\ldb \gamma,\ldb \alpha,\beta\rdb _l\rdb _{r,L}
\end{equation}
We now check these identities systematically. By convention $\sigma$ permutes
factors in tensor products of $D_B A$ and $\tau$ permutes
factors in  tensor products of $A$ (so no signs occur in $\tau$).
\begin{itemize}
\item[\eqref{gensid11}]  This is \eqref{ref-3.2-46}.
\item[\eqref{gensid21}] We compute
\begin{align*}
\ldb \alpha,a\beta\rdb _l&=\tau_{(23)}((\alpha\otimes 1)(a\beta)-
(1\otimes a\beta)\alpha )\\
&=\tau_{(23)}((1\otimes 1\otimes a\cdot -)((\alpha\otimes 1)(\beta)-
(1\otimes \beta)\alpha ))\\
&=(1\otimes a\cdot -\otimes 1)\tau_{(23)}((\alpha\otimes 1)(\beta)-
(1\otimes \beta)\alpha )\\
&=a\,\ldb \alpha,\beta\rdb_l 
\end{align*}
\end{itemize}
\item[\eqref{gensid22}] We compute
\begin{align*}
\ldb \alpha,a\beta\rdb _r&=\tau_{(12)}((1\otimes \alpha)(a\beta)-
(a\beta\otimes 1)\alpha )\\
&=\tau_{(12)}((1\otimes a\cdot- \otimes 1)((1\otimes \alpha)(\beta)-
(\beta\otimes 1)\alpha ))+\tau_{(12)}\epsilon\\
&=( a\cdot -\otimes 1\otimes 1 )\tau_{(12)}((1\otimes \alpha)(\beta)-
(\beta\otimes 1)\alpha )+\tau_{(12)}\epsilon
\end{align*}
where $\epsilon$ is a map $A\r A^{\otimes 3}$ satisfying for $c\in A$:
\[
\epsilon(c)=\beta(c)'\otimes \alpha(a)'\otimes \alpha(a)''\beta(c)''
\]
and thus
\[
\tau_{(12)}\epsilon(c)=(\alpha(a)\beta)(c)
\]
Here $\alpha(a)$ is to be interpreted as an element of $A\otimes A\subset D_B A
\otimes D_B A$ and $\beta\in D_B A$ acts on $D_B A
\otimes D_B A$ through the outer bimodule structure.

Thus we obtain
\begin{align*}
\ldb \alpha,a\beta\rdb _r&=
( a\cdot -\otimes 1\otimes 1)\ldb \alpha,\beta\rdb _r
+\alpha(a)\beta\\
&=a\,\ldb \alpha,\beta\rdb _r+\alpha(a)\beta\\
&=a\,\ldb \alpha,\beta\rdb _r+\ldb\alpha,a\rdb \beta
\end{align*}
\item[\eqref{gensid31}\eqref{gensid32}] These are similar to 
\eqref{gensid21}\eqref{gensid22}.
\item[\eqref{gensid41}] We compute the individual terms. 
\begin{align*}
\ldb a,\ldb \alpha,\beta\rdb _l\rdb _L&=
\ldb a,\ldb \alpha,\beta\rdb '_l\rdb \otimes \ldb \alpha,\beta\rdb_l ''\\
&=-\tau_{(12)}\ldb \alpha,\beta\rdb '_l(a)\otimes \ldb \alpha,\beta\rdb ''_l\\
&=-\tau_{(12)}\tau_{(23)}((\alpha\otimes 1)\beta-(1\otimes \beta)\alpha)(a)\\
&=-\tau_{(123)}((\alpha\otimes 1)\beta-(1\otimes \beta)\alpha)(a)
\end{align*}
\[
\sigma_{(123)}\ldb \alpha,\ldb \beta,a\rdb \rdb_L=\tau_{(123)}(\alpha\otimes 1)\beta(a)
\]
\begin{align*}
\sigma_{(132)}\ldb \beta,\ldb a,\alpha\rdb \rdb_L&=-\tau_{(132)}(\beta\otimes 1)
\tau_{(12)}\alpha(a)\\
&=-\tau_{(132)}\tau_{(132)}(1\otimes \beta)\alpha(a)\\
&=-\tau_{(123)}(1\otimes \beta)\alpha(a)
\end{align*}
The sum of these three terms is indeed zero.
\item[\eqref{gensid51}] This is the most tedious computation. We compute
again the individual terms.
\begin{align*}
\ldb\alpha,\ldb\beta,\gamma\rdb _l\rdb _{l,L}&=\ldb\alpha,\ldb\beta,\gamma\rdb _l'\rdb _l\otimes 
\ldb\beta,\gamma\rdb ''\\
&=\tau_{(23)}((\alpha\otimes 1)\ldb\beta,\gamma\rdb '_l-(1\otimes \ldb\beta,\gamma\rdb '_l)
\alpha)\otimes \ldb\beta,\gamma\rdb ''_l\\
&=\tau_{(23)}((\alpha\otimes 1)\ldb\beta,\gamma\rdb _l
-(1\otimes \ldb\beta,\gamma\rdb _l)\alpha)\\
&=\tau_{(23)}(
(\alpha\otimes 1\otimes 1)\tau_{(23)}
(
(\beta\otimes 1)\gamma-(1\otimes \gamma)\beta
)
-
1\otimes  \tau_{(23)}
(
(\beta\otimes 1)\gamma
-
(1\otimes \gamma)\beta
)\alpha
)\\
&
=\tau_{(23)}\tau_{(34)}(
(\alpha\otimes 1 \otimes 1)(
(\beta \otimes 1)\gamma-(1\otimes \gamma)\beta
)
-
1\otimes (
(\beta\otimes 1)\gamma-(1\otimes \gamma)\beta
)\alpha
)\\
&=
\tau_{(234)}
(
(\alpha\otimes 1 \otimes 1)(\beta \otimes 1)\gamma
-
(\alpha\otimes 1 \otimes 1)(1\otimes \gamma)\beta
-
(1\otimes \beta \otimes 1)(1\otimes \gamma)\alpha
\\&\qquad+
(1\otimes 1\otimes \gamma)(1\otimes \beta)\alpha
)
\end{align*}
\begin{align*}
\sigma_{(132)} \ldb\beta,\ldb\gamma,\alpha\rdb _r\rdb _L
&=\tau_{(13)(24)}(\beta\ldb \gamma,\alpha\rdb _r'\otimes \ldb\gamma,\alpha\rdb _r'')\\
&=\tau_{(13)(24)}((\beta\otimes 1 \otimes 1)\ldb\gamma,\alpha\rdb _r)\\
&=\tau_{(13)(24)}
(\beta\otimes 1 \otimes 1)\tau_{(12)}(
(1\otimes \gamma)\alpha-(\alpha\otimes 1)\gamma
)\\
&=
\tau_{(13)(24)}\tau_{(132)}
(1\otimes \beta \otimes 1)
(
(1\otimes \gamma)\alpha-(\alpha\otimes 1)\gamma
)
\\
&=\tau_{(234)}
(
(1\otimes \beta \otimes 1)(1\otimes \gamma)\alpha-
(1\otimes \beta \otimes 1)(\alpha\otimes 1)\gamma
)
\end{align*}
\begin{align*}
\sigma_{(132)}\ldb\gamma,\ldb\alpha,\beta\rdb _l\rdb _{r,L}
&=
\tau_{(1432)}(
\tau_{(12)}
(
(1\otimes \gamma)\ldb\alpha,\beta\rdb _l'-(\ldb\alpha,\beta\rdb _l'\otimes 1)\gamma
)
\otimes 
\ldb\alpha,\beta\rdb _l''
)\\
&=\tau_{(1432)}\tau_{(12)}
((1\otimes \gamma\otimes 1)\ldb\alpha,\beta\rdb _l-\tau_{(34)}
(\ldb\alpha,\beta\rdb _l\otimes 1)\gamma
)\\
&=\tau_{(1432)}\tau_{(12)}
((1\otimes \gamma\otimes 1)
\tau_{(23)}((\alpha\otimes 1)\beta-(1\otimes \beta)\alpha)
\\
&\qquad -\tau_{(34)}
\tau_{(23)}
(
((\alpha\otimes 1)\beta-(1\otimes \beta)\alpha)
\otimes 1)\gamma
)\\
&=\tau_{(1432)}\tau_{(12)}\tau_{(243)}
((1\otimes 1\otimes \gamma)
((\alpha\otimes 1)\beta-(1\otimes \beta)\alpha)
\\
&\qquad -
(
((\alpha\otimes 1)\beta-(1\otimes \beta)\alpha)
\otimes 1)\gamma
)\\
&=\tau_{(234)}(
(1\otimes 1\otimes \gamma)(\alpha\otimes 1)\beta
-
(1\otimes 1\otimes \gamma)(1\otimes \beta)\alpha
-
(\alpha\otimes 1\otimes 1)(\beta\otimes 1)\gamma
\\&\qquad
+
(1\otimes \beta\otimes 1)(\alpha\otimes 1)\gamma
)
\end{align*}
And again the sum of the three terms is zero.
\end{proof}
\begin{remarks} If we equip $D_B A$ with the associated single bracket
  $\{-,-\}$ then it follows from the above theorem and Corollary
  \ref{ref-2.4.4-25} that $D_B A$ is a Loday algebra. It is easy to
  see that the Loday bracket on $D_B A$ is compatible with the map
  $D_B A\r \Der_B A$ when we equip $\Der_B A$ with the usual
  commutator bracket.
\end{remarks}
\subsection{Gauge elements}
\label{ref-3.3-53}
We now assume in addition that $B$ is commutative semi-simple. I.e.
\[
B=k e_1\oplus\cdots \oplus k e_n
\]
such that $e_i^2=e_i$. We will define some special elements $E_i$ of $D_{A/B}$ 
which we call ``gauge elements''. This terminology will be explained in 
\S\ref{ref-7.9-112}.
\[
E_i(a)=ae_i\otimes e_i-e_i\otimes e_ia
\]
We will also put $E=\sum_i E_i$. Clearly $E_i=e_i E e_i$.
\begin{propositions} For $D\in D_B A$ we have
\begin{equation}
\label{ref-3.9-54}
\ldb E_i,D\rdb =De_i\otimes e_i-e_i\otimes e_i D
\end{equation}
\end{propositions}
\begin{proof}
Since $D_B A\r D_BA \otimes D_B A:D\r De_i\otimes e_i-e_i\otimes e_i D$
is a graded derivation (of degree zero) it suffices to prove \eqref{ref-3.9-54} for
$D=a$, $a\in A$ and $D=\delta$, $\delta\in D_{A/B}$. 

We compute
\[
\ldb E_i,a\rdb =E_i(a)=ae_i\otimes e_i-e_i\otimes e_i a
\]
and
\begin{align*}
\ldb E_i,\delta\rdb _l(a)
&=\tau_{(23)}((E_i\otimes 1)\delta(a)-(1\otimes \delta)E_i(a))\\
&=\tau_{(23)}(\delta(a)'e_i\otimes e_i\otimes \delta(a)''
-
e_i\otimes e_i\delta(a)'\otimes \delta(a)''
+
e_i\otimes \delta(e_ia)'\otimes \delta(e_ia)'')\\
&=(\delta e_i)(a)\otimes e_i
\\
&=(\delta e_i\otimes e_i)(a)
\end{align*}
(where in the third line we use the $B$-linearity of $\delta$).
\begin{align*}
\ldb E_i,\delta\rdb _r(a)
&=\tau_{(12)}((1\otimes E_i)\delta(a)-(\delta\otimes 1)E_i(a))\\
&=\tau_{(12)}(\delta(a)'\otimes \delta(a)''e_i\otimes e_i
-
\delta(a)'\otimes e_i\otimes e_i\delta(a)''
-
\delta(ae_i)'\otimes \delta(ae_i)''\otimes e_i)
\\
&=-e_i\otimes (e_i\delta)(a)\\
&=-(e_i\otimes e_i\delta)(a)
\end{align*}
Taking the sum of these two expressions and letting $a$ vary we obtain 
\[
\ldb E_i,\delta\rdb =\delta e_i\otimes e_i-e_i\otimes e_i \delta\qed
\]
\def\qed{}\end{proof}
\subsection{Morita invariance}
In this section we show that $D_B A/[D_B A,D_B A]$, with its Schouten bracket,  is invariant under Morita equivalence. The fact that $D_B A$ is invariant under
Morita equivalence was already proved in \cite{CB3} but for 
the convenience of the reader we restate the proof.
 We will only consider the case when the
Morita equivalence is given by an idempotent (since this is the only 
case we will need). It is well-known that this 
implies the general case.
\begin{lemmas} 
\label{ref-3.4.1-55} Let $M$ be an $A$ bimodule and let $e\in A$ be an
idempotent such that $AeA=A$. Then 
\begin{equation}
\label{ref-3.10-56}
e(T_A M)e=T_{eAe} (eMe)
\end{equation}
and furthermore $T_A Me T_A M=T_AM$. Hence $T_{eAe} (eMe)$ is Morita equivalent
to $T_A M$.
\end{lemmas}
\begin{proof}
Since $AeA=A$ we have $Ae\otimes_{eAe} eA\cong A$. Thus for $A$-bimodules
$M$, $N$ we obtain 
\begin{align*}
e(M\otimes_A N)e&= eA\otimes_A M \otimes_A Ae \otimes_{eAe}  eA\otimes_A M\otimes_A Ae\\
&=eMe\otimes_{eAe} eNe 
\end{align*}
This implies \eqref{ref-3.10-56}. The second assertion is clear since
$T_A Me T_A M$ contains $A=AeA$.
\end{proof}
\begin{lemmas} Assume $e\in B$ is an idempotent such that $BeB=B$. Then 
\[
eD_{A/B}e=D_{eAe/eBe}
\]
\end{lemmas}
\begin{proof}
There is an obvious map
\[
c:eD_{A/B}e\r D_{eAe/eBe}
\]
We have to construct its inverse $c^{-1}$. Write
\[
1=\sum_i p_ieq_i
\]
with $p_i,q_i\in B$. Then any element $a\in A$ can be written as
\[
\sum_{i,j,k} p_i(eq_i a p_j e)  q_j
\]
Let $\delta\in D_{eAe/eBe}$. We put
\[
c^{-1}(\delta)(a)=\sum_{i,j,k} 
p_i\delta(eq_i a p_j e) q_j
\]
It is easy to see that this is a well-defined element of $eD_{A/B}e$ and
that $c^{-1}$ is indeed a two sided inverse to $c$.
\end{proof}
Using lemma \ref{ref-3.4.1-55} we obtain that  there is an isomorphism
\[
D_{eBe}(eAe)=eD_B A e
\]
\begin{propositions} 
\label{ref-3.4.3-57} Assume $e\in B$ is an idempotent such that $BeB=B$.
Then the Schouten bracket on $D_BA$ restricted to $e(D_BA) e=D_{eBe}(eAe)$
coincides with the Schouten bracket on $D_{eBe}(eAe)$.
\end{propositions} 
\begin{proof}
Since $e(D_BA) e=T_{eAe}(eD_{A/B}e)$ it suffices to check that the
Schouten bracket on $D_{eAe/eBe}$ and the restricted Schouten bracket
on $eD_{A/B}e$ coincide. Since $\delta\in eD_{A/B}e$ restricts to
a derivation $eAe\r eAe\otimes eAe$ it is easy to see that both Schouten
brackets are given by the same formulas.
\end{proof}

\subsection{Hamiltonian vector fields}
\label{ref-3.5-58}
Assume that $A$ is a equipped with a $B$-linear double bracket.
If $a\in A$ then we write $H_a=\ldb a,-\rdb$. We call $H_a$ the
\emph{Hamiltonian vector field corresponding to $a$}. Using this notation
we may write
\begin{equation}
\label{ref-3.11-59}
\ldb a,b\rdb=H_a(b)
\end{equation}
\begin{propositions} 
  \label{ref-3.5.1-60} The following are equivalent
\begin{enumerate}
\item $\ldb-,-\rdb$ is a double Poisson bracket.
\item The following identity holds for all $a,b\in A$:
\[
\ldb H_a, H_b\rdb_l= H_{\ldb a,b \rdb'} \otimes \ldb a,b \rdb''
\]
\item The following identity holds for all $a,b\in A$:
\[
\ldb H_a, H_b\rdb_r= \ldb a,b \rdb' \otimes H_{{\ldb a,b \rdb'}''}
\]
\item The following identity holds for all $a,b\in A$:
\[
\ldb H_a, H_b\rdb= H_{{\ldb a,b \rdb}}
\]
where we use the convention $H_{x'\otimes x''}=H_{x'}\otimes x''+x'\otimes H_{x''}$.
\end{enumerate}
\end{propositions}
\begin{proof} We first prove the equivalence of (1) and (2). We have
  to rewrite the expression for the associated triple bracket
\begin{equation}
\label{ref-3.12-61}
\ldb a,b,c\rdb=\ldb a, \ldb b,c \rdb \rdb_L+
\tau_{123}\ldb b, \ldb c,a \rdb \rdb_L+
\tau_{132}\ldb c, \ldb a,b \rdb \rdb_L
\end{equation}
For the first term we use
\[
\ldb a, \ldb b,c \rdb \rdb_L=(H_a\otimes 1) H_b(c)
\]
For the second term we use
\[
\ldb b, \ldb c,a \rdb \rdb_L=-(H_b\otimes 1)((H_a (c))^\circ)
\]
and hence
\begin{align*}
\tau_{123}\ldb b, \ldb c,a \rdb \rdb_L &=-\tau_{123}(H_b\otimes 1)\tau_{12}(H_a (c))\\
&=-\tau_{123}\tau_{132}(1\otimes H_b)H_a (c)\\
&=-(1\otimes H_b)H_a (c)
\end{align*}
For the third term we use
\begin{align*}
 \ldb h,x'\otimes x''\rdb_L&=\ldb h,x'\rdb \otimes x''\\
&=-\ldb x',h\rdb^\circ \otimes x''\\
&=-(H_{x'}h)^\circ \otimes x''
\end{align*} 
and thus
\[
\ldb c, \ldb a,b \rdb \rdb_L=-(H_{\ldb a,b\rdb'} c)^\circ \otimes \ldb a,b\rdb''
\]
and hence
\begin{align*}
\tau_{132}\ldb c, \ldb a,b \rdb \rdb_L&=-\tau_{132}\tau_{12}(H_{\ldb a,b\rdb'} c)
\otimes \ldb a,b\rdb''\\
&=-\tau_{23} (H_{\ldb a,b\rdb'} c
\otimes \ldb a,b\rdb'')
\end{align*}
So we get that  $\ldb a, b, c\rdb$ is equal to 
\begin{align*}
(H_a\otimes 1) H_b(c)-(1\otimes H_b)H_a (c)-\tau_{23} (H_{\ldb a,b\rdb'} c
\otimes \ldb a,b\rdb'')
&=
\ldb H_a,H_b\rdb\tilde{}_l-\tau_{23} (H_{\ldb a,b\rdb'} c
\otimes \ldb a,b\rdb'')\\
&=\tau_{23}(\ldb H_a,H_b\rdb_l-H_{\ldb a,b\rdb'}(-)
\otimes \ldb a,b\rdb'')(c)
\end{align*}
which completes the proof of the equivalence of (1) and (2).

To prove the implication (2)$\Rightarrow$(3) we use \eqref{ref-3.2-46}.
Assuming (2) we obtain
\begin{align*}
\ldb H_a,H_b\rdb_r&=-\ldb H_b,H_a\rdb_l^\circ\\ 
&=-\ldb b,a\rdb'' \otimes H_{\ldb b,a\rdb'}\\
&=\ldb a,b\rdb' \otimes H_{\ldb a,b\rdb''}
\end{align*}
The implication (3)$\Rightarrow$(2) is proved in the same way. 
(4) is the sum of (2) and (3).  To go back we regard (4) as an
identity in $D_B A\otimes A \oplus A\otimes D_B A$. The projection on the
two terms yields (2) and (3).
\end{proof}

\section{The relation between p{{oly-vector}} fields and 
brackets}
\subsection{Generalities}
\label{ref-4.1-62}
We assume that $A$ is a finitely generated $B$-algebra.
\begin{propositions} 
There is a well defined linear map
\begin{equation}
\label{ref-4.1-63}
\mu:(D_B A)_n \r \{\text{$B$-linear $n$-brackets on $A$}\}:Q\mapsto \ldb-,\cdots,-\rdb_Q
\end{equation}
which on $Q=\delta_1\cdots\delta_n$ is given by 
\[
\ldb -,\cdots,-\rdb_Q=\sum_{i=0}^{n-1} (-1)^{(n-1)i}
\tau_{(1\cdots n)}^i \circ \ldb -,\cdots,-\rdb\,\tilde{}_Q
\circ\tau_{(1\cdots n)}^{-i}
\]
where
\begin{equation}
\label{ref-4.2-64}
\ldb a_1,\cdots,a_n\rdb\,\tilde{}_Q=\delta_n(a_n)'\delta_1(a_1)''
\otimes\delta_1(a_1)' \delta_2(a_2)''\otimes \cdots
\otimes
\delta_{n-1}(a_{n-1})'\delta_n(a_n)'' 
\end{equation}
This map factors through $D_BA/[D_BA,D_BA]$.
\end{propositions}
\begin{proof} This follows easily from the following alternative  formula
\[
\ldb a_1,\cdots,a_n\rdb_{\delta_1\cdots\delta_n}=\sum_{i=0}^{n-1} (-1)^{(n-1)i}
\ldb a_1,\cdots,a_n\rdb\,\tilde{}_{\delta_{n-i+1}\cdots \delta_n\cdots \delta_1\cdots\delta_{n-i}}\qed
\]
\def\qed{}\end{proof} Slightly generalizing \cite{CQ} let us say that
$A/B$ is \emph{smooth} if $A/B$ is left and right flat and
$\Omega_{A/B}$ is a projective $A$-bimodule (in addition to $A/B$
being finitely generated).
\begin{propositions} \label{ref-4.1.2-65} If $A/B$ is smooth
then $\mu$ is an isomorphism.
\end{propositions}
\begin{proof} To prove this it will be convenient to work in slightly greater
generality. 

Let  $M$ be an $A$-bimodule. We put
$M^\ast=\Hom_{A^e}(M,A\otimes A)$ where we use the outer bimodule structure on
$A^{\otimes 2}$. We view $M^\ast$ as an $A$-bimodule through the inner bimodule
structure on $A^{\otimes 2}$.

We will consider $M^{\otimes n}$ as a $(A^e)^{\otimes n}$-modules where the 
the $i$'th copy of $A^e$ acts on the $i$'th copy of $M$. We will also consider
an $(A^e)^{\otimes n}$-module $[A^{\otimes (n+1)}]$ which
is equal to $A^{\otimes n+1}$ as vector space and
 where the $i$'th copy of $A^e$ act as follows
\[
(a'\otimes a'')(a_1\otimes \cdots \otimes a_{n+1})=a_1\otimes 
\cdots \otimes a_ia''\otimes a'a_{i+1}  \otimes \cdots \otimes a_{n+1}
\]
All these bimodule structures commute with the outer bimodule
structure on $[A^{\otimes (n+1)}]$.
 There is a morphism of $A$-bimodules 
\[
\Psi:(M^\ast)^{\otimes_A n}\r \Hom_{(A^e)^{\otimes n}}(M^{\otimes n},[A^{\otimes (n+1)}])
\]
given by
\[
\phi_1\otimes \cdots \otimes \phi_n\mapsto
\bigl(m_1\otimes \cdots \otimes m_n\mapsto 
\phi_1(m_1)''\otimes \phi_1(m_1)'\phi_2(m_2)''
\otimes\cdots \otimes \phi_{n-1}(m_{n-1})' \phi_n(m_n)''\otimes \phi_n(m_n)'\bigr)
\]
In case $M$ is a finitely generated projective bimodule then this is an isomorphism.
To prove this one may assume $M=A\otimes_k A$, in which case it is easy.

Let $\{A^{\otimes n}\}$ be the $(A^{e})^{\otimes n }$-module which
is $A^{\otimes n}$ as a vector space and where the $i$'th copy
of $A^e$ for $i=1,\ldots,n-1$ acts as on $[A^{\otimes n+1}]$ but
where the $n$'th copy acts by the outer bimodule structure.  

If $N$ is an $A$-bimodule then we have that $N\otimes_{A^e}A=N/[A,N]$. Therefore
we denote an element of $N\otimes_{A^e}A$ by $\bar{n}$ where $n\in N$. 

The map
\[
\overline{a_1\otimes\cdots \otimes a_{n+1}}\r a_{n+1}a_1\otimes\cdots \otimes a_{n}
\]
gives an isomorphism
\[
[A^{\otimes (n+1)}]\otimes_{A^e} A\cong \{A^{\otimes n}\}
\]
We define $\psi$ as the composition
\begin{multline*}
(M^\ast)^{\otimes_A n}\otimes_{A^e} A
\xrightarrow{\Psi\otimes 1} 
\Hom_{(A^e)^{\otimes n}}(M^{\otimes n},[A^{\otimes (n+1)}])\otimes_{A^e} A
\\\xrightarrow{\text{can.}} \Hom_{(A^e)^{\otimes n}}(M^{\otimes n},[A^{\otimes (n+1)}]\otimes_{A^e} A)
\cong \Hom_{(A^e)^{\otimes n}}(M^{\otimes n},\{A^{\otimes n}\})
\end{multline*}
Explicitly:
\[
\overline{\phi_1\otimes \cdots \otimes \phi_n}\mapsto
\bigl(m_1\otimes \cdots \otimes m_n\mapsto 
\phi_n(m_n)'\phi_1(m_1)''\otimes \phi_1(m_1)'\phi_2(m_2)''
\otimes\cdots \otimes \phi_{n-1}(m_{n-1})' \phi_n(m_n)'' \bigr)
\]
It is clear that $\psi$ will also be an isomorphism if $M$ is finitely
generated projective.
The cyclic group $C_n$ acts on $(M^\ast)^{\otimes_A n}\otimes_{A^e} A$ by
\[
\sigma_{(1\cdots n)}(\overline{\phi_1\otimes\cdots \otimes \phi_n})=(-1)^{n-1}\overline{\phi_n\otimes \phi_1\otimes
\cdots \phi_{n-1}}
\]
An easy verification shows that the following diagram is commutative. 
\[
\begin{CD}
(M^\ast)^{\otimes_A n}\otimes_{A^e}A @>\psi>> \Hom_{(A^e)^{\otimes n}}(M^{\otimes n},\{A^{\otimes n}\})\\
@V \sigma_{(1\cdots n)} VV @VV(-1)^{n-1}\tau_{(1\cdots n)}\circ -\circ \tau^{-1}_{(1\cdots n)}  V\\
(M^\ast)^{\otimes_A n}\otimes_{A^e}A @>\psi>> \Hom_{(A^e)^{\otimes n}}(M^{\otimes n},\{A^{\otimes n}\})\end{CD}
\]
Let $\inv$ and $\coinv$ denote respectively the signed invariants and coinvariants
for the action of $C_n$. We view $T_A M^\ast$ as a graded ring with $M^\ast$ in degree
$1$. In that case we have
\[
\left(T_A M^\ast/[T_A M^\ast, T_A M^\ast]\right)_n=\coinv \left((M^\ast)^{\otimes_A n}
\otimes_{A^e} A\right)
\]
where $[-,-]$ means signed commutators. 

We define $\mu$ as the composition of the maps.
\begin{multline*}
\left(T_A M^\ast/[T_A M^\ast, T_A M^\ast]\right)_n=\coinv \left((M^\ast)^{\otimes_A n}
\otimes_{A^e} A\right)\xrightarrow{\psi} \coinv \Hom_{(A^e)^{\otimes n}}(M^{\otimes n},\{A^{\otimes n}\})\\\xrightarrow[\cong]{\text{trace}} \inv  \Hom_{(A^e)^{\otimes n}}(M^{\otimes n},\{A^{\otimes n}\})
\end{multline*}
If $M$ is finitely generated projective then $\mu$ is an isomorphism.
$\overline{\phi_1\otimes_A \cdots \otimes_A \phi_n}$ is mapped under $\mu$ to 
\begin{equation}
\label{ref-4.3-66}
\sum_i (-1)^{(n-1)i}
\tau_{(1\cdots n)}^i \circ \Phi
\circ\tau_{(1\cdots n)}^{-i}
\end{equation}
where 
\[
\Phi(m_1\otimes\cdots \otimes m_n)= \phi_n(m_n)'\phi_1(m_1)''\otimes \phi_1(m_1)'\phi_2(m_2)''
\otimes\cdots \otimes \phi_{n-1}(m_{n-1})' \phi_n(m_n)''
\]
Now consider the case $M=\Omega_{A/B}$. In that case there is an
isomorphism 
\[
\inv \Hom_{(A^e)^{\otimes n}}(\Omega^{\otimes n},\{A^{\otimes n}\})
\cong \{\text{$B$-linear $n$-brackets on $A$}\}:
\]
which maps $\theta\in \Hom_{(A^e)^{\otimes n}}(\Omega^{\otimes n},\{A^{\otimes n}\})$ to
the bracket
\[
\ldb a_1,\ldots,a_n\rdb=\theta(da_1\otimes\cdots\otimes da_n)
\]
Composing this identification with the map $\mu$ defined by \eqref{ref-4.3-66} gives us
precisely \eqref{ref-4.1-63}.
\end{proof}
\subsection{Compatibility}
\begin{propositions}
\label{ref-4.2.1-67}
For $Q\in (D_B A)_n$   the following identity  holds
\[
\ldb a_1,\ldots,a_{n}\rdb_{Q}=(-1)^{\frac{n(n-1)}{2}}\ldb a_1,\ldots,\ldb a_{n-1}, \{Q,a_n\}
\rdb_L\cdots \rdb_L
\]
\end{propositions}
\begin{proof}
It suffices to prove this for $Q=\delta_1\cdots \delta_n$ with $\delta_i\in D_BA$. We compute
\begin{multline*}
\{\delta_1\cdots \delta_n,a_n\}=(-1)^{n-1}\delta_1(a_n)'\delta_2\cdots \delta_n \delta_1(a_n)''
+(-1)^{n-2}\delta_2(a_n)'\delta_3\cdots \delta_n\delta_1 \delta_2(a_n)''+\cdots\\
+
\delta_n(a_n)'\delta_1\cdots \delta_{n-1} \delta_n(a_n)''
\end{multline*}
We concentrate on the last term. The other terms are obtained by
cyclically permuting the $\delta$'s.
We find
\begin{multline*}
\ldb a_1,\ldots,\ldb a_{n-1}, \delta_n(a_n)'\delta_1\cdots \delta_{n-1} \delta_n(a_n)''\rdb_L\cdots \rdb_L\\
=
(-1)^{\frac{n(n-1)}{2}}
\delta_n(a_n)'
\delta_1(a_1)''\otimes \delta_1(a_1)'
\delta_2(a_2)''\otimes \delta_2(a_2)'
\cdots
\delta_{n-1}(a_{n-1})''\otimes \delta_{n-1}(a_{n-1})'
\delta_n(a_n)''\\
  =(-1)^{\frac{n(n-1)}{2}}\ldb
  a_1,\ldots,a_{n}\rdb\,\tilde{}_{\delta_1\cdots\delta_{n}\delta_n}
\end{multline*}
We find
\[
  \ldb a_1,\ldots,\ldb a_{n-1}, \{\delta_1\cdots \delta_n,a_n\}
  \rdb_L\cdots \rdb_L
  =
  (-1)^{\frac{n(n-1)}{2}}
  \ldb a_1,\ldots,a_{n}\rdb_{\delta_1\cdots\delta_n}
\]
which is what we want.
\end{proof}

\begin{propositions}
Let $P\in ( D_B A)_2$.
We have the following identity for $a,b,c\in A$:
\[
-(1/2)\ldb a,\ldb b,\{\{ P, P\},c\}\rdb\rdb_L=
\ldb a,\ldb b,c \rdb_P \rdb_{P,L}+\tau_{(123)}\ldb b,\ldb c,a\rdb_P\rdb_{P,L}
+\tau_{(132)}\ldb c,\ldb a,b\rdb_P\rdb_{P,L}
\]
\end{propositions}
\begin{proof}
By (the graded version of) \eqref{ref-2.19-26} we have
\[
\{\{P,P\},c\}=2\{P,\{P,c\}\}
\]
and by \eqref{ref-2.17-24}
\[
\ldb b, \{P,\{P,c\}\}\rdb=-\{P,\ldb b,\{P,c\}\rdb\}+\ldb \{P,b\},\{P,c\}\rdb
\]
We now apply $\ldb a,-\rdb_L$ to the individual terms.
\begin{align*}
\ldb a,\{P,\ldb b,\{P,c\}\rdb\}\rdb_L&=
\ldb a,\{P, \ldb b,\{ P,c\}\rdb'\}\rdb\otimes \ldb b,\{P,c\}\rdb''\\
&=\ldb \{P,a\}, \ldb b,\{P,c\}\rdb'\rdb\otimes \ldb b,\{P,c\}\rdb''\\
&=\ldb a,\ldb b,c\rdb_P\rdb_{P,L}
\end{align*}
where in the second line we have used the graded version of 
\eqref{ref-2.17-24} in the form of the formula
\begin{equation}
\label{ref-4.4-68}
\ldb \{P,a\},x\}=\ldb a,\{P,x\}\rdb
\end{equation}
 In the third line we have once more used this identity.

 Recall for the double bracket on $D_B A$ is odd. This explains the signs.

A similar computation yields
\begin{align*}
\ldb a, \ldb \{P,b\},\{P,c\}\rdb\rdb_L
&=
-\sigma_{(123)}\ldb \{P,b\},\ldb \{P,c\},a\rdb\rdb_L
-\sigma_{(132)}\ldb \{P,c\},\ldb a,\{P,b\}\rdb\rdb_L
\\
&=-\tau_{(123)}\ldb b,\ldb c,a\rdb_P\rdb_{P,L}-\tau_{(132)}\ldb c,\ldb a,b
\rdb_P\rdb_{P,L}
\end{align*}
Collecting everything proves the proposition.
\end{proof}
Summarizing we obtain
\begin{theorems}
\label{ref-4.2.3-69}  Let $P\in (D_BA)_2$. Then one has
\[
\ldb a,b,c \rdb_{1/2 \{P,P\}}=
\ldb a,\ldb b,c \rdb_P \rdb_{P,L}+\tau_{(123)}\ldb b,\ldb c,a\rdb_P\rdb_{P,L}
+\tau_{(132)}\ldb c,\ldb a,b\rdb_P\rdb_{P,L}
\]
\end{theorems}
\subsection{The trace map}
\label{ref-4.3-70}
 As above let $e\in B$ be an idempotent
such that $BeB=B$. We have an associated trace map
\[
\Tr:D_BA/[D_BA,D_BA]\r eD_BAe/[eD_BAe,eD_BAe]=
D_{eBe}(eAe)/[D_{eBe}(eAe),D_{eBe}(eAe)] 
\]
respecting $\{-,-\}$ by Propositions \ref{ref-2.5.7-34} and
 \ref{ref-3.4.3-57}.

Furthermore one has
\[
\ldb-,\cdots,-\rdb_Q =\ldb-,\cdots,-\rdb_{\Tr(Q)}
\]
since $\Tr(Q)=\sum_i eq_i Q p_ie$ (with the notations of \S\ref{ref-2.5-28}) and 
hence
\begin{align*}
\ldb-,\cdots,-\rdb_{\Tr(Q)}&=\ldb-,\cdots,-\rdb_{\sum_i eq_i Q p_ie}\\
&=\ldb-,\cdots,-\rdb_{\sum_i  p_ieq_i Q}\\
&=\ldb-,\cdots,-\rdb_{Q}
\end{align*}
\subsection{Differential  double Poisson algebras}
\begin{definitions} 
We say that $A$ is a \emph{differential double Poisson} algebra (DDP) over $B$ if it is equipped 
an element $P\in (D_BA)_2$ (a \emph{differential double Poisson bracket}) such that
\begin{equation}
\label{ref-4.5-71}
\{P,P\}=0 \quad\operatorname{mod}\ [D_BA,D_BA]
\end{equation}
\end{definitions}
If $A,P$ is a differential double Poisson algebra then by Theorem
\ref{ref-4.2.3-69} $A$ is a double Poisson algebra with double
bracket $\ldb-,-\rdb_P$. From Proposition \ref{ref-4.1.2-65} it follows that 
in the smooth case the notions of differential double Poisson algebra and double
Poisson algebra are equivalent. This is not true in the non-smooth case
as the following example shows.
\begin{examples} 
Let $A=k[\epsilon]/(\epsilon^2)$, $B=k$. 
According to Example \ref{ref-2.3.3-15} $A$
has a double Poisson bracket given by $\ldb \epsilon, \epsilon\rdb=\epsilon\otimes 1
-1\otimes \epsilon$. 

On the other hand it easy to check that every element of
$\Der(A,A\otimes A)$ sends $\epsilon$ to $k\epsilon \otimes \epsilon$.
Using \eqref{ref-4.2-64} we deduce that if $P\in (D_B A)_2$ then $\ldb
\epsilon, \epsilon\rdb_P=0$. So $\ldb-,-\rdb$ is not differential.
\end{examples}
\begin{examples} Assume that $A=k[t]$ and $B=k$. Then the double
  Poisson bracket $\ldb t,t\rdb=t\otimes 1 - 1 \otimes t$ is obtained
  from
\[
P=t\frac{\partial\ }{\partial t} \frac{\partial\ }{\partial t}
\]
where $\partial/\partial t$ is the double derivation defined by
$\partial t/\partial t=1\otimes 1$.
\end{examples}
\begin{propositions} 
  \label{ref-4.4.4-72} If $P\in (D_BA)_2$ is a differential double
  Poisson bracket then $\mu\in \oplus_i e_i A e_i$ is a moment map
  (cfr.\ Definition \ref{ref-2.6.4-40}) for $\ldb-,-\rdb_P$ if and only if
\[
\{P,\mu_i\}=-E_i
\]
\end{propositions}
\begin{proof}
By Proposition \ref{ref-4.2.1-67} and \eqref{ref-4.4-68} we have
\[
\ldb \mu_i ,a\rdb_P=-\ldb \{P,\mu_i\}, a\rdb=-\{P,\mu_i\}(a) 
\]
Thus $\mu$ is indeed a moment map if and only if $\{P,\mu_i\}=-E_i$.
\end{proof}
It seems logical to call a differential double Poisson algebra equipped
with a moment map a differential Hamiltonian algebra.
\section{Double quasi-Poisson algebras}
We now we introduce a twisted version of double Poisson algebras.
For simplicity we assume throughout 
that $B=ke_1\oplus \cdots \oplus ke_n$ is semi-simple.
\subsection{General definitions}
\begin{definitions}
A \emph{double quasi-Poisson bracket} on $A$ (over $B$) is a $B$-linear bracket $\ldb-,-\rdb
$ such that 
\[
\ldb-,-,-\rdb=\frac{1}{12}\sum_i \ldb-,-,-\rdb_{E_i^3}
\]
We say that $A$ is a \emph{double quasi-Poisson algebra} over $A$ if $A$ is equipped
with a double quasi-Poisson bracket.
\end{definitions}
\begin{propositions} If $A,\ldb-,-\rdb$  is a double quasi-Poisson algebra
then $A,\{-,-\}$ is a left Loday algebra. 
\end{propositions} 
\begin{proof}
According to Corollary \ref{ref-2.4.4-25} we have to show  
\[
\{-,-,-\}_{E_i^3}=0
\]
This identity is immediate from the definition.
\end{proof}
In a similar way we obtain
\begin{lemmas} \label{ref-5.1.3-73} If $A,\ldb-,-\rdb$ is a double quasi-Poisson 
algebra then $\{-,-\}$ induces a Poisson structure on $A$.
\end{lemmas}
\begin{proof} This is proved as Lemma \ref{ref-2.6.2-39}.
\end{proof}
\begin{definitions}
  Let $A,\ldb-,-\rdb$ be a double quasi-Poisson algebra. A
  \emph{multiplicative moment map} for $A$ is an element
  $\Phi=(\Phi_i)_i\in \oplus_i e_i A e_i$ such that for all $a\in A$
  we have
\[
\ldb \Phi_i,a\rdb=\frac{1}{2} (\Phi_i E_i+E_i \Phi_i)(a)
\]
A double quasi-Poisson algebra equipped with a moment map is said to be
a \emph{quasi-Hamiltonian algebra}.
\end{definitions}
\begin{propositions} 
\label{ref-5.1.5-74} Let $A,\ldb-,-\rdb,\Phi$ be a 
quasi-Hamiltonian algebra.  Fix $q\in B^\ast$
  and put $\bar{A}=A/(\Phi-q)$.  Then the associated Poisson structure
\[
p:A/[A,A]\r \Der_B(A,A)/\Inn_B(A,A)
\]
descends to a Poisson structure on $\bar{A}/[\bar{A},\bar{A}]$
\end{propositions}
\begin{proof}
Left to the reader.
\end{proof}
\subsection{Differential versions}
\label{ref-5.2-75}
\begin{definitions}
  We say that $A$  is a \emph{differential double quasi-Poisson} algebra (DDQP-algebra) 
 over  $B$-algebra $A$ if $A$ is equipped with an element $P\in (D_BA)_2$ 
(a \emph{differential double quasi-Poisson bracket}) such that
\begin{equation}
\label{ref-5.1-76}
\{P,P\}=\frac{1}{6} \sum_{i=1}^n E_i^3 \quad \operatorname{mod}\ [D_BA,D_BA]
\end{equation}
\end{definitions}
It follows from Theorem \ref{ref-4.2.3-69} that a DDQP-algebra is a 
double quasi-Poisson algebra. For smooth algebras the two notions are equivalent by
Proposition \ref{ref-4.1.2-65}.
\begin{propositions} If $P\in (D_BA)_2$ is a double quasi-Poisson bracket then $\Phi\in \oplus_i e_i A
  e_i$ is a multiplicative moment map for $\ldb-,-\rdb_P$ if and only if
\[
\{P,\Phi_i\}=-\frac{1}{2}(E_i\Phi_i+\Phi_iE_i)
\]
in $D_{A/B}$.
\end{propositions}
\begin{proof} This is similar to the proof of \ref{ref-4.4.4-72}.
\end{proof}
A differential quasi-Hamiltonian algebra is a quasi-Hamiltonian algebra where
the double bracket coming from an element of $(D_B A)_2$.
\subsection{Calculus on fusion algebras}
In this section the notations are as in \S\ref{ref-2.5-28}.
Our aim is to show that if $A$ is double quasi-Poisson algebra or
a quasi-Hamiltonian algebra over $B$ then the same is true for the fused
algebra $A^f$. Why this is to be expected will be explained in 
\S\ref{ref-7.10-114}. The methods in this section are basically
translations of the methods in \cite[\S5]{AKM}.

The non-quasi-versions of these methods
are  easy and have been treated in Corollary \ref{ref-2.5.6-33}
and Proposition \ref{ref-2.6.6-42}. The quasi-case is more tricky notation
wise. For this reason we will restrict ourselves to the differential
case. I have no doubt that the general case also works but I have not checked it. 

Extending derivations yields a canonical map
\[
D_{A/B} \r D_{\bar{A}/\bar{B}}
\]
and hence a corresponding map
\[
\bar{(-)}:D_B A\r D_{\bar{B}} \bar{A}
\]
We will often identify $D_B A$ with its image in $D_{\bar{B}} \bar{A}$. It
is easy to see that $\bar{(-)}$ is compatible with the Schouten bracket.

By composition we define a map
\[
(-)^f:D_BA \r 
D_{B^f}(A^f): P\mapsto \Tr(\bar{P})
\]
where we compute $\Tr$ using the decomposition $1=1\cdot \epsilon \cdot 1+e_{21}
\epsilon e_{12}$.
It follows from \S\ref{ref-4.3-70} that $(-)^f$ is compatible
with Schouten brackets. 

For convenience we now define some operators in $D_{\bar{A}/\bar{B}}$.
In order to avoid confusing notations we define $F_i\in D_{\epsilon
  B\epsilon} (\epsilon \bar{A} \epsilon)$ for $i\neq 2$ by
$F_i(a)=ae_i\otimes e_i-e_i\otimes e_i a$ Note that $E_i^f=F_i$ for
$i>2$ but this is not case for $i=1$.

In this section we prove the following two results.
\begin{theorems}  \label{ref-5.3.1-77}
Assume that $A,P$ is a differential double
  quasi-Poisson algebra over $B$. Then $A^f,P^{f\!\!f}$ with
\begin{equation}
\label{ref-5.2-78}
P^{f\!\!f}=P^f-\frac{1}{2}E_1^f E_2^f
\end{equation}
is a differential double
  quasi-Poisson algebra.
\end{theorems}
\begin{theorems}
\label{ref-5.3.2-79} Assume that $A,P,\Phi$ is a differential quasi-Hamiltonian
  algebra over $B$. Then $A^f,P^{f\!\!f},\Phi^{f\!\!f}$ with $P^{f\!\!f}$ as in
\eqref{ref-5.2-78} and
 with
\[
\Phi^{f\!\!f}_i=
\begin{cases}
\Phi^f_1\Phi^f_2&\text{if $i=1$}\\
\Phi^f_i&\text{if $i>2$}
\end{cases}
\]
is a differential quasi-Hamiltonian algebra over $B^f$.
\end{theorems}
The proof of these theorems needs some preparation. We  put
\begin{align*}
E&=\bar{E}_1\\
\hat{E}&=e_{12} \bar{E}_2e_{21}
\end{align*}
\begin{lemmas} \label{ref-5.3.3-80} We have for $a\in \epsilon \bar{A}\epsilon $
\[
(E+\hat{E})(a)=ae_1\otimes e_1-e_1\otimes e_1a
\]
\end{lemmas}
\begin{proof}
Both sides of the equation are derivations in $a$. Hence it suffices
to check the identity on generators for $\epsilon \bar{A}\epsilon$. It
is easy to check that these generators are given by 
\begin{align*}
t&\qquad \text{for $t\in \epsilon A \epsilon$} \\
e_{12} u &\qquad \text{for $u\in e_2 A \epsilon$}\\
v e_{21} &\qquad \text{for $v\in \epsilon A e_2$}\\
e_{12} w e_{21} &\qquad \text{for $w\in e_2 A e_2$}
\end{align*}
We compute
\begin{align*}
(E+\hat{E})(t)&=\bar{E}_1(t)+(e_{12}\bar{E}_2e_{21})(t)\\
&=te_1\otimes e_1-e_1\otimes e_1t
\end{align*}
\begin{align*}
(E+\hat{E})(e_{12}u)&=e_{12} \bar{E}_1(u)+e_{12}(e_{12}\bar{E}_2e_{21})(u)\\
&=        e_{12}  ue_1\otimes e_{1}  
-e_{1}\otimes e_{12}u \\
\end{align*}
\begin{align*}
(E+\hat{E})(ve_{21})&=\bar{E}_1(v)e_{21}+(e_{12}\bar{E}_2e_{21})(v)e_{21}\\
&=       - e_1\otimes e_1 ve_{21}+              ve_{21}\otimes e_{1}
\end{align*}
\begin{align*}
(E+\hat{E})(e_{12}we_{21})&=e_{12}\bar{E}_1(w)e_{21}
+e_{12}(e_{12}\bar{E}_2e_{21})(w)e_{21}\\
&=    e_{12} w e_{21} \otimes e_1-e_1\otimes e_{12}we_{21}
\end{align*}
In each of the cases we find the correct result.
\end{proof}
We now compute the Schouten brackets between the operators $E,\hat{E}$.
\begin{align*}
\ldb E,E\rdb &=\ldb \bar{E}_1,\bar{E}_1\rdb\\
&=\overline{\ldb E_1,E_1\rdb}\\
&=\overline{E_1\otimes e_1-e_1\otimes E_1\strut}
\qquad (\text{by \eqref{ref-3.9-54}})\\
&=E\otimes e_1-e_1\otimes E \\*[0.05in]
\ldb E,\hat{E}\rdb&=e_{12} \ldb \bar{E}_1,\bar{E}_2\rdb e_{21}\\
&=e_{12} \,\overline{\ldb E_1,E_2\rdb}\, e_{21}\\
&=e_{12}\, \overline{(E_2 e_1\otimes e_1-e_1\otimes e_1 E_2)\strut}\, e_{21}\\
&=0\\*[0.05in]
\ldb \hat{E},\hat{E}\rdb&=e_{12} 
(e_{12}\ast\ldb \bar{E}_2,\bar{E}_2\rdb \ast e_{21})e_{21}\\
&=e_{12} (e_{12}\ast\overline{\ldb E_2,E_2\rdb}\ast e_{21}) e_{21}\\
&=e_{12} (e_{12}\ast (\bar{E}_2\otimes e_2-\bar{E}_2\otimes e_2)
\ast e_{21}) e_{21}\\
&=e_{12}\bar{E}_2 e_{21}\otimes e_1
-e_1\otimes e_{12}\bar{E}_2 e_{21}\\
&=\hat{E}\otimes e_1-e_1\otimes \hat{E}
\end{align*}
We also need the following Schouten bracket.
\begin{align*}
\ldb E\hat{E},E\hat{E}\rdb&=\ldb E\hat{E},E\rdb\bar{E}-
E\ldb E\hat{E},\hat{E}\rdb\\
&=(\ldb E,E\rdb\ast\hat{E})\hat{E}-E(E\ast\ldb \hat{E},\hat{E}\rdb)\\
&=E\hat{E}\otimes \hat{E}+\hat{E}\otimes E\hat{E}
+E\hat{E}\otimes E+E\otimes E\hat{E}
\end{align*}
Hence 
\[
\{ E\hat{E},E\hat{E}\}=2E \hat{E}^2+2E^2\hat{E} \quad
\operatorname{mod}\ [D_{\bar{B}} \bar{A},D_{\bar{B}} \bar{A}]
\]
For $P\in (D_BA)_2$ we compute
\begin{align*}
\ldb E\hat{E},\bar{P}\rdb&=-\ldb E,\bar{P}\rdb\ast\hat{E} +E\ast\ldb \hat{E},\bar{P}\rdb\\
&=-\overline{\ldb E_1,P\rdb}
\ast \hat{E}
+E\ast (e_{12}\ast \overline{\ldb E_2 ,P\rdb} \ast e_{21})\\
&=-(\bar{P}e_1\otimes e_1-e_1\otimes e_1\bar{P})
\ast \hat{E}
+
E\ast (e_{12}\ast (\bar{P}e_2\otimes e_2
- e_2\otimes e_2\bar{P})\ast e_{21})\\
&=-\bar{P}\hat{E}\otimes e_1 +\hat{E}\otimes e_1 \bar{P}
+\bar{P}e_{21}\otimes Ee_{12}
-e_{21}\otimes Ee_{12}\bar{P}
\end{align*}
Hence 
\[
\{ E\hat{E},\bar{P}\}=0 \quad
\operatorname{mod}\ [D_{\bar{B}} \bar{A},D_{\bar{B}} \bar{A}]
\]
\begin{proof}[Proof of Theorem \ref{ref-5.3.1-77}]
  Note that $E_1^f=E$, $E_2^f=\hat{E}$ and $E_1^fE_2^f=E\hat{E}$ (here
  and below we view $E,\hat{E},E\hat{E}$ as elements of
  $D_{\epsilon B\epsilon}(\epsilon\overline{A}\epsilon)$).

 The result
of Lemma \ref{ref-5.3.3-80} may be rewritten as 
\begin{equation}
\label{ref-5.3-81}
E_1^f+E_2^f=F_1
\end{equation}
 Using the fact that 
$\epsilon \bar{E}_2 =0$ we have
\[
(E^n_2)^f=\Tr(\bar{E}^n_2)=e_{12} \bar{E}^n_2 e_{21}=\hat{E}^n
\]
A similar computation yields.
\[
(E^n_1)^f=\Tr(\bar{E}^n_1)=E^n
\]
Finally we have for $i>2$
\[
(E^n_i)^f=F_i^n
\]
Applying $(-)^f$ to the identity
\[
\{P,P\}=\frac{1}{6}\sum_i E^3_i
\]
yields
\[
\{P^f,P^f\}=\frac{1}{6}E^3+\frac{1}{6}\hat{E}^3+\frac{1}{6}\sum_{i>2} F^3_i
\]
We  compute (modulo commutators)
\begin{align*}
\left\{P^f-\frac{1}{2} E^f_1E^f_2,P^f-\frac{1}{2} E^f_1E^f_2\right\}&=
\left\{P^f-\frac{1}{2} E\hat{E},P^f-\frac{1}{2} E\hat{E}\right\}\\
&=\Tr \left\{\bar{P} -\frac{1}{2} E\hat{E},\bar{P} -\frac{1}{2} E\hat{E}\right\}\\
&=\Tr \{\bar{P},\bar{P}\}+\frac{1}{2} \Tr (E\hat{E}^2)
+\frac{1}{2} \Tr(E^2\hat{E})\\
&=\{P^f,P^f\}+\frac{1}{2}E \hat{E}^2+\frac{1}{2} E^2\hat{E}\\
&=\frac{1}{6} (E+\hat{E})^3+\frac{1}{6}\sum_{i>2} F^3_i\\
&=\frac{1}{6} F_1^3+\frac{1}{6}\sum_{i>2} F^3_i\\
&=\frac{1}{6} \sum_{i\neq 2} F_i^3
\end{align*}
This finishes the proof.
\end{proof}

\begin{proof}[Proof of Theorem \ref{ref-5.3.2-79}]
We need to prove
\[
\{P^{f\!\!f},\Phi^{f\!\!f}_i\}=-\frac{1}{2}(F_i\Phi_i+\Phi_iF_i)
\]
Since the case $i>2$ is easy we assume $i=1$. In that case we have to prove.
\begin{equation}
\label{ref-5.4-82}
\left\{P^{f}-\frac{1}{2}E_1^fE_2^f,\Phi^f_1\Phi_2^f\right\}=
-\frac{1}{2}(F_1\Phi^f_1\Phi_2^f +\Phi^f_1\Phi_2^f F_1)
\end{equation}
We compute the left hand side of this equation. We have
\[
\left\{P^{f}-\frac{1}{2}E_1^fE_2^f,\Phi^f_1\Phi_2^f\right\}=
\Tr \left\{\bar{P}-\frac{1}{2}E\hat{E}, \Phi\hat{\Phi}\right\}
\]
where $\Phi=\bar{\Phi}_1$, $\hat{\Phi}=e_{12}\bar{\Phi}_2e_{21}$.

We compute \def\hP{\bar{P}}
\[
\{\hP,\Phi\hat{\Phi}\}=\{\hP,\Phi\}\hat{\Phi}+\Phi\{\hP,\hat{\Phi}\}
\]
where
\begin{align*}
\{\hP,\Phi\}&=\overline{\{P,\Phi_1\}}\\
&=-\frac{1}{2}\overline{(E_1\Phi_1+\Phi_1E_1)\strut}\\
&=-\frac{1}{2}(E\Phi+\Phi E)
\end{align*}
and
\begin{align*}
\{\hP,\hat{\Phi}\}&=\{\hP,e_{12} \bar{\Phi}_2 e_{21}\}\\
&=e_{12} \,\overline{\{P,\Phi_2\}} \,e_{21}\\
&=-\frac{1}{2}e_{12} \overline{(E_2\Phi_2+\Phi_2 E_2)\strut}\, e_{21}\\
&=-\frac{1}{2}(\hat{E}\hat{\Phi}+\hat{\Phi}\hat{E})
\end{align*}
Taking things together we find
\begin{equation}
\label{ref-5.5-83}
\begin{split}
\{\hP,\Phi\hat{\Phi}\}&=-\frac{1}{2}(E\Phi+\Phi E)\hat{\Phi}
-\frac{1}{2}\Phi(\hat{E}\hat{\Phi}+\hat{\Phi}\hat{E})\\
&=-\frac{1}{2}E\Phi\hat{\Phi}-\frac{1}{2}\Phi E\hat{\Phi}
-\frac{1}{2}\Phi\hat{E}\hat{\Phi}-\frac{1}{2}\Phi\hat{\Phi}\hat{E}
\end{split}
\end{equation}
Next we compute $\{E\hat{E},\Phi\hat{\Phi}\}$. We need the following
preliminary results.
\begin{align*}
\ldb E,\Phi\rdb &= \overline{\ldb E_1,\Phi_1\rdb}\\
&=\Phi\otimes e_1-e_1\otimes \Phi\\*[0.05in]
\ldb E,\hat{\Phi}\rdb&=e_{12}\,\overline{\ldb E_1,\Phi_2\rdb}\, e_{21}\\
&=e_{12}\,\overline{(\Phi_2\otimes e_1-e_1\otimes \Phi_2)}\, e_{21}\\
&=0\\*[0.05in]
\ldb \hat{E},\Phi\rdb &=e_{12}\ast \overline{\ldb E_2,\Phi_1\rdb}\, \ast e_{21}\\
&=e_{12} \ast\,\overline{ (\Phi_1\otimes e_2-e_2\otimes \Phi_1)}\,\ast e_{21}\\
&=0\\*[0.05in]
\ldb \hat{E},\hat{\Phi}\rdb &=e_{12} (e_{12}\ast\overline{\ldb E_2,\Phi_2\rdb}
\ast e_{21})e_{21}\\
&=e_{12} (e_{12}\ast \,\overline{(\Phi_2\otimes e_2-e_2\otimes \Phi_2)}
\,\ast e_{21})e_{21}\\
&=\hat{\Phi}\otimes e_1-e_1\otimes \hat{\Phi}
\end{align*}
We then compute
\begin{align*}
\ldb E\hat{E},\Phi\hat{\Phi}\rdb &=\Phi\ldb E\hat{E},\hat{\Phi}\rdb
+\ldb E\hat{E},\Phi\rdb \hat{\Phi}\\
&=\Phi(E\ast \ldb \hat{E},\hat{\Phi}\rdb-\ldb E,\hat{\Phi}\rdb \ast
\hat{E})+(E\ast \ldb \hat{E},\Phi\rdb-\ldb E,\Phi\rdb \ast \hat{E})\hat{\Phi}\\
&=\Phi(E\ast (\hat{\Phi}\otimes e_1-e_1\otimes \hat{\Phi}))
-((  \Phi\otimes e_1-e_1\otimes \Phi             ) \ast \hat{E})\hat{\Phi}\\
&=\Phi\hat{\Phi}\otimes E-\Phi\otimes E\hat{\Phi}
-\Phi\hat{E}\otimes \hat{\Phi}+\hat{E}\otimes \Phi\hat{\Phi}
\end{align*}
and hence
\begin{equation}
\label{ref-5.6-84}
\{ E\hat{E},\Phi\hat{\Phi}\}=\Phi\hat{\Phi}E-\Phi E\hat{\Phi}
-\Phi\hat{E}\hat{\Phi}+\hat{E}\Phi\hat{\Phi}
\end{equation}
Combining \eqref{ref-5.5-83} and \eqref{ref-5.6-84} we obtain
\begin{align*}
\left\{\hP-\frac{1}{2} E\hat{E} ,\Phi\hat{\Phi}\right\}&=
-\frac{1}{2}E\Phi\hat{\Phi}-\frac{1}{2}\Phi E\hat{\Phi}
-\frac{1}{2}\Phi\hat{E}\hat{\Phi}-\frac{1}{2}\Phi\hat{\Phi}\hat{E}
\\
&\qquad-\frac{1}{2}\Phi\hat{\Phi} E+\frac{1}{2}\Phi E\hat{\Phi}
+\frac{1}{2}\Phi\hat{E}\hat{\Phi}-\frac{1}{2}\hat{E} \Phi\hat{\Phi}\\
&=-\frac{1}{2}E\Phi\hat{\Phi}
-\frac{1}{2}\Phi\hat{\Phi}E
-\frac{1}{2}\hat{E} \Phi\hat{\Phi}
-\frac{1}{2} \Phi\hat{\Phi}\hat{E}
\end{align*}
and hence we obtain
\begin{align*}
\operatorname{LHS}\eqref{ref-5.4-82} 
&=\Tr \left\{\bar{P}-\frac{1}{2}E\hat{E}, \Phi\hat{\Phi}\right\}\\
&=-\frac{1}{2}E_1^f\Phi_1^f\Phi_2^f-\frac{1}{2}\Phi_1^f\Phi_2^fE_1^f
-\frac{1}{2}E_2^f \Phi_1^f\Phi_2^f-\frac{1}{2}\Phi_1^f\Phi_2^fE_2^f \\
&=  -\frac{1}{2} F_1 \Phi_1^f\Phi_2^f -\frac{1}{2} \Phi_1^f\Phi_2^f  F_1\\
&=\operatorname{RHS}\eqref{ref-5.4-82}\qed
\end{align*}
\def\qed{}\end{proof}
\section{Quivers}
\label{ref-6-85}
\subsection{Generalities}
Below $Q=(Q,I,h,t)$ is a finite quiver with vertex set
$I=\{1,\ldots,n\}$ and edge set $Q$. The maps $t,h:Q\r I$ associate
with every edge its start and end. We extend the definitions of $h,t$
to paths in $Q$.  By $e_i$ we denote the idempotent associated to the
vertex $i$ and we put $B=\oplus_i ke_i$.  We let $\bar{Q}$ be the
double of $Q$. $\bar{Q}$ is obtained from $Q$ by adjoining for every
arrow $a$ an opposite arrow $a^\ast$. We define $\epsilon:\bar{Q}\r
\{\pm 1\}$ as the function which is $1$ on $Q$ and $-1$ on
$\bar{Q}-Q$. By $kQ$ we denote the path algebra of $Q$ (with
multiplication given by concatenation of paths).  Note that $kQ/B$ is
smooth so we don't have to make a difference between differential and
ordinary notions (see \S\ref{ref-5.2-75}) in the case of quivers.

\subsection{Vector fields and the Schouten bracket}
Let $A=kQ$.  For $a\in Q$ we define the element $\frac{\partial\
}{\partial a}\in D_BA$ which on $b\in Q$ acts as
\[
\frac{\partial b}{\partial a}
=
\begin{cases}
e_{t(a)}\otimes e_{h(a)}&\text{if $a=b$}\\
0&\text{otherwise}
\end{cases}
\]
It is clear that $D_{A/B}$ is generated by $\left(\frac{\partial\
  }{\partial a}\right)_{a\in Q}$ as an $A$-bimodule. Hence $D_BA$ is the
tensor algebra over $A$ generated by $\left(\frac{\partial\ }{\partial
    a}\right)_a$.
\begin{propositions}
 Let  $a,b\in Q$. Then
\begin{align*}
\ldb a,b\rdb &=0 \\
\ldbgg \frac{\partial\ }{\partial a},b\rdbgg&=
\begin{cases}
e_{t(a)}\otimes e_{h(a)}&\text{if $a=b$}\\
0&\text{otherwise}
\end{cases}\\
\ldbgg \frac{\partial\ }{\partial a},\frac{\partial\ }{\partial b}\rdbgg&=0
\end{align*}
\end{propositions}
\begin{proof} Only the third equality is not immediately obvious. But
  a quick check of the definitions reveals that $ \ldb \frac{\partial\
  }{\partial a},\frac{\partial\ }{\partial b}\rdb(c)=0$ for any $c\in
  Q$.
\end{proof}
\begin{propositions}  
\begin{enumerate}
\item For $\delta\in D_{A/B}$ we have the equality in $D_BA$:
\begin{equation}
\label{ref-6.1-86}
\delta=\sum_{a\in Q} \delta(a)''\frac{\partial\ }{\partial a} \delta(a)'
\end{equation}
\item For $i=1,\ldots, n$ we have the equality:
\begin{equation}
\label{ref-6.2-87}
E_i=\sum_{a\in Q, h(a)=i} \frac{\partial\ }{\partial a}a
-\sum_{a\in Q, t(a)=i} a\frac{\partial\ }{\partial a}
\end{equation}
\end{enumerate}
\end{propositions}
\begin{proof}
\begin{enumerate}
\item
  Let $b\in Q$. Evaluated on $b$ \eqref{ref-6.1-86} can be
  rewritten as
\[
\delta(b)=\delta(b)''*(e_{t(b)}\otimes e_{h(b)})* \delta(b)'
\]
The right hand side of this equation is equal to
$e_{t(b)}\delta(b)'\otimes \delta(b)''e_{h(b)}=\delta(e_{t(b)}be_{h(b)})=\delta(b)$.
\item If we substitute $\delta=E_i$ in \eqref{ref-6.1-86} then we obtain
\begin{align*}
E_i&=\sum_{a\in Q}e_i \frac{\partial\ }{\partial a} ae_i
-e_ia \frac{\partial\ }{\partial a} e_i\\
&=\sum_{a\in Q, h(a)=i} \frac{\partial\ }{\partial a}a
-\sum_{a\in Q, t(a)=i} a\frac{\partial\ }{\partial a}
\end{align*}
\end{enumerate}
\end{proof}
\begin{remarks}
  The expression for $E_i$ can be conveniently rewritten as follows.
  Put $E=\sum_i E_i$. Then
\[
E=\sum_{a\in Q} \left[\frac{\partial\ }{\partial a},a\right]
\]
\end{remarks}

\subsection{Hamiltonian structure}
\begin{theorems}
\label{ref-6.3.1-88}
 $A=k\bar{Q}$ has a Hamiltonian structure given by 
\begin{equation}
\label{ref-6.3-89}
P=\sum_{a\in Q} \frac{\partial\ }{\partial a}
\frac{\partial\ }{\partial a^\ast}
\end{equation}
\[
\mu=\sum_{a\in Q} [a,a^\ast]
\]
\end{theorems}
\begin{proof} The fact that $\{P,P\}=0$ is trivial. For the moment map
property we compute
\begin{align*}
\ldb P,a\rdb&=-(e_{t(a)}\otimes e_{h(a)})\ast \frac{\partial\ }{\partial a^\ast}\\
\ldb P,a^\ast\rdb &=\frac{\partial\ }{\partial a}\ast (e_{h(a)}\otimes e_{t(a)})
\end{align*}
whence
\begin{align*}
\{ P,a\}&= -\frac{\partial\ }{\partial a^\ast}\\
\{ P,a^\ast\}&=\frac{\partial\ }{\partial a}
\end{align*}
Thus
\begin{align*}
\{P,\mu\}&=\sum_{a\in Q} [\{P,a\},a^\ast]+[a,\{P,a^\ast\}]\\
&=\sum_{a\in Q} -\left[\frac{\partial\ }{\partial a^\ast},a^\ast\right]+\left[a,\frac{\partial\ }{\partial a}\right]\\
&=-E\qed
\end{align*}
\def\qed{}\end{proof}
\subsection{The necklace Loday algebra}
\label{ref-6.4-90} A simple computation yields the following formula
for the induced double Poisson bracket on $k\bar{Q}$. Let $a\in Q$.
\begin{align*}
\ldb a,a^\ast \rdb_P&=e_{h(a)}\otimes e_{t(a)}\\
\ldb a^\ast,a \rdb_P&=-e_{t(a)}\otimes e_{h(a)}
\end{align*}
and all other double brackets are zero. 

The double bracket on paths is pictorially given as follows:
\[
\psfrag{a}[][]{$a$}
\psfrag{ast}[][]{$a^\ast$}
\psfrag{x}[][]{$x$}
\psfrag{y}[][]{$y$}
\psfrag{z}[][]{$z$}
\psfrag{t}[][]{$t$}
\left\{\!\!\left\{
\raisebox{-1.5cm}{\includegraphics[height=3cm]{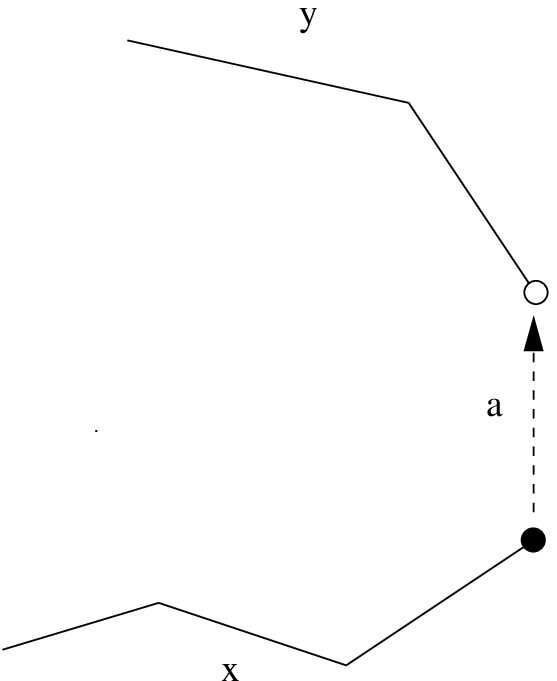}}\qquad,\qquad
\raisebox{-1.5cm}{\includegraphics[height=3cm]{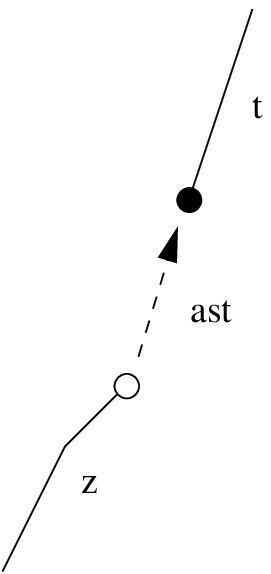}}
\right\}\!\!\right\}=
\]
\[
\psfrag{a}[][]{$a$}
\psfrag{ast}[][]{$a^\ast$}
\psfrag{x}[][]{$x$}
\psfrag{y}[][]{$y$}
\psfrag{z}[][]{$z$}
\psfrag{t}[][]{$t$}
\sum_{a\in \bar{Q}}\epsilon(a)
\quad\raisebox{-1.5cm}{\includegraphics[height=3cm]{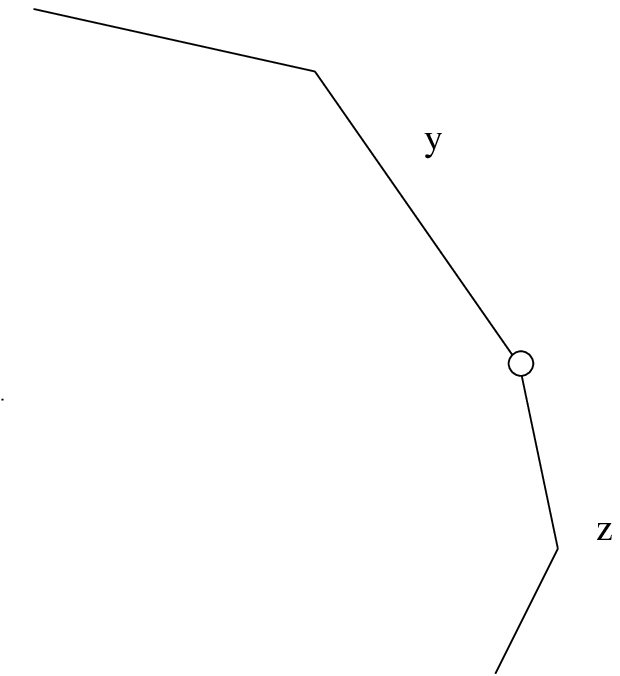}}\qquad
\otimes 
\quad\raisebox{-1.5cm}{\includegraphics[height=3cm]{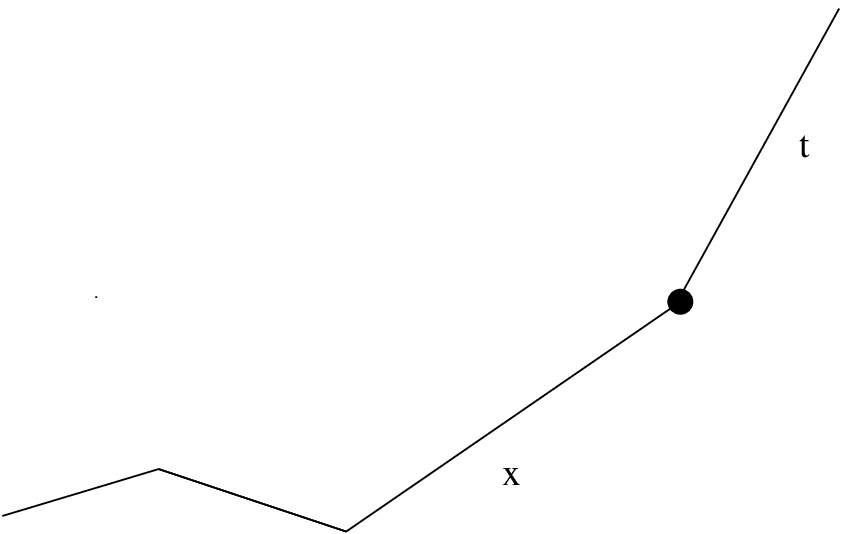}}\qquad
\]
where the black and white dots represent identical vertices.

The double Poisson structure on $k\bar{Q}$ induces a left Loday
algebra structure $\{-,-\}$ on $k\bar{Q}$. 
\begin{propositions} For oriented paths $x,y$ in $\bar{Q}$ we have
$\{x,y\}=0$ if $x$ is not closed. Otherwise the bracket can be pictorially
represented as follows
\[
\psfrag{a}[][]{$a$}
\psfrag{ast}[][]{$a^\ast$}
\left\{
\raisebox{-1.5cm}{\includegraphics[height=3cm]{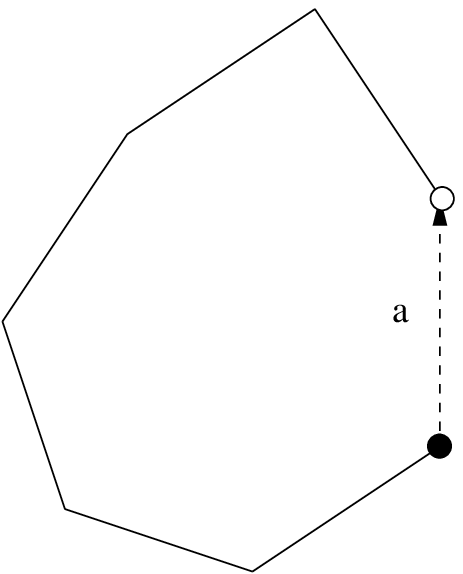}}\qquad,\qquad
\raisebox{-1.5cm}{\includegraphics[height=3cm]{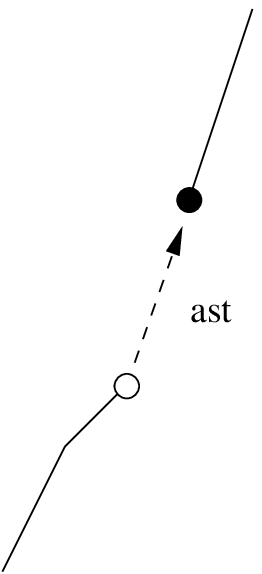}}
\right\}=
\]
\[
\sum_{a\in \bar{Q}}\epsilon(a)\quad\raisebox{-1.5cm}{\includegraphics[height=3cm]{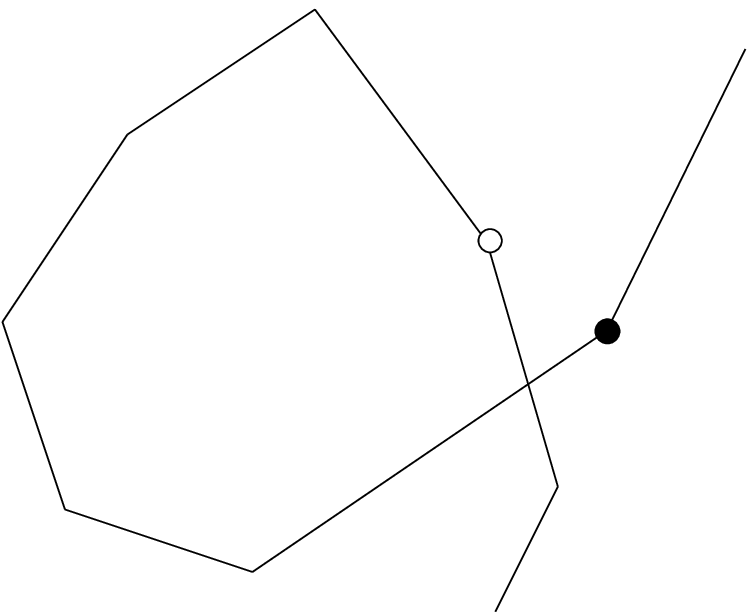}}
\]
\end{propositions}
If we restrict this bracket to closed paths we obtain the so-called
\emph{necklace Lie algebra} structure on $k\bar{Q}/[k\bar{Q},k\bar{Q}]$
\cite{LebBock,Ginzburg1,Kosymp}.
\subsection{Quasi-Hamiltonian structure for a very simple quiver}
In this section we consider the   quiver $Q$ given by
\[
\xymatrix{
1 \ar@/^/[rr]^a & &2\ar@/^/[ll]^{a^\ast}
}
\]
We let $A$ be the path algebra of $kQ$ with $e_1+aa^\ast$ and
$e_2+a^\ast a$ inverted. By inverted we mean that we introduce
elements $I$, $J$ such that $I=Ie_1=e_1I$, $J=Je_2=e_2J$ and
$I(e_1+aa^\ast)=(e_1+aa^\ast)I=e_1$ and $J(e_2+a^\ast a)=(e_2+a^\ast
a)J=e_2$. Below we denote use the notation $(e_1+aa^\ast)^{-1}$
and $(e_2+a^\ast a)^{-1}$ for $I$ and $J$ (commiting an abuse of notation).
\begin{theorems}
\label{ref-6.5.1-91} $A$ has a quasi-Hamiltonian structure given by
\begin{align}
P&=\frac{\partial\ }{\partial a}
\frac{\partial\ }{\partial a^\ast}+
\frac{1}{2}\left(a\frac{\partial\ }{\partial a}\frac{\partial\ }{\partial a^\ast}a^\ast - a^\ast\frac{\partial\ }{\partial a^\ast} \frac{\partial\ }{\partial a}a\right)\label{ref-6.4-92}\\
&=\frac{1}{2}(1+a^\ast a) \frac{\partial\ }{\partial a}\frac{\partial\ }{\partial a^\ast}
-\frac{1}{2} (1+a a^\ast)\frac{\partial\ }{\partial a^\st}\frac{\partial\ }{\partial a}\qquad\text{\rm mod $[-,-]$} \label{ref-6.5-92a}
\\*[0.1in]
\Phi&=(1+aa^\ast)(1+a^\ast a)^{-1}
\end{align}
\end{theorems}
Note that the partial derivatives have odd degree. So their commutator has
a plus sign. This explains why the $1$'s in \eqref{ref-6.5-92a} do not cancel. 
\begin{proof}
We first consider the quasi-Poisson structure. For simplicity we introduce
the following elements of $D_{A/B}$.
\begin{align*}
U&=a\frac{\partial\ }{\partial a}\\
V&=\frac{\partial\ }{\partial a}a\\
U^\ast&=a^\ast\frac{\partial\ }{\partial a^\ast}\\
V^\ast&=\frac{\partial\ }{\partial a^\ast}a^\ast
\end{align*}
Then $P$ becomes
\[
P=\frac{\partial\ }{\partial a}
\frac{\partial\ }{\partial a^\ast}+\frac{1}{2}(UV^\ast-U^\ast V)
\]
We have to prove 
\begin{equation}
\label{ref-6.7-93}
\{P,P\}=\frac{1}{6}(E_1^3+E_2^3)\qquad \text{mod $[-,-]$}
\end{equation}
By \eqref{ref-6.2-87} we have
\begin{align*}
E_1&=\frac{\partial\ }{\partial a^\ast}a^\ast -a
\frac{\partial\ }{\partial a}=V^\ast-U\\
E_2&=\frac{\partial\ }{\partial a}a -a^\ast
\frac{\partial\ }{\partial a^\ast}=V-U^\ast
\end{align*}
We compute
\begin{align*}
\ldb UV^\ast,UV^\ast\rdb&=\ldb UV^\ast,U\rdb V^\ast-U\ldb UV^\ast,V^\ast \rdb\\
&=(\ldb U ,U\rdb\ast V^\ast)V^\ast-U(U\ast \ldb V^\ast,V^\ast\rdb)
\\
&=(( e_1\otimes U-U\otimes e_1) \ast V^\ast)V^\ast-
U(U\ast  ( -e_1\otimes V^\ast+V^\ast\otimes e_1                    )          )\\
&=
-V^\ast\otimes UV^\ast - UV^\ast\otimes  V^\ast+U\otimes UV^\ast+UV^\ast\otimes U
\end{align*}
Here $\ldb U,U\rdb$ and $\ldb V^\ast,V^\ast\rdb$ have been computed using Lemma
\ref{ref-6.5.2-94} below.
Hence we obtain
\[
\{ UV^\ast,UV^\ast\}=-2U(V^\ast)^2+2U^2V^\ast\qquad \text{mod $[-,-]$}
\]
Similarly
\[
\{U^\ast V,U^\ast V\}=-2U^\ast V^2+2(U^\ast)^2 V\qquad \text{mod $[-,-]$}
\]
It is also clear from Lemma \ref{ref-6.5.2-94} that $\ldb UV^\ast , U^\ast V\rdb=0$. Hence
\[
\{UV^\ast , U^\ast V\}=0
\]
We need some more computations
\begin{align*}
\ldbgg UV^\ast,\frac{\partial\ }{\partial a}
\frac{\partial\ }{\partial a^\ast}\rdbgg &=
\ldbgg UV^\ast,\frac{\partial\ }{\partial a}
\rdbgg \frac{\partial\ }{\partial a^\ast}
-
\frac{\partial\ }{\partial a}\ldbgg UV^\ast,
\frac{\partial\ }{\partial a^\ast}\rdbgg
\\
&=\Bigl(\ldbgg U,\frac{\partial\ }{\partial a}
\rdbgg \ast V^\ast\Bigr)\frac{\partial\ }{\partial a^\ast}
-
\frac{\partial\ }{\partial a}\Bigl(U\ast \ldbgg V^\ast,
\frac{\partial\ }{\partial a^\ast}\rdbgg\Bigr)\\
&=\Bigl(\Bigl(-\frac{\partial\ }{\partial a}\otimes e_1\Bigr)
\ast V^\ast\Bigr)\frac{\partial\ }{\partial a^\ast}
-
\frac{\partial\ }{\partial a}\Bigl(U\ast 
\Bigl(
-e_1\otimes \frac{\partial\ }{\partial a^\ast}
\Bigr)\Bigr)\\
&=-\frac{\partial\ }{\partial a} V^\ast\otimes \frac{\partial\ }{\partial a^\ast}
+\frac{\partial\ }{\partial a}\otimes U\frac{\partial\ }{\partial a^\ast}
\end{align*}
Hence 
\begin{align*}
\Bigl\{ UV^\ast,\frac{\partial\ }{\partial a}
\frac{\partial\ }{\partial a^\ast}\Bigr\}&=-\frac{\partial\ }{\partial a} V^\ast\frac{\partial\ }{\partial a^\ast}
+\frac{\partial\ }{\partial a} U\frac{\partial\ }{\partial a^\ast}\\
&=-\frac{\partial\ }{\partial a} \frac{\partial\ }{\partial a^\ast}a^\ast
\frac{\partial\ }{\partial a^\ast}
+\frac{\partial\ }{\partial a} a\frac{\partial\ }{\partial a} 
\frac{\partial\ }{\partial a^\ast}
\end{align*}

Similarly, computing modulo commutators we obtain
\begin{align*}
\Bigl\{ U^\ast V,\frac{\partial\ }{\partial a}
\frac{\partial\ }{\partial a^\ast}\Bigr\}&=
-\Bigl\{ U^\ast V,\frac{\partial\ }{\partial a^\ast}
\frac{\partial\ }{\partial a}\Bigr\}\\
&=\frac{\partial\ }{\partial a^\ast} \frac{\partial\ }{\partial a}a
\frac{\partial\ }{\partial a}
-\frac{\partial\ }{\partial a^\ast} a^\ast\frac{\partial\ }{\partial a^\ast} 
\frac{\partial\ }{\partial a}
\end{align*}
If follows that
\[
\Bigl\{UV^\ast-U^\ast V,\frac{\partial\ }{\partial a}
\frac{\partial\ }{\partial a^\ast}\Bigr\}=0\qquad \text{mod $[-,-]$}
\]
Combining everything we find
\[
\{P,P\}=\frac{1}{2}(-U(V^\ast)^2+U^2V^\ast-U^\ast V^2+(U^\ast)^2V)\qquad \text{mod
$[-,-]$}
\]
On the other hand we have
\[
\begin{split}
E_1^3&=(V^\ast)^3-U^3-3U(V^\ast)^2+3U^2 V^\ast\\
E_2^3&=V^3-(U^\ast)^3-3U^\ast V^2+3(U^\ast)^2 V
\end{split}
\qquad \text{mod $[-,-]$}
\]
and also
\[
\begin{split}
U^3&=V^3\\
(U^\ast)^3&=(V^\ast)^3
\end{split}
\qquad \text{mod $[-,-]$}
\]
It follows
\[
\frac{1}{6}(E_1^3+E_2^3)=\frac{1}{2}(-U(V^\ast)^2+U^2 V^\ast-U^\ast V^2+(U^\ast)^2 V)\qquad \text{mod $[-,-]$}
\]
So it follows that \eqref{ref-6.7-93} is indeed true.

\medskip

Now we prove that $\Phi$ is a multiplicative moment map.  We have 
$\Phi=\Phi_1+\Phi_2$ where 
\begin{align*}
\Phi_1&=e_1+aa^\ast\\
\Phi_2&=(e_2+a^\ast a)^{-1}
\end{align*}
We have to prove
\[
\{P,\Phi_i\}=-\frac{1}{2}(E_i\Phi_i+E_i\Phi_i)
\]
 We will first consider $\Phi_1$.
We compute
\begin{align*}
\ldb UV^\ast,aa^\ast\rdb&=\ldb UV^\ast,a\rdb a^\ast+a \ldb UV^\ast,a^\ast\rdb \\
&=-(\ldb U,a\rdb\at V^\ast) a^\ast +a(U\ast \ldb V^\ast,a^\ast\rdb)\\
&=
-((e_1\otimes a)\ast V^\ast)a^\ast+a(U\ast (a^\ast\otimes e_1))\\
&=-V^\ast\otimes aa^\ast+aa^\ast\otimes U
\end{align*}
So
\begin{align*}
\{UV^\ast,aa^\ast\}&=-V^\ast aa^\ast+aa^\ast U\\
&=-\frac{\partial\ }{\partial a^\ast}a^\ast aa^\ast+aa^\ast a\frac{\partial\ }{\partial a}
\end{align*}
Similarly
\begin{align*}
\ldb U^\ast V,aa^\ast\rdb
&=\ldb U^\ast V,a\rdb a^\ast+a\ldb U^\ast  V,a^\ast\rdb\\
&=(U^\ast \ast \ldb V,a\rdb)a^\ast-a(\ldb U^\ast,a^\ast\rdb\ast V)\\
&=(U^\ast \ast (a\otimes e_2))a^\ast-a((e_2\otimes a^\ast)\ast V)\\
&=a\otimes U^\ast a^\ast-aV\otimes a^\ast
\end{align*}
which yields
\begin{align*}
\{ U^\ast V,aa^\ast\}&=aU^\ast a^\ast-aVa^\ast\\
&=aa^\ast\frac{\partial\ }{\partial a^\ast} a^\ast-a
\frac{\partial }{\partial a}aa^\ast
\end{align*}
We obtain
\begin{align*}
\{UV^\ast-U^\ast V,aa^\ast\}&=-\frac{\partial\ }{\partial a^\ast}a^\ast aa^\ast+aa^\ast a\frac{\partial\ }{\partial a}
-
aa^\ast\frac{\partial\ }{\partial a^\ast} a^\ast+a
\frac{\partial }{\partial a}aa^\ast\\
&=-E_1 aa^\ast-aa^\ast E_1
\end{align*}
We also have
\begin{align*}
\ldbgg \frac{\partial\ }{\partial a}\frac{\partial\ }{\partial a^\ast},
aa^\ast
\rdbgg &=
\ldb \frac{\partial\ }{\partial a}\frac{\partial\ }{\partial a^\ast},
a
\rdbgg a^\ast+
a\ldbgg \frac{\partial\ }{\partial a}\frac{\partial\ }{\partial a^\ast},
a^\ast
\rdbgg\\
&=-\Bigl((e_1\otimes e_2)\ast \frac{\partial\ }{\partial a^\ast}\Bigr)a^\ast
+a\Bigl( \frac{\partial\ }{\partial a}\ast (e_2\otimes e_1)\Bigr)\\
&=-\frac{\partial\ }{\partial a^\ast}\otimes a^\ast
+a\otimes \frac{\partial\ }{\partial a}
\end{align*}
which yields
\[
\Bigl\{ \frac{\partial\ }{\partial a}\frac{\partial\ }{\partial a^\ast},
aa^\ast
\Bigr\}=-E_1
\]
Taking everything together we obtain
\begin{align*}
\{P,\Phi_1\}&=\{P,aa^\ast\}\\
&=-E_1-\frac{1}{2}(E_1 aa^\ast+aa^\ast E_1)\\
&=-\frac{1}{2}(E_1\Phi_1+\Phi_1E_1)
\end{align*}
Now we consider $\Phi_2$. Applying the automorphism of order two
$
e_1\leftrightarrow e_2$,
$a\leftrightarrow a^\ast$
has the effect
$
P\mapsto -P\ \text{(up to commutators)}$ and
$E_1\mapsto E_2$,
$E_2\mapsto E_1$.

Hence we obtain the following identity for $\Psi=e_2+a^\ast a$:
\[
-\{P,\Psi\}=-\frac{1}{2}(E_2\Psi+\Psi E_2)
\]
Since $\Phi_2=\Psi^{-1}$ and $\{P,-\}$ is a derivation we obtain
\begin{align*}
\{P,\Phi_2\}&=-\Psi^{-1}\{P,\Psi\}\Psi^{-1}\\
&=-\frac{1}{2}(\Psi^{-1} E_2+E_2\Psi^{-1} )\\
&=-\frac{1}{2}(\Phi_2 E_2+E_2\Phi_2 )\qed
\end{align*}
\def\qed{}
\end{proof}

\begin{lemmas}
\label{ref-6.5.2-94}
Let $Q$ be an arbitrary quiver, $a\in Q$ and $h(a)\neq t(a)$ (i.e.\
$a$ is not a loop). Put $e=e_{t(a)}$, $f=e_{h(a)}$. 
If $X$ is in the subring of $D_BA$ generated by $a$ and 
$\frac{\partial\ }{\partial a}$ then
\begin{align*}
\ldb U,X\rdb &= - Xe\otimes e+e\otimes eX\\
\ldb V,X\rdb &= Xf\otimes f-f\otimes fX\\
\end{align*}
\end{lemmas}
\begin{proof}
By
\eqref{ref-6.2-87} we obtain
\begin{align*}
\ldb E_{t(a)},X\rdb&=-\ldb U,X\rdb\\
\ldb E_{h(a)},X\rdb&=\ldb V,X\rdb
\end{align*}
Using \eqref{ref-3.9-54} this implies what we want.
\end{proof}
\subsection{Fusion for quivers}
\label{ref-6.6-95}
We now discuss what happens if we perform fusion on  path
algebras. This will be used in the next section. Put $A=kQ$.
 It is clear that $\bar{A}$ is generated over $B$ by $a\in
Q$ and $e_{12}$, $e_{21}$, subject to the
relations $e_{12}e_{21}=e_1$, $e_{21}e_{12}=e_2$.  Then it is not hard
to see that $A^f$ is freely generated over $B^f$ by
\begin{align*}
a&\qquad\text{$h(a)\neq 2$, $t(a)\neq 2$}\\
ae_{21} &\qquad\text{$h(a)=2$, $t(a)\neq 2$}\\
e_{12}a &\qquad\text{$h(a) \neq 2$, $t(a)=2$}\\
e_{12}ae_{21} &\qquad\text{$h(a)=2$, $t(a)=2$}
\end{align*}
Now let $Q^f$ be the quiver obtained from $Q$ by ``fusing'' 
vertices $1$ and $2$. I.e. $Q^f$ has the same edges as $Q$ and
 vertices $I^f=I-\{2\}$. The maps $t,h$ are redefined as follows. 
\begin{align*}
h^f(a)=&
\begin{cases}
1&\text{if $h(a)= 2$}\\
h(a)&\text{otherwise}
\end{cases}\\
t^f(a)=&
\begin{cases}
1&\text{if $t(a)= 2$}\\
t(a)&\text{otherwise}
\end{cases}
\end{align*}
The following result is easy to prove
\begin{propositions} The map
\[
(kQ)^f\r k(Q^f)
\]
which is defined by (for $a\in Q$)
\begin{align*}
a\mapsto a&\qquad\text{$h(a)\neq 2$, $t(a)\neq 2$}\\
ae_{21}\mapsto a &\qquad\text{$h(a)=2$, $t(a)\neq 2$}\\
e_{12}a \mapsto a &\qquad\text{$h(a) \neq 2$, $t(a)=2$}\\
e_{12}ae_{21} \mapsto a&\qquad\text{$h(a)=2$, $t(a)=2$}
\end{align*}
is an isomorphism.
\end{propositions}
We need a slight extension of this result. 
\begin{propositions}
Let $S\subset \bigcup_{i,j} e_{i}
A e_j$ and let $A_S$ be the algebra obtained from $A$ by
formally adjoining for all $s\in S$ an element $s^{-1}$ which satisfies the
relations $s^{-1}s=e_{h(s)}$, $ss^{-1}=e_{t(s)}$. Then one has
\[
(A_S)^f=A^f_{S^f}
\]
where 
\[
S^f=\{s^f\mid s\in S\}
\]
\end{propositions}
\begin{proof} Left to the reader. 
\end{proof}
\subsection{Quasi-Hamiltonian structure for general quivers}
In this section we prove the following result. 
\begin{theorems} 
\label{ref-6.7.1-96}Let $A$ be obtained from $k\bar{Q}$  by inverting all elements
$(1+aa^\ast)_{a\in \bar{Q}}$. Fix an arbitrary total ordering
  on $\bar{Q}$.  Then $A$ has a quasi-Hamiltonian structure given by
\[
P=\frac{1}{2}\left(\sum_{a\in \bar{Q}}
\left(\epsilon(a)
(1
+
a^\ast a)\frac{\partial\ }{\partial a} \frac{\partial\ }{\partial a^\ast}\right)
-
\sum_{a<b\in\bar{Q}}\left(\frac{\partial }{\partial  a^\ast} a^\ast-a\frac{\partial }{\partial  a}\right)\left(
\frac{\partial }{\partial  b^\ast} b^\ast-b\frac{\partial }{\partial  b}
\right)\right)
\]
\[
\Phi=\prod_{a\in \bar{Q}}(1+aa^\ast)^{\epsilon(a)}
\]
In the definition of $\Phi$ the product is taken with respect to the
chosen ordering on~$\bar{Q}$.
\end{theorems}
\begin{proof} We will deduce this result from Theorem \ref{ref-6.5.1-91} using
fusion (as discussed in \S\ref{ref-6.6-95}).

Let $Q^{\text{sep}}$ be the quiver with the same edges as $\bar{Q}$ but
with vertices $(v_a)_{a\in \bar{Q}}$. The head and tail of an edge are
defined by
\begin{align*}
t(a)&=v_a\\
h(a)&=t(a^\ast)=v_{a^\ast}
 \end{align*}
Thus $Q^{\text{sep}}$ is a disjoint union of little quivers of the form
\[
\xymatrix{
v_a \ar@/^/[rr]^a & &v_{a^\ast}\ar@/^/[ll]^{a^\ast}
}
\]

 Let $A^{\text{sep}}$ be the algebra obtained from $kQ^{\text{sep}}$
 by inverting all $(1+aa^\ast)_{a\in Q^{\text{sep}}}$.
 By Theorem \ref{ref-6.5.1-91} $kQ^{\text{sep}}$ has a quasi-Hamiltonian
 structure $(P^{\text{sep}}, \Phi^{\text{sep}})$ where
\[
P^{\text{sep}}=\frac{1}{2}\sum_{a\in Q^{\text{sep}}}
\epsilon(a)
(1
+
a^\ast a)\frac{\partial\ }{\partial a} \frac{\partial\ }{\partial a^\ast}  
\]
and 
\[
\Phi^{\text{sep}}_{v_a}=(e_{v_a}+aa^\ast)^{\epsilon(a)}
\]
To obtain $k\bar{Q}$ from $kQ^{\text{sep}}$ we need to fuse certain vertices.
More precisely for a vertex $i\in I$ we need to fuse the vertices
$v_a$ such that $t(a)=i$. The fusing process depends on the order
in which we perform it. To fix this we fix a total ordering of all
edges in $\bar{Q}$. We put the same total ordering on the vertices $v_a$.

By Theorems \ref{ref-5.3.1-77}, \ref{ref-5.3.2-79} and \eqref{ref-6.2-87} we 
see that
fusing the vertices $v_a$ with $t(a)=i$ has the effect of adding
\[
-\frac{1}{2}\sum_{a<b\in\bar{Q}, t(a)=t(b)=i} F_a F_b
\]
to $P^{\text{sep}}$ where
\[
F_a=
\frac{\partial\ }{\partial  a^\ast} a^\ast-a\frac{\partial\ }{\partial  a}
\]
and to replace $(\Phi^{\text{sep}}_{v_a})_{t(a)=i}$ by  the product 
\[
\Phi_i=\prod_{a\in \bar{Q},t(a)=i} (e_{i}+aa^\ast)^{\epsilon(a)}
\]
where the order on the product is given by the ordering of the edges.
Performing this for all vertices in $\bar{Q}$ proves the theorem.
\end{proof}
\begin{remarks} The total ordering on the edges  of $\bar{Q}$ 
actually contains too much information. It is sufficient to order for
every vertex $i$ the edges starting in $i$.
\end{remarks}
\begin{remarks} It follows from the formulas for $P$ and $\Phi$ that
  $k\bar{Q}$ has always double quasi-Poisson structure. However in order
to have a quasi-Hamiltonian structure we need to invert 
the elements $(1+aa^\ast)_{a\in \bar{Q}}$.
\end{remarks} 
\subsection{Preprojective algebras and multiplicative preprojective
algebras}
Fix $\lambda\in B$. The algebra
\[
\Pi^\lambda=k\bar{Q}/\Bigl(\sum_{a\in Q} [a,a^\ast]-\lambda\Bigr)
\]
is the so-called ``deformed preprojective algebra'' It was introduced
by Crawley-Boevey and Holland in \cite{CH}. A multiplicative version 
was
introduced by Crawley-Boevey and Shaw in \cite{CBShaw}. Let
 Fix $q\in B^\ast$ and put
\[
\Lambda^q=k\bar{Q}_{(1+aa^\ast)_{a\in \bar{Q}}}\biggl/\biggl(\prod_{a\in \bar{Q}}
(1+aa^\ast)^{\epsilon(a)}-q\biggr)\biggr.
\]
The product is taken with respect to an arbitrary ordering
of $\bar{Q}$ but it is shown in \cite{CBShaw} that the resulting
algebra is independent of this ordering. 

Combining Propositions \ref{ref-2.6.5-41},\ref{ref-5.1.5-74},\ref{ref-6.3.1-88} and 
Theorem \ref{ref-6.5.1-91} we obtain:
\begin{propositions} \label{ref-6.8.1-97} Both the ordinary deformed
  preprojective algebra and the deformed multiplicative preprojective
  algebra have a Poisson structure (as in Definition \ref{ref-2.6.1-38}).
\end{propositions}

\section{Representation spaces}
\label{ref-7-98}
\subsection{General principles}
 We put $I=\{1,\ldots, n\}$.
Let $\alpha=(\alpha_1,\ldots,\alpha_n)\in
\NN^n$ and put $|\alpha|=\sum_i \alpha_i$.  Define the
function
\[
\phi:[\,1\,.\,.\,|\alpha|\,]\r I
\]
 by the property
\[
\phi(p)=i\iff \alpha_1+\cdots+\alpha_{i-1}+1\le p\le \alpha_1+\cdots+\alpha_{i}
\]
 Throughout we assume that  $B=k e_1\oplus\cdots \oplus k e_n$ 
is semi-simple. As usual $A$ is a finitely generated $B$ algebra.  

We view an element $X$
of $M_{|\alpha|}(k)$ as a block matrix $(X_{uv})_{uv}$ with
$u,v=1,\ldots, n$ and $X_{uv}\in M_{\alpha_u\times \alpha_v}(k)$.
We will also consider $B$ as being diagonally embedded in $M_{|\alpha|}(k)$
where $e_i$ is the identity matrix in $M_{\alpha_i\times \alpha_i}(k)$.

We define $\Rep(A,\alpha)$ as the affine scheme representing the functor
\[
R\mapsto \Hom_B(A,M_{|\alpha|}(R))
\]
from commutative $k$-algebras to sets. The coordinate ring of
$\Rep(A,\alpha)$ is generated by symbols $a_{pq}$ for $a\in A$,
$p,q=1,\ldots,|\alpha|$ which are linear in $a$ and satisfy the relations
\begin{align*}
a_{pq}b_{qr}&=(ab)_{pr}\\
(e_i)_{pq}&=\delta_{pq}\delta_{\phi(p),i}
\end{align*}
A map $f\in \Hom_B(A,M_{|\alpha|}(R))$ corresponds to the point $x\in
\Rep(A,\alpha)(R)$ if the following relation holds for $a\in A$,
$p,q=1,\ldots, |\alpha|$
\[
a_{pq}(x)=f(a)_{pq}
\]
Below we identify $\Rep(A,\alpha)(R)$ and $\Hom_B(A, M_{|\alpha|}(R))$.

For $a\in A$ it will be convenient to introduce the
$M_{|\alpha|}(k)$-valued function $X(a)$ on $\Rep(A,\alpha)$ by the
rule $X(a)_{ij}=a_{ij}$.  The defining relations on $\Rep(A,\alpha)$
may then be written as
\begin{align*}
X(ab)&=X(a)X(b)\\
X(e_i)&=e_i
\end{align*}

Put $\Gl_{\alpha}=\prod_i \Gl_{\alpha_i}$. $\Gl_\alpha$ acts by conjugation
on $M_{|\alpha|}$. This induces an action on $\Rep(A,\alpha)$. To work 
out what this action is let $x\in \Rep(A,\alpha)(R)=\Hom_B(A,M_{|\alpha| }(R))$. We have for $a\in A$.
\[
a_{ij}(x)=x(a)_{ij}
\]
and hence for $g\in \Gl_{\alpha}(R)$
\begin{align*}
(g\cdot a_{ij})(x)&=a_{ij}((g^{-1}-g)\circ x)\\
&=(g^{-1}x(a)g)_{ij}\\
&=(g^{-1})_{iu} x(a)_{uv} g_{vj}\\
&=(g^{-1})_{iu} a_{uv}(x) g_{vj}
\end{align*}
In terms of the $X(a)$ we may write:
\[
g\cdot X(a)=g^{-1} X(a) g
\]
where the ``$\cdot$'' means that we apply the action of $g$ entry wise.

Let $M_\alpha=\prod_i M_{\alpha_i}$. We consider $M_\alpha$ as being diagonally embedded in $M_{|\alpha|}$. $M_\alpha$ is the Lie algebra of $\Gl_\alpha$.
The derivative of the $\Gl_\alpha$-action on $\Rep(A,\alpha)$ 
yields an $M_\alpha$ action which has the following formula for $v\in M_\alpha(k)$: 
\begin{equation}
\label{ref-7.1-99}
v\cdot X(a)=[X(a),v]
\end{equation}

We now indicate how some of the possible properties of $A$ we have introduced 
induce standard geometrical properties on $\Rep(A,\alpha)$.

\subsection{Functions}
We have already seen that $a\in A$ induces  functions $(a_{ij})_{ij}$
on $\Rep(A,\alpha)$.
\subsection{Differential forms}
If $\omega=f_1df_2\cdots df_n\in (\Omega_B A)_n$ then we define
\begin{equation}
\label{ref-7.2-100}
\omega_{ij}=f_{1,ia_1}df_{2,a_1a_2}\cdots df_{n,a_{n-1}j}
\end{equation}
$(\omega_{ij})_{ij}$ is a matrix valued differential form on $\Rep(A,\alpha)$. 
If we write it as $X(\omega)$ then \eqref{ref-7.2-100} may be rewritten as
\[
X(\omega)=X(f_1)dX(f_2)\cdots dX(f_n)
\]
\subsection{P{{oly-vector}} fields}
If $\delta\in D_{A/B}$ then we define corresponding vector fields  $\delta_{ij}\in$
$\Rep(A,\alpha)$ by the
rule
\begin{equation}
\label{ref-7.3-101}
\delta_{ij}(a_{uv})=\delta(a)_{uj}'\delta(a)''_{iv}
\end{equation}
If $\delta=\delta_1\cdots \delta_n\in (D_BA)_n$ then we put
\[
\delta_{ij}=\delta_{1,ia_1}\delta_{2,a_1a_2}\cdots \delta_{n,a_{n-1}j}
\in \bigwedge_{\Oscr(\Rep(A,\alpha))} \Der(\Oscr(\Rep(A,\alpha)))
\]
or in the standard matrix notation
\[
X(\delta)=X(\delta_1)\cdots X(\delta_n)
\]
\subsection{Brackets}
\label{ref-7.5-102}
We have the following result. 
\begin{propositions} Assume that $\ldb -,-\rdb: A\times A\r A\otimes A$ is
a $B$-linear double bracket on $A$.
Then there is a unique antisymmetric biderivation 
\[
\{ -,-\}:\Oscr(\Rep(A,\alpha))\times \Oscr(\Rep(A,\alpha))\r
\Oscr(\Rep(A,\alpha))
\]
with the property
\begin{equation}
\label{ref-7.4-103}
\{ a_{ij},b_{uv}\}=\ldb a,b\rdb'_{uj} \ldb a,b\rdb''_{iv}
\end{equation}
for $a,b\in A$.
\end{propositions}
\begin{proof} It is a routine verification that \eqref{ref-7.4-103}
is compatible with the defining relations of $\Oscr(\Rep(A,\alpha))$.
The antisymmetry of $\{-,-\}$ may be checked on the generators
$(a_{ij})_{ij}$ where it follows from the corresponding property of $\ldb-,-\rdb$.\end{proof}
The following proposition gives the connection between the double Jacobi identity
in $A$ and the Jacobi identity on $\Rep(A,\alpha)$.
\begin{propositions} \label{ref-7.5.2-104} The following identity holds for $a,b,c\in A$.
\begin{multline}
\label{ref-7.5-105}
\{ a_{pq},\{b_{rs},c_{uv}\}\}
+
\{ b_{rs},\{c_{uv},a_{pq}\}\}
+
\{ c_{uv},\{a_{pq},b_{rs}\}\}
=\\
\ldb a,b,c\rdb'_{uq}\ldb a,b,c\rdb''_{ps}\ldb a,b,c\rdb'''_{rv}
-
\ldb a,c,b\rdb'_{rq}\ldb a,c,b\rdb''_{pv}\ldb a,c,b\rdb'''_{us}
\end{multline}
In particular, if $A,\ldb-,-\rdb$  is a double Poisson algebra then 
$\Oscr(\Rep(A,\alpha)),\{-,-\}$ is a Poisson algebra.
\end{propositions}
\begin{proof} We compute
\begin{align*}
\{ a_{pq},\{b_{rs},c_{uv}\}\}&=\{ a_{pq},
\ldb b,c\rdb'_{us} \ldb b,c\rdb''_{rv}
\}\\
&=\{ a_{pq},
\ldb b,c\rdb'_{us} \}\ldb b,c\rdb''_{rv}
+
\ldb b,c\rdb'_{us} \{ a_{pq},\ldb b,c\rdb''_{rv}\}\\
&=
\ldb a,
\ldb b,c\rdb' \rdb'_{uq} \ldb a,
\ldb b,c\rdb' \rdb''_{ps}
\ldb b,c\rdb''_{rv}
+
\ldb b,c\rdb'_{us} \ldb a,\ldb b,c\rdb''\rdb'_{rq}
\ldb a,\ldb b,c\rdb''\rdb''_{pv}\\
&=
\ldb a,
\ldb b,c\rdb' \rdb'_{uq} \ldb a,
\ldb b,c\rdb' \rdb''_{ps}
\ldb b,c\rdb''_{rv}
-
 \ldb a,\ldb c,b\rdb'\rdb'_{rq}
\ldb a,\ldb c,b\rdb'\rdb''_{pv}\ldb c,b\rdb''_{us}
\end{align*}
and hence
\begin{align*}
\{ b_{rs},\{c_{uv},a_{pq}\}\}&=
\ldb b,
\ldb c,a\rdb' \rdb'_{ps} \ldb b,
\ldb c,a\rdb' \rdb''_{rv}
\ldb c,a\rdb''_{uq}
-
 \ldb b,\ldb a,c\rdb'\rdb'_{us}
\ldb b,\ldb a,c\rdb'\rdb''_{rq}\ldb a,c\rdb''_{pv}\\
\{ c_{uv},\{a_{pq},b_{rs}\}\}&=
\ldb c,
\ldb a,b\rdb' \rdb'_{rv} \ldb c,
\ldb a,b\rdb' \rdb''_{uq}
\ldb a,b\rdb''_{ps}
-
 \ldb c,\ldb b,a\rdb'\rdb'_{pv}
\ldb c,\ldb b,a\rdb'\rdb''_{us}\ldb b,a\rdb''_{rq}
\end{align*}
Taking the sum yields \eqref{ref-7.5-105}.
\end{proof}
\begin{examples}
\label{ref-7.5.3-106}
Recall that if $\frak{g}$ is a Lie algebra then the functions on $\frak{g}^\ast$
carry a Poisson bracket defined by
\[
\{\ev_v,\ev_w\}=\ev_{[v,w]}
\]
where $v,w\in \frak{g}$ and $\ev_v$ is the evaluation of an element of $\frak{g}^\ast$ at $v$. Clearly $\ev_v$ defines a set generating functions for
$\Oscr(\frak{g}^\ast)$. 

Since $M_n(k)$ can be identified with its dual through the trace pairing it follows
that the functions on $M_n$ have a canonical Poisson bracket. On the other
hand $M_n(k)=\Rep(k[t],M_n(k))$.
 It is then easy
to show that this Poisson bracket on $\Oscr(M_n(k))$
comes from the double Poisson bracket on $k[t]$
given by
\[
\ldb t,t\rdb =t\otimes 1-1\otimes t
\]
which we considered in Example \ref{ref-2.3.3-15}.
\end{examples}
\subsection{The Schouten bracket}
The idea is that constructions in $A$ are compatible with the corresponding
constructions on $\Rep(A,\alpha)$. This is usually clear. For the Schouten
bracket it requires some work.
\begin{propositions} Let $P,Q\in D_B A$. Then
\begin{equation}
\label{ref-7.6-107}
\{ P_{ij},Q_{uv}\} =\ldb P,Q\rdb_{uj}'\ldb P,Q\rdb_{iv}''
\end{equation}
where  $\ldb-,-\rdb$ denotes the Schouten bracket on $D_B A$ and $\{-,-\}$ is
the Schouten bracket between p{{oly-vector}} fields on $\Rep(A,\alpha)$.
\end{propositions}
\begin{proof} We claim that the correctness of \eqref{ref-7.6-107}  is
  multiplicative in both arguments. To check this put first $Q=RS$ and
  assume that $P$, $R$, $S$ are homogeneous. Assume that
  \eqref{ref-7.6-107} holds for $Q=R,S$. We compute
\begin{align*}
\{ P_{ij},(RS)_{uv}\}&=\{ P_{ij},R_{uw}S_{wv}\}\\
&=\{ P_{ij},R_{uw}\} S_{wv}+(-1)^{|R|(|P|-1)}R_{uw}\{ P_{ij},S_{wv}\}
\\
&=\ldb P,R\rdb'_{uj} \ldb P,R\rdb''_{iw} S_{wv}
+(-1)^{|R|(|P|-1)}R_{uw} \ldb P,S\rdb'_{wj}\ldb P,S\rdb''_{iv}\\
&=\ldb P,R\rdb'_{uj}(\ldb P,R\rdb''S)_{iv}
+(-1)^{|R|(|P|-1)}(R \ldb P,S\rdb')_{uj}\ldb P,S\rdb''_{iv}\\
&=\ldb P,RS\rdb'_{uj}\ldb P, RS\rdb''_{iv}
\end{align*}
We now check multiplicativity in the other argument. Put $P=UV$ and
assume that \eqref{ref-7.6-107} holds with $P=U,V$.
\begin{align*}
\{ (UV)_{ij},Q_{uv}\} &=\{ U_{ik}V_{kj},Q_{uv}\}\\
&=U_{ik} \{ V_{kj},Q_{uv}\}+ (-1)^{|V|(|Q|-1)}\{ U_{ik},Q_{uv}\}
V_{kj}\\
&=U_{ik}\ldb V,Q\rdb'_{uj}\ldb V,Q\rdb''_{kv}
+ (-1)^{|V|(|Q|-1)}
\ldb U,Q\rdb'_{uk} \ldb U,Q\rdb''_{iv}
V_{kj}\\
&=\ldb V,Q\rdb'_{uj}(U\ldb V,Q\rdb'')_{iv}+ (-1)^{|V|(|Q|-1)}
(\ldb U,Q\rdb'V)_{uj}\ldb U,Q\rdb''_{iv}\\
&=\ldb UV,Q\rdb'_{uj}\ldb UV,Q\rdb''_{iv}
\end{align*}
If follows that we have check \eqref{ref-7.6-107} only on elements of $(D_BA)_i$
with $i=0,1$.

If $P,Q\in A$ then there is nothing to prove. So assume $P=\delta\in D_{A/B}$ and
$Q=a\in A$. Then we need to prove
\[
\delta_{ij}(a_{uv})=\delta(a)'_{uj}\delta(a)''_{iv}
\]
but this is precisely \eqref{ref-7.3-101}.

The case $P\in A$ and $Q\in D_B A$ follows from the previous case by antisymmetry
of both $\{-,-\}$ and $\ldb -,-\rdb$. Hence we concentrate on the final
case $P=\delta\in D_{A/B}$ and $Q=\Delta\in D_{A/B}$. Let $a$ be an arbitrary element
of $A$. We will show 
\[
\{ \delta_{ij},\Delta_{uv}\}(a_{pq}) 
=(\ldb \delta,\Delta\rdb_{uj}'\ldb \delta,\Delta\rdb_{iv}'')(a_{pq})
\]
This equation translates into
\begin{equation}
\label{ref-7.7-108}
\delta_{ij}\Delta_{uv}(a_{pq})-\Delta_{uv}\delta_{ij}(a_{pq})=\
\ldb \delta,\Delta\rdb'_{l,uj}(a_{pq}) \ldb \delta,\Delta\rdb_{l,iv}''
+\ldb \delta,\Delta\rdb'_{r,uj} \ldb \delta,\Delta\rdb_{r,iv}''(a_{pq}) 
\end{equation}
The lefthand side of this equation is obtained from the fact that 
the Schouten bracket of vector fields is the commutator. The righthand side
is obtained by writing $\ldb-,-\rdb=\ldb-,-\rdb_l+\ldb-,-\rdb_r$ and observing that $\ldb-,-\rdb_l$ takes values in $D_BA\otimes A$
and $\ldb-,-\rdb_r$ takes values in $A\otimes D_BA$.

We compute
\begin{align*}
\delta_{ij}\Delta_{uv}(a_{pq})&=\delta_{ij}(\Delta(a)'_{pv}\Delta(a)''_{uq})\\
&=\delta_{ij}(\Delta(a)'_{pv})\Delta(a)''_{uq}
+\Delta(a)'_{pv}\delta_{ij}(\Delta(a)''_{uq})\\
&=
\delta(\Delta(a)')'_{pj}\delta(\Delta(a)')''_{iv}\Delta(a)''_{uq}
+\Delta(a)'_{pv}\delta(\Delta(a)'')'_{uj}\delta(\Delta(a)'')''_{iq}
\end{align*}
and in the same way
\begin{align*}
\Delta_{uv}\delta_{ij}(a_{pq})&=
\Delta(\delta(a)')'_{pv}\Delta(\delta(a)')''_{uj}\delta(a)''_{iq}
+\delta(a)'_{pj}\Delta(\delta(a)'')'_{iv}\Delta(\delta(a)'')''_{uq}
\end{align*}
We deduce
\begin{multline*}
\delta_{ij}\Delta_{uv}(a_{pq})-\Delta_{uv}\delta_{ij}(a_{pq})\\
=
[\delta,\Delta]\,\tilde{}_l(a)'_{pj}[\delta,\Delta]\,\tilde{}_l(a)''_{iv}
[\delta,\Delta]\,\tilde{}_l(a)'''_{uq}
+
[\delta,\Delta]\,\tilde{}_r(a)'_{pv}[\delta,\Delta]\,\tilde{}_r(a)''_{uj}
[\delta,\Delta]\,\tilde{}_r(a)'''_{iq}
\end{multline*}
Now we look at the righthand side of \eqref{ref-7.7-108}.
\begin{align*}
\ldb \delta,\Delta\rdb'_{l,uj}(a_{pq}) \ldb \delta,\Delta\rdb_{l,iv}''
&=\ldb \delta,\Delta\rdb'_{l}(a)'_{pj}\ldb \delta,\Delta\rdb'_{l}(a)''_{uq}
 \ldb \delta,\Delta\rdb_{l,iv}''\\
&=[\delta,\Delta]\,\tilde{}_l(a)'_{pj}[\delta,\Delta]\,\tilde{}_l(a)''_{iv}
[\delta,\Delta]\,\tilde{}_l(a)'''_{uq}
\end{align*}
and similarly
\begin{align*}
\ldb \delta,\Delta\rdb'_{r,uj} \ldb \delta,\Delta\rdb_{r,iv}''(a_{pq}) 
&=
\ldb \delta,\Delta\rdb'_{r,uj} \ldb 
\delta,\Delta\rdb_{r}''(a)'_{pv}\ldb\delta,\Delta\rdb_{r}''(a)''_{iq}\\
&=[\delta,\Delta]\,\tilde{}_r(a)'_{pv}[\delta,\Delta]\,\tilde{}_r(a)''_{uj}
[\delta,\Delta]\,\tilde{}_r(a)'''_{iq}
\end{align*}
which finishes the proof.
\end{proof}
\subsection{Invariant functions} We leave it to the reader to check the following property.
Let $a\in A$, $\omega\in (\Omega_B A)_n$, $\delta\in (D_B A)_n$. Then
$\Tr X(a)$, $\Tr X(\omega)$, $\Tr X(\delta)$ depend only on the value
of $a,\omega,\delta$, modulo commutators.  For simplicity we write $\tr(-)=
\Tr X(-)$.

The famous Artin, Le Bruyn, Procesi theorem reformulated in this language
reads:
\begin{theorems} $\Oscr(\Rep(A,\alpha))^{\Gl_\alpha}$  is the ring generated
by the functions $\tr(a)$ for $a\in A$.
\end{theorems}
The following result was proved by Crawley-Boevey \cite{CB4}.
\begin{propositions}
\label{ref-7.7.2-109} If $A$ is equipped with a Poisson structure 
(see \ref{ref-2.6-35}) with associated Lie bracket $\{-,-\}$ then 
 $\Oscr(\Rep(A,\alpha))^{\Gl_\alpha}$ has a unique Poisson structure with the property
\[
\{\tr(a),\tr(b)\}=\tr\{\bar{a},\bar{b}\}
\]
\end{propositions}
Traces are also compatible with the Schouten bracket.
\begin{propositions} 
\label{ref-7.7.3-110}
For $P,Q\in D_B A$ one has
\[
\{\tr(P),\tr(Q)\}=\tr \{P,Q\}
\]
\end{propositions}
\begin{proof}
This is an easy computation.
\begin{align*}
\{\tr(P),\tr(Q)\}&=\{P_{ii},Q_{jj}\}\\
&=\ldb P,Q\rdb_{ji}'\ldb P,Q\rdb_{ij}''\\
&=\{P,Q\}_{ii}\\
&=\tr \{P,Q\}
\end{align*}
\end{proof}
\begin{corollarys}
The map
\[
\tr: D_B A/[D_B A,D_B A]\r \bigwedge_{\Oscr(\Rep(A,a))} \Der(\Oscr(\Rep(A,a)))
\]
is a Lie algebra homomorphism if both sides are equipped with the
Schouten bracket.
\end{corollarys}
\subsection{Compatibility}
\label{ref-7.8-111}
Assume that $P\in (D_B A)_2$. Then $P$ induces a double bracket
$\ldb-,-\rdb_P$ on $A$ and hence a corresponding antisymmetric $\{-,-\}_P$
biderivation on $\Oscr(\Rep(A,\alpha))$. On the other hand $\tr(P)$
also induces an antisymmetric biderivation on $\Oscr(\Rep(A,\alpha))$.
We claim that these are the same.  More precisely we want to show
for $f,g\in \Oscr(\Rep(A,\alpha))$ that
\[
\{f,g\}_P=\tr(P)(f,g)
\]
It suffices to
check this for $P=\delta\Delta$ with $\delta,\Delta\in D_{A/B}$.
Recall that we have for $a,b\in A$
\begin{align*}
\ldb a,b\rdb_P &=-\ldb a,\{\delta\Delta,b\}\rdb\\
&=-(\delta b)' (\Delta a)''\otimes (\Delta a)'(\delta b)''+
(\Delta b)' (\delta a)'' \otimes (\delta a)'(\Delta b)''
\end{align*}
Hence we compute
\begin{align*}
\{a_{ij},b_{uv}\}_P&=\ldb a,b\rdb_{uj}' \ldb a,b\rdb''_{iv}\\
&=-(\delta b)'_{uw} (\Delta a)''_{wj} (\Delta a)'_{iz} (\delta b)''_{zv}
+
(\Delta b)'_{uw} (\delta a)''_{wj} (\delta a)'_{iz} (\Delta b)''_{zw}\\
&=-\delta_{zw}(b_{uv}) \Delta_{wz}(a_{ij}) +\Delta_{zw} (b_{uv}) 
\delta_{wz}(a_{ij})\\
&=(\delta_{wz}\wedge \Delta_{zw})(a_{ij},b_{uv})\\
&=\tr(\delta\Delta) (a_{ij},b_{uv})
\end{align*}
\subsection{Base change}
\label{ref-7.9-112}
\begin{propositions}
\label{ref-7.9.1-113}
  Let $f_{ij}\in M_{\alpha}=\Lie(\Gl_\alpha)$ be the elementary matrix which is
  $1$ in the $(i,j)$-entry and zero everywhere else. 
Then $(E_p)_{ij}$ acts as $f_{ji}$ on $\Oscr(\Rep(A,\alpha))$
if  $\phi(i)=\phi(j)=p$ and else as zero.
\end{propositions}
\begin{proof}
  Consider first the case  $\phi(i)=\phi(j)=p$.
  The formula \eqref{ref-7.1-99} yields
\[
f_{ji}a_{uv}=[X(a),f_{ji}]_{uv}=a_{uj}\delta_{iv}-\delta_{uj} a_{iv}
\]
(here $\delta$ is the Kronecker delta). 
\begin{align*}
(E_p)_{ij}(a_{uv})&=E_p(a)'_{uj}E_p(a)''_{iv}\\
&=(ae_p)_{uj}(e_p)_{iv}-
(e_p)_{uj}(e_pa)_{iv}\\
&=a_{uj}\delta_{iv}-
\delta_{uj}a_{iv}
\end{align*}
where we have used $(ae_p)_{uj}=a_{uw}(e_p)_{wj}=a_{uw}\delta_{wj}=a_{uj}$.
If it is not true that $\phi(i)=\phi(j)=p$ then a similar computation shows
that $(E_p)_{ij}$ acts as zero.
\end{proof}
\begin{remarks} The previous proposition explains why we have called the
elements $(E_p)_p\in (D_BA)_1$ ``gauge elements'' in \S\ref{ref-3.3-53}. 
They correspond to gauge transformations on $\Oscr(\Rep(A,\alpha))$.
\end{remarks} 
\subsection{Fusion}
\label{ref-7.10-114}
In this section the notations are as in \S\ref{ref-2.5-28}. 
\begin{lemmas}
Assume that $\alpha_1=\alpha_2$ and let $\alpha^f=(\alpha_1,\alpha_3,\ldots,
\alpha_n)$. Consider $\Gl_{\alpha^f}$ as being embedded in $\Gl_\alpha$ where
the embedding on the first factor is diagonal and on the other factors the 
identity.

There is a canonical isomorphism between $\Rep(A,\alpha)$ and
$\Rep(A^f,\alpha^f)$ such that the induced $\Gl_{\alpha^f}$ action on
$\Rep(A,\alpha)$ is obtained by restriction from the
$\Gl_\alpha$-action.
\end{lemmas}
\begin{proof} Left to the reader.
\end{proof}
\subsection{Hamiltonian structure}
We have shown in Proposition \ref{ref-7.5.2-104} that if $A$ is a
double Poisson algebra then $\Oscr(\Rep(A,\alpha))$ is a Poisson
algebra. In this section we discuss the Hamiltonian structure.

  If $G$ is an algebraic group, with Lie algebra $\frak{g}$, acting on
  an (affine) Poisson
variety $X$ then a moment map for this action is by definition an invariant map 
$\psi:X\r \frak{g}^\ast$ such that for all $v\in \frak{g}$ and $f\in \Oscr(X)$ we have
\begin{equation}
\label{ref-7.8-115}
\{\langle v,-\rangle\circ \psi, f\}=v(f)
\end{equation}

\begin{propositions} Let $A,\ldb-,-\rdb ,\mu$ be a
Hamiltonian algebra. Then
\[
X(\mu_p)_p:\Rep(A,\alpha)\r M_\alpha
\]
is a moment map for $\Rep(A,\alpha)$ equipped with the associated bracket
$\{-,-\}$ (as in \S\ref{ref-7.5-102}).
\end{propositions}
\begin{proof}

We verify this \eqref{ref-7.8-115} in the case $X=\Rep(A,\alpha)$ and $\psi=X(\mu_p)_p$.
It suffices to check \eqref{ref-7.8-115} with $v=f_{ji}$ 
with $\phi(i)=\phi(j)=p$
and $f=a_{uv}$, $a\in A$.
Then \eqref{ref-7.8-115} becomes
\[
\sum_p \{\Tr(f_{ji}X(\mu_p)_p), a_{uv}\}=(E_{p})_{ij}(a_{uv})
\]
We compute the left hand side of this equation
\begin{align*}
\smash{\sum_p} \{\Tr(f_{ji}X(\mu_p)_p), a_{uv}\}&=\{\mu_{p,ij}, a_{uv}\}\\
&=\ldb\mu_{p}, a\rdb'_{uj}\ldb\mu_{p}, a\rdb''_{iv}\\
&=E_{p}(a)'_{uj}E_{p}(a)''_{iv}\\
&=(E_p)_{ij}(a_{uv})
\end{align*}
This finishes the proof.
\end{proof}
\subsection{Quasi-Poisson structure}
Let $\frak{g}$ be a Lie algebra equipped with an invariant
non-degenerate symmetric bilinear form $(-,-)$. Let $(f_a)_a$,
$(f^a)_a$ be dual bases of $\frak{g}$. Then there is a canonical
invariant element $\phi\in \wedge^3\frak{g}$ given by
\[
\frac{1}{12} c^{abc} f_a\wedge f_b \wedge f_c
\]
where 
\[
c^{abc}=(f^a,[f^b,f^c])
\]
If $G$ acts on an affine variety $X$ then we have an induced $3$-vector field
$\phi_X$ on~$X$. Following \cite{AKM} an element $P\in \bigwedge^2_{\Oscr(X)} 
\Der(\Oscr(X))$ is said to be a \emph{quasi-Poisson bracket} if 
\[
\{P,P\}=\phi_X
\]

Now we compute  $\phi$  for $M_\alpha$ with the trace pairing.
In that case $(f_a)_a=(f_{ij})_{ij}$, $(f^a)_a=(f_{ji})_{ij}$. Hence
\begin{align*}
c^{ij,kl,mn}&=\Tr(f_{ji}[f_{lk},f_{nm}])\\
&=\Tr(f_{ji}f_{lk}f_{nm}-f_{ji}f_{nm}f_{lk})\\
&=\delta_{il}\delta_{kn}\delta_{jm}-\delta_{in}\delta_{lm}\delta_{jk}
\end{align*}
We can now compute $\phi$.
\begin{align*}
12\phi&=
(\delta_{il}\delta_{kn}\delta_{jm}-\delta_{in}\delta_{lm}\delta_{jk})
f_{ij}\wedge f_{kl}\wedge f_{mn}\\
&=f_{ij}\wedge f_{ki}\wedge f_{jk}-f_{ij}\wedge f_{jl}\wedge f_{li}\\
&=2f_{ij}\wedge f_{ki}\wedge f_{jk}
\end{align*}

From Proposition \ref{ref-7.9.1-113} we obtain. 
\begin{propositions}
\label{ref-7.12.1-116}
  The three vector field on $\Rep(A,\alpha)$ induced by $\phi$ is given
  by 
\[
\frac{1}{6}\sum_i \tr(E^3_i) 
\]
\end{propositions}
We obtain 
\begin{theorems} 
\label{ref-7.12.2-117}
Assume that $A,P$ is a differential double quasi-Poisson algebra. 
Then $\tr(P)$ is a quasi-Poisson bracket on $\Rep(A,\alpha)$. 
\end{theorems}
\begin{proof} This follows by applying taking the trace of the defining property
\[
\{P,P\}=\frac{1}{6} \sum_{i=1}^n E_i^3 \quad \operatorname{mod}\ [D_BA,D_BA]
\]
(see \S\ref{ref-5.2-75}) together with Propositions \ref{ref-7.7.3-110}
and \ref{ref-7.12.1-116}.
\end{proof}
\begin{remarks} By a somewhat tedious verification using
Proposition \ref{ref-7.5.2-104}
it follows
that Theorem \ref{ref-7.12.2-117} is also true in the non-differential case.
We omit this.
\end{remarks}
\subsection{Quasi-Hamiltonian structure}
Let $G,X,\frak{g},(-,-)$ be as in the beginning of the previous section.

For $v\in \frak{g}$
 let $v^L$, $v^R$ be the associated  left and right invariant vector fields. 
According to the conventions in \cite[p2,3]{AKM}, if $g$ is a function on $G$ then
\begin{align}
v^L(g)(z)&=\frac{d\ }{dt}g(z\exp(tv))_{t=0} \label{ref-7.9-118}\\
v^R(g)(z)&=\frac{d\ }{dt}g(\exp(tv)z)_{t=0} \label{ref-7.10-119}
\end{align}
 If $v\in \frak{g}$
then $v_X$ is by definition the vector field on $X$ defined by 
\[
v_X(g)(x)=\frac{d\ }{dt}g(\exp(-tv)x)_{t=0}
\]
for a function $g$ on $X$.

\begin{definitions}\cite{AKM}
Assume that $\Oscr(X),P$ is a quasi-Poisson algebra
Let $(f^a)_a$ be $(f_a)_a$ be a pair of dual bases
for $\frak{g}$. An $\operatorname{Ad}$-equivariant map 
\[
\Phi:X\r G
\]
is a multiplicative moment map if for all functions $g$ on $G$ we have
\[
\{g\circ \Phi,-\}=\frac{1}{2} f^a_X\left((f_a^L+f_a^R)(g)\circ \Phi\right)
\]
\end{definitions}
We can now prove the following result.
\begin{propositions} Let $A,P$ be a double differential 
quasi-Poisson algebra and let $\Phi=(\Phi_p)_p\in \oplus_p e_p A e_p$ be a multiplicative moment map. Then
\[
X(\Phi_p)_p:\Rep(A,\alpha)\r M_\alpha
\]
is a multiplicative moment map for $\Rep(A,\alpha)$ equipped with the Poisson
bracket
$\tr(P)$.
\end{propositions}
\begin{proof}
 As dual bases (for the trace pairing on $M_\alpha$)
we choose $(f_{ij})_{ij}$ and $(f_{ji})_{ij}$.

We apply \eqref{ref-7.9-118} 
with $v=f_{ij}$ and $g=g_{uv}$ where
$g_{uv}$ is the projection on the $uv$'th entry of $M_{\alpha}$ and
$u,v$ are such that $\phi(u)=\phi(v)=q$. This yields
\begin{align*}
f^L_{ij}(g_{uv})(z)&=g_{uv}(zf_{ij})\\
&=\delta_{jv}z_{ui}
\end{align*}
and hence
\[
f^L_{ij}(g_{uv})=\delta_{jv} g_{ui}
\]
Similarly 
\begin{align*}
f^R_{ij}(g_{uv})(z)&=g_{uv}(f_{ij}z)\\
&=\delta_{iu}z_{jv}
\end{align*}
and hence
\[
f^R_{ij}(g_{uv})=\delta_{iu} g_{jv}
\]
From this computation we obtain
(with $X=\Rep(A,\alpha)$)
\begin{align*}
\frac{1}{2} (f_{ji})_X\left(((f_{ij})^L+(f_{ij})^R)(g_{uv})\circ \Phi\right)&=
\frac{1}{2} (f_{ji})_X((\delta_{jv} g_{ui}+ \delta_{iu} g_{jv})\circ \Phi)\\
&=\frac{1}{2}((f_{vi})_X \Phi_{q,ui}+  (f_{ju})_X    \Phi_{q,jv})\\
&=\frac{1}{2}(E_{q,iv}\Phi_{q,ui}+E_{q,uj}\Phi_{q,jv})\\
&=\frac{1}{2}(\Phi_q E_q+E_q \Phi_q)_{uv}
\end{align*}
Thus for $a\in A$:
\begin{equation}
\label{ref-7.11-120}
\frac{1}{2} (f_{ji})_X\left(((f_{ij})^L+(f_{ij})^R)(g_{uv})\circ \Phi\right)
(a_{rs})=
\frac{1}{2}(\Phi_q E_q+E_q \Phi_q)(a)'_{rv}(\Phi_q E_q+E_q \Phi_q)(a)''_{us}
\end{equation}
On the other hand
\begin{equation}
\label{ref-7.12-121}
\{g_{uv}\circ X(\Phi),a_{rs}\}=\{\Phi_{q,uv},a_{rs}\}=
\ldb\Phi_q,a\rdb'_{rv}\ldb\Phi_q,a\rdb''_{us}
\end{equation}
We obtain that \eqref{ref-7.11-120} is indeed equal to \eqref{ref-7.12-121} from the
defining identity for a multiplicative moment map.
\[
\ldb \Phi_q,a\rdb=\frac{1}{2}(\Phi_q E_q + E_q \Phi_q)(a)\qed
\]
\def\qed{}\end{proof}
\subsection{Interpretation for quivers}
It follows from Proposition \ref{ref-6.8.1-97} together with Proposition
\ref{ref-7.7.2-109} that if $A$ is either a deformed preprojective algebra or
a deformed multiplicative preprojective algebra then 
$\Oscr(\Rep(A,\alpha))^{\Gl(\alpha)}$ has a Poisson
structure. The explicit formulas for the 
Poisson bracket may be obtained from \eqref{ref-6.3-89} and \eqref{ref-6.4-92}
provided we can interpret the partial derivatives that occur.

It is easy to see that
$\Rep(kQ,\alpha)$ is the polynomial algebra with generators $a_{ij}$ for
$a\in Q$ and  $\phi(i)=h(a)$, $\phi(j)=t(a)$. It is convenient to
set $a_{ij}=0$ if this last condition is not satisfied.
\begin{lemmas}  We have
\[
\left(\frac{\partial\ }{\partial a}
\right)_{ij}=
\begin{cases}
\displaystyle \frac{\partial\ }{\partial a_{ji}}&\text{\rm if $\phi(i)=h(a)$, $\phi(j)=t(a)$}\\
0&\text{\rm otherwise}
\end{cases}
\]
\end{lemmas}
\begin{proof} By the definition of $\frac{\partial\ }{\partial a}$ it
  follows that $e_p\frac{\partial\ }{\partial a}=0$ for $p\neq h(a)$
  and $\frac{\partial\ }{\partial a}e_q=0$ for $q\neq t(a)$. From this
  it follows that 
if $\phi(i)\neq h(a)$ or $\phi(j)\neq t(a)$ 
then
  $\left(\frac{\partial\ }{\partial a} \right)_{ij}=0$. So let us
assume that $\phi(i)=h(a)$ and $\phi(j)=t(a)$.

We have for $a,b\in Q$
\[
\left(\frac{\partial\ }{\partial a}\right)_{ij}(b_{uv})
=\left(\frac{\partial b}{\partial a}\right)'_{uj}
\left(\frac{\partial b}{\partial a}\right)''_{iv}
\]
If $a\neq b$ then we obtain 
\[
\left(\frac{\partial\ }{\partial a}\right)_{ij}(b_{uv})=0
\]
So assume $b=a$. Then
\begin{equation}
\label{ref-7.13-122}
\left(\frac{\partial\ }{\partial a}\right)_{ij}(a_{uv})
=(e_{t(a)})_{uj}(e_{h(a)})_{iv} 
\end{equation}
if $\phi(u)\neq t(a)$ or $\phi(v)\neq h(a)$
then both sides of  \eqref{ref-7.13-122} are zero. So let us assume
that $\phi(u)=t(a)$ or $\phi(v)= h(a)$. Then \eqref{ref-7.13-122} becomes
\begin{align*}
\left(\frac{\partial\ }{\partial a}\right)_{ij}(a_{uv})
&=\delta_{uj}\delta_{iv}\\
&=\frac{\partial a_{uv}}{\partial a_{ji}}\qed
\end{align*}
\def\qed{}\end{proof}
In the case of the deformed preprojective algebra we obtain the
classical result that the Poisson bracket corresponds to the bi-vector field
\[
\sum_{a\in Q} \frac{\partial\ }{\partial a_{ij}}
\frac{\partial\ }{\partial a^\ast_{ji}}
\]
For the deformed multiplicative preprojective algebra we obtain 
(using \eqref{ref-6.4-92}) a similar but more complicated Poisson bracket.

\appendix
\section{Relation to the theory of bi-symplectic forms}
In this appendix we relate our theory of double Poisson brackets to
the theory of bi-symplectic forms introduced in \cite{CBEG}.  The
analogous, but more involved theory for double quasi-Poisson brackets 
will be deferred to a separate note. 

We assume as usual that $A/B$ is finitely generated.  Let $\Omega_B A$
be the tensor algebra over $A$ of $\Omega_{A/B}$. This is a
DG-algebra. Assume that $A$ is equipped with a $B$-linear
bi-symplectic form $\omega$ (see Definition \ref{ref-A.3.1-134} below). We prove:
\begin{enumerate}
\item The Lie bracket on $A/[A,A]$ associated to $\omega$
  \cite[Prop.\ 4.4.1]{CBEG} comes from a $B$-linear 
  double differential Poisson bracket $P$ on $A$ (and hence, by Proposition
  \ref{ref-2.4.4-25}, from the structure of a left Loday algebra on $A$).
\item  
The algebras $\Omega_B A$
  and $D_B A$ become isomorphic DG-algebras if we equip $D_B A$ with
  the differential $-\{P,-\}$.
\end{enumerate}
The formalism we will outline is remarkably similar to the commutative
case. For example in Theorem \ref{ref-A.6.1-146} below we prove that
the condition $d\omega=0$ for a bi-symplectic form is precisely equivalent
to the condition $\{P,P\}=0$ for the corresponding double Poisson bracket.

\subsection{Differentials and double derivations}
We recall some definitions from \cite{CBEG}. We also give some
properties which we will need afterward.

Let $\delta\in D_{A/B}$. Then
we may define double derivations
\begin{align*}
i_\delta&:\Omega_B A\r \Omega_BA \otimes \Omega_B A\\
L_\delta&:\Omega_B A\r \Omega_BA \otimes \Omega_B A
\end{align*}
in the usual way: for $a\in A$ define
\begin{align*}
i_\delta(a)&=0& i_\delta(da)&=\delta(a)\\
L_\delta(a)&=\delta(a) & L_\delta(da)&=d(\delta(a))
\end{align*}
where here and below we use the convention that $d$ acts on tensor products
by means of  the usual Leibniz rule. 

If $C$ is a graded $k$-algebra and $c=c_1\otimes c_2$ then we put
\[
{}^\circ c=(-1)^{|c_1||c_2|}c_2c_1
\]
and if  $\phi:C\r C^{\otimes 2}$ is a linear map then we define
\[
{}^\circ \phi:C\r C: c\mapsto {}^\circ(\phi(c))
\]
If $\delta$ is a double derivation then ${}^\circ \delta$ vanishes on
commutators.

We apply this with  $C=\Omega_B A$. 
Following \cite{CBEG} we put
\begin{align*}
\imath_\delta&={}^\circ\! i_\delta\\
\Lscr_\delta&={}^\circ\! L_\delta
\end{align*}
Now we discuss some commutation relations between these operators.
Checking on the generators $a\in A$ and $da\in \Omega_{A/B}$ of
$\Omega_B A$ we find the usual Cartan formula \cite[eq.\
(2.7.2)]{CBEG}
\[
L_\delta=di_\delta+i_\delta d
\]
from which one obtains by applying the operation ${}^\circ(-)$ \cite[Lemma
2.8.8(i)]{CBEG}
\begin{equation}
\label{ref-A.1-123}
\Lscr_\delta=d\imath_\delta+\imath_\delta d
\end{equation}
As the $(i_\delta)_\delta$ are double derivations one can take their
Schouten brackets.  One has for $\delta,\Delta\in D_{A/B}$.
\begin{equation}
\label{ref-A.2-124}
\ldb i_\delta, i_\Delta\rdb_l=\ldb i_\delta, i_\Delta\rdb_r=0
\end{equation}
To see this note that both $\ldb i_\delta, i_\Delta\rdb\tilde{}_l$ and
$\ldb i_\delta, i_\Delta\rdb\tilde{}_r$ are derivations $\Omega_B A\r
(\Omega_B)^{\otimes_B A}$ of degree $-2$. Hence they must vanish on
$A$ and $\Omega_{A/B}$. This means that they must vanish on the whole
of $\Omega_B A$.

We will also use the following identities 
\begin{equation}
\label{ref-A.3-125}
\begin{aligned}
\ldb i_\delta,L_\Delta\rdb_l&=i_{\ldb\delta,\Delta\rdb'_l}\otimes {\scriptstyle \ldb\delta,\Delta\rdb''_l}
\\
\ldb i_\delta,L_\Delta\rdb_r&={\scriptstyle \ldb\delta,\Delta\rdb'_r}\otimes i_{\ldb\delta,\Delta\rdb''_r}
\end{aligned}
\end{equation}
which are proved by checking them on the generators of $\Omega_B A$. 

Let $C$ be a $k$-algebra and $\delta,\Delta$ be double $C$-derivations.
The straightforward proof of the next formula is left to the
reader.
\begin{equation}
\label{ref-A.4-126}
\delta \circ {}^\circ \!\Delta-\tau_{12}\circ\Delta \circ {}^\circ \!\delta
={}^{\circ,l} \ldb \delta,\Delta\rdb_l +
{}^{\circ,r} \ldb \Delta,\delta\rdb_r
\end{equation}
where
${}^{\circ,l}\!(\epsilon'\otimes\epsilon'')={}^{\circ}\!\epsilon'\otimes
\epsilon''$ and
${}^{\circ,r}\!(\mu'\otimes\mu'')=\mu'\otimes
{}^{\circ}\!\mu''$.

From the graded version of \eqref{ref-A.4-126} together with
\eqref{ref-A.2-124} we obtain the following formula
\begin{equation}
\label{ref-A.5-127}
i_\delta \imath_{\Delta}+\sigma_{12} i_\Delta \imath_{\delta}=0
\end{equation}
Combining the  graded version of \eqref{ref-A.4-126} with \eqref{ref-A.3-125}
we find
\begin{equation}
\label{ref-A.6-128}
i_\delta\Lscr_{\Delta}-\sigma_{12} L_{\Delta} \imath_{\delta}=
\imath_{\ldb\delta,\Delta\rdb'_l}\otimes {\scriptstyle \ldb\delta,\Delta\rdb''_l}
+ {\scriptstyle \ldb\delta,\Delta\rdb'_r}\otimes\imath_{\ldb\delta,\Delta\rdb''_r}
\end{equation}

Finally assume that $\delta$ is an inner double $A$-derivation of the
form $[b,-]$ for $b\in B\otimes_k B$. Then for
all $\eta\in \Omega_B A$ one has
\begin{equation}
\label{ref-A.7-129}
\Lscr_\delta \eta=0
\end{equation}
To prove this note first that
\begin{equation}
\label{ref-A.8-130}
L_\delta \eta=b\eta-\eta b
\end{equation}
As usual this is checked on generators. \eqref{ref-A.7-129} follows immediately
from \eqref{ref-A.8-130}. 
\subsection{The Koszul bracket}
Assume that $\ldb-,-\rdb$ is a $B$-linear double bracket on~$A$. One has
the following result.
\begin{propositions}
There is a unique double bracket $\ldb-,-\rdb^{\Omega_B A}$ of degree $-1$
on $\Omega_B A$ commuting with $d$ which satisfies for $a,b\in A$
\[
\ldb da,b\rdb^{\Omega_B A} =\ldb a,b\rdb
\]
If $\ldb-,-\rdb$ is Poisson then so is $\ldb-,-\rdb^{\Omega_B A}$. 
\end{propositions}
\begin{proof} The various asserted properties may be checked on
generators. We leave the full proof to the reader.
\end{proof}
Following the commutative case we call $\ldb-,-\rdb^{\Omega_B A}$ the \emph{Koszul
bracket} associated to $\ldb-,-\rdb$.

For use below we give some formulas which are easy consequences of the
definition
\begin{equation}
\label{ref-A.9-131}
\begin{aligned}
\ldb a,b\rdb^{\Omega_B A}&=0\\
\ldb da,b\rdb^{\Omega_B A}&=\ldb a,db\rdb^{\Omega_B A}=\ldb a,b\rdb\\
\ldb da,db\rdb^{\Omega_B A}&=d\ldb a,b\rdb
\end{aligned}
\end{equation}

\begin{propositions} \label{ref-A.2.2-132} Assume that $\ldb-,-\rdb$ is a double bracket on
  $A$. Then there is a well defined map of graded $A$-algebras
\[
\Sigma:\Omega_B A\r D_B A: da\mapsto H_a
\]
(see \S\ref{ref-3.5-58} for notation). If $\ldb-,-\rdb$ is Poisson then
this map intertwines the Koszul bracket on the left with the Schouten bracket
on the right.
\end{propositions}
\begin{proof} That $\Sigma$ is well defined is easy to see by checking
on generators. To prove that $\Sigma$ is compatible with the brackets
we have to show that the following analogues of \eqref{ref-A.9-131}
hold in $D_BA$.
\begin{equation}
\label{ref-A.10-133}
\begin{aligned}
\ldb a,b\rdb^{D_BA}&=0\\
\ldb H_a,b\rdb^{D_B A}&=\ldb a,H_b\rdb^{D_B A}=\ldb a,b\rdb\\
\ldb H_a,H_b\rdb^{D_B A}&=H_{\ldb a,b\rdb}
\end{aligned}
\end{equation}
where to avoid confusion we have denoted the Schouten bracket by $\ldb
-,-\rdb^{D_B A}$. The first  equation is obvious. The second equation
follows from \eqref{ref-3.11-59}. The third equation follows
from Proposition \ref{ref-3.5.1-60}.
\end{proof}
\subsection{Bi-symplectic forms}
 Following \cite{CBEG} we
put $\DR_B(A)=\Omega_B A/[\Omega_B A, \Omega_B A]$. The differential
on $\Omega_B A$ descends to a differential on $\DR_B(A)$.

\begin{definitions} \label{ref-A.3.1-134} \cite{CBEG} 
An element $\omega\in \DR_B(A)_2$ is \emph{bi-non-degenerate} if 
the map of $A$-bimodules
\[
\imath(\omega):D_{A/B}\r \Omega_{A/B}:\delta\mapsto \imath_\delta \omega
\]
is an isomorphism. If in addition $\omega$ is closed in $\DR_B(A)$ then
we say that $\omega$ is \emph{bi-symplectic}.
\end{definitions}
Assume that $\omega\in \DR_B(A)_2$ is a bi-non-degenerate.
Let $a\in A$. Following \cite{CBEG} we define the Hamiltonian vector
field  $H_a\in D_{A/B}$ corresponding to $a$ via
\[
\imath_{H_a} \omega =da
\]
and we put
\[
\ldb a,b\rdb_\omega=H_a(b)
\]
Since  $H_a(b)=i_{H_a}(db)$ we may also write
\begin{equation}
\label{ref-A.11-135}
\ldb a,b\rdb_\omega=i_{H_a}\imath_{H_b}(\omega)
\end{equation}
\begin{lemmas} $\ldb a,b\rdb_\omega$ is a double bracket on $A$. 
\end{lemmas}
\begin{proof} It is clear that  $\ldb a,b\rdb_\omega$ a derivation in
its second argument. So we only need to prove anti-symmetry.  This follows
  immediately from \eqref{ref-A.11-135} and \eqref{ref-A.5-127}.
\end{proof}

According to \cite[Prop.\ 4.4.1]{CBEG} if $\omega$ is bi-symplectic
then the associated simple bracket $\{-,-\}_\omega$ induces a Lie
algebra structure on $A/[A,A]$. We will prove the following stronger
result.
\begin{propositions} \label{ref-A.3.3-136} Assume that $\omega\in (\Omega_B A)_2$ is
  bi-symplectic.  Then $\ldb -,-\rdb_\omega$ is a double Poisson
  bracket on $A$.
\end{propositions}
\begin{proof}
 To simplify the notations we write $\ldb-,-\rdb$ for
  $\ldb-,-\rdb_\omega$. 
According to Proposition \ref{ref-3.5.1-60} we have to prove
\[
\ldb H_a, H_b\rdb_l= H_{\ldb a,b \rdb'} \otimes \ldb a,b \rdb''
\]
Appying $\imath(\omega)\otimes 1$ this is equivalent to
\[
\imath_{\ldb H_a, H_b\rdb_l'}(\omega)\otimes \ldb H_a, H_b\rdb_l''=
d\ldb a, b\rdb_l'\otimes \ldb a, b\rdb_l''
\]
This would follow by projecting the following identity in
$\Omega_{A/B} \otimes A\oplus A\otimes \Omega_{A/B}$ onto 
the first factor
\[
\imath_{\ldb H_a, H_b\rdb_l'}(\omega)\otimes \ldb H_a, H_b\rdb_l''
+\ldb H_a, H_b\rdb_r'\otimes \imath_{\ldb H_a, H_b\rdb_r''}(\omega)=
d\ldb a, b\rdb
\]

Applying \eqref{ref-A.6-128} we find that this is equivalent to
\begin{equation}
\label{ref-A.12-137}
(i_{H_a} \Lscr_{H_b}-\tau_{12} L_{H_b} \imath_{H_a})(\omega)=d\ldb a,b\rdb
\end{equation}
Now
\[
\imath_{H_a}(\omega)=da
\]
and 
\[
\Lscr_{H_b}(\omega)=d\imath_{H_b}(\omega)+\imath_{H_b}d\omega=ddg=0
\]
(this is the place where the assumption that $\omega$ is closed is used). 

It follows that the left hand side of \eqref{ref-A.12-137} is equal
to 
\[
-\tau_{12} L_{H_b}(da)=-\tau_{12}d\ldb b,a\rdb=d\ldb a, b\rdb
\]
which is indeed equal to the righthand side of \eqref{ref-A.12-137}.
\end{proof}

\subsection{Duality yoga for bimodules}
Below $A$ is a $k$-algebra and $M$ is an $A$-bimodule. Later we take
$M=\Omega_B A$ but in this section this is not necessary. 

If $m\in M$ then there is a double derivation of degree $-1$
\[
i_m:T_A(M^\ast)\r T_A(M^\ast)\otimes T_A(M^\ast)
\]
which on $\sigma\in M^\ast$ is given by
\[
i_m(\sigma)=\sigma(m)''\otimes \sigma(m)'
\]
Similarly for $\sigma\in M^\ast$ there is an associated double derivation
of degree $-1$ 
\[
i_\sigma:T_A(M)\r T_A(M)\otimes T_A(M)
\]
with the property 
\[
i_\sigma(m)=\sigma(m)
\]
for $m\in M$. As before we put $\imath_m={}^\circ\! i_m$ and
$\imath_\sigma={}^\circ\! i_\sigma$.
\begin{propositions}  
\label{ref-A.4.1-138} 
Assume there is an element $\omega\in M\otimes_A M$
such that the induced map
\[
\imath(\omega):M^\ast \r M: \sigma\mapsto \imath_\sigma \omega
\]
is an isomorphism. Define $P=-(\imath(\omega)^{-1}\otimes
\imath(\omega)^{-1})(\omega)$. Then one has for $m\in M$
\[
\imath(P)=\imath(\omega)^{-1}
\]
where $\imath(P)$ is the map
\[
\imath(P):M\r M^\ast: m\mapsto \imath_m P
\]
\end{propositions}
\begin{proof}  Put $\theta=\imath(\omega)$ and $\psi=\imath(\omega)^{-1}$.  We
have to prove $\psi=\imath(P)$.
We have explicitly for $\sigma \in M^\ast$
\begin{equation}
\label{ref-A.13-139}
\theta(\sigma)=\sigma(\omega')'' \omega'' \sigma(\omega')'
-\sigma(\omega'')''\omega' \sigma(\omega'')'
\end{equation}
and hence for $\tau\in M^\ast$
\[
\tau(\theta(\sigma))=
\sigma(\omega')'' \tau(\omega'') \sigma(\omega')'
-\sigma(\omega'')''\tau(\omega') \sigma(\omega'')'
\]
We deduce
\begin{equation}
\label{ref-A.14-140}
\tau(\theta(\sigma))=-\sigma(\theta(\tau))^\circ
\end{equation}
 Put
$\tau=\psi(n)$ and $\sigma=\psi(m)$. We deduce from \eqref{ref-A.14-140}:
\begin{equation}
\label{ref-A.15-141}
\psi(n)(m)=-\psi(m)(n)^\circ
\end{equation}
Taking into account that $P=-\psi(\omega')\otimes\psi(\omega'')$
we find
\[
\imath_m P=-\psi(\omega')(m)' \psi(\omega'')\psi(\omega')(m)''+
\psi(\omega'')(m)'\psi(\omega')\psi(\omega'')(m)''
\]
Hence
\begin{align*}
\theta(\imath_m P)&=-\psi(\omega')(m)' \omega''\psi(\omega')(m)''
+\psi(\omega'')(m)'\omega'\psi(\omega'')(m)''\\
&=-\theta(\psi(-)(m)^\circ)
\end{align*}
where in the last line we have used \eqref{ref-A.13-139}. Hence for $n\in M$
\[
(\imath_m P)(n)=-\psi(n)(m)^\circ=\psi(m)(n)
\]
where we have used \eqref{ref-A.15-141}.
\end{proof}
\subsection{Compatibility of brackets}
\label{ref-A.5-142}
Let $\omega\in (\Omega_BA)_2$ be a bi-non-degenerate form
Putting $M=\Omega_{A/B}$ and $M^\ast=D_{A/B}$ we define  $P\in (D_B A)_2$ 
\[
P=-(\imath(\omega)^{-1}\otimes \imath(\omega)^{-1})(\omega)
\]
as in the previous section. 
\begin{propositions} One has 
\[
\ldb -,-\rdb_P=\ldb-,-\rdb_\omega
\]
\end{propositions}
\begin{proof}  Using Proposition \ref{ref-A.4.1-138} one has
\[
\ldb f,g\rdb_\omega=H_f(g)=\imath(P)(df)(g)=(\imath_{df}P)(g)
\]
The result now follows from Lemma \ref{ref-A.5.2-143} below.
\end{proof}
\begin{lemmas}
\label{ref-A.5.2-143}
 Let $P$ be an arbitrary element of $(D_BA)_2$. Then
\[
\ldb f,g\rdb_P=(\imath_{df}P)(g)
\]
\end{lemmas}
\begin{proof} It is sufficient to check this for $P=\delta\Delta$ with
$\delta,\Delta\in D_{A/B}$. We have
\[
\ldb f,g \rdb_{\delta\Delta}=
\Delta(g)'\delta(f)''\otimes \delta(f)' \Delta(g)''
-
\delta(g)'\Delta(f)''\otimes \Delta(f)' \delta(g)''
\]
and
\begin{align*}
i_{df}(\delta\Delta)&=i_{df}(\delta)\Delta-\delta i_{df}(\Delta)\\
&=\delta(f)''\otimes \delta(f)'\Delta-\delta\Delta(f)''\otimes \Delta(f)'
\end{align*}
Thus
\[
\imath_{df}(\delta\Delta)=\delta(f)'\Delta
\delta(f)''-\Delta(f)'\delta \Delta(f)''
\]
and hence
\[
\imath_{df}(\delta\Delta)(g)=\Delta(g)'\delta(f)''\otimes \delta(f)' \Delta(g)''
-
\delta(g)'\Delta(f)''\otimes \Delta(f)' \delta(g)''
\]
finishing the proof. 
\end{proof}
The following result will be used below.
\begin{lemmas} \label{ref-A.5.3-144} Let $\ldb-,-\rdb$ be the double bracket associated to
a bi-non-degenerate form $\omega$ and let $P$ be the element of $(D_B A)_2$
associated to $\omega$ (as in \S\ref{ref-A.5-142}). Let $\Sigma:\Omega_B A\r D_B A$ as in
Proposition \ref{ref-A.2.2-132}. Then $\Sigma$ is an isomorphism of
$A$-algebras and furthermore
\begin{equation}
\label{ref-A.16-145}
\Sigma(\omega)=-P
\end{equation}
\end{lemmas}
\begin{proof}
Let $\overline{\imath(P)}$ be the $A$ algebra morphism obtained by extending
$\imath(P)$. We claim that $\Sigma=\imath(P)$. For $a\in A$ we need
$\Sigma(da)=\imath(P)(da)$. By definition $\Sigma(da)=H_a$ and
$\imath(P)(da)=\imath_{da}(P)=H_a$ by Lemma \ref{ref-A.5.2-143}. In particular
$\Sigma$ is an isomorphism of algebras.
We also deduce
\[
\Sigma(\omega)=(\imath(P)\otimes \imath(P))(\omega)=(\imath(\omega)^{-1}
\otimes \imath(\omega)^{-1})(\omega)=-P\qed
\]
\def\qed{}\end{proof}
\subsection{The relation between Poisson brackets and bi-symplectic forms}
In this section we prove the following result.
\begin{theorems}
\label{ref-A.6.1-146}
Assume that $\omega\in (\Omega_B A)_2$ is a bi-non-degenerate form and let
$P$ be the corresponding element of $D_B A$. Then the following are
equivalent: 
\begin{enumerate}
\item $\omega$ is bi-symplectic, i.e. 
\[
d\omega=0
\]
in $\DR_B(A)$.
\item $P$ is differential Poisson, i.e.
\[
 \{P,P\}=0
\]
in $D_BA/[D_B A,D_BA]$. 
\end{enumerate}
\end{theorems}
\begin{proof}
 Let $\ldb-,-\rdb=
\ldb-,-\rdb_\omega=\ldb-,-\rdb_P$ be the corresponding bracket on $A$ and
let 
\[
\Sigma:\Omega_B A\r D_B A
\]
be as in Proposition \ref{ref-A.2.2-132}.  Below we assume
that either $d\omega=0$ or $\{P,P\}=0$. It follows by Proposition
\ref{ref-A.3.3-136} and Theorem \ref{ref-4.2.3-69} that in both these cases
$\ldb-,-\rdb$ is Poisson. Hence by Proposition \ref{ref-A.2.2-132}
$\Sigma$ intertwines the Koszul bracket and the Schouten bracket.

Assume first that $P$ is Poisson. 
 We claim that the following
diagram is commutative
\begin{equation}
\label{ref-A.17-147}
\begin{CD}
\Omega_B A@>\Sigma >> D_B A\\
@VdVV @VV -\{P,-\} V\\
\Omega_B A @>\Sigma >> D_B A
\end{CD}
\end{equation}
It is sufficient to check this on generators. Let $a\in A$. We have
\[
(\Sigma\circ d)(a)=\Sigma(da)=H_a
\]
and
\[
-(\{P,-\}\circ \Sigma)(a)=-\{P,a\}
\]
We need to see $\{P,a\}=-H_a$. 
Evaluating on $b\in A$ it is sufficient to prove
\begin{equation}
\label{ref-A.18-148}
\{P,a\}(b)=-\ldb a,b\rdb_P
\end{equation}
Using Proposition \ref{ref-4.2.1-67} we find  that the left hand side 
of \eqref{ref-A.18-148} is equal to 
\[
-\ldb b,\{ P, a\} \rdb^\circ=\ldb b,a\rdb^\circ_P=-\ldb a,b\rdb_P
\]
This is equal to the right hand side of \eqref{ref-A.18-148}.

Now we consider the generators of the type $da$, $a\in A$. We have
\[
(\Sigma\circ d)(da)=0
\]
and 
\[
-(\{P,-\}\circ \Sigma)(da)=-\{P,H_a\}=-\{P,\{P,a\}\}=-(1/2)\{\{P,P\},a\}\}=0
\]
Applying the identiy $\Sigma\circ d=-\{P,-\}\circ d$ to $\omega$ we conclude
\[
\Sigma(d\omega)=-\{P,\Sigma(\omega)\}=\{P,P\}=0
\]
where we have used Lemma \ref{ref-A.5.3-144}.  Since $\Sigma$ is an
isomorphism it follows $d\omega=0$.

\medskip

Now assume that $d\omega=0$. We  claim $d=\{\omega,-\}^{\Omega_B A}$. 
Let $a\in A$. First we need to check $da=\{\omega,a\}^{\Omega_B A}$. Applying
$\Sigma$ this is equivalent to $H_a=-\{ P,a\}$. This we have checked above.

Now we need to check $\{\omega,da\}^{\Omega_B A}=0$. This is the following verification. 
\[
\{\omega,da\}^{\Omega_B A}=d\{\omega,a\}^{\Omega_B A}-\{d\omega,a\}^{\Omega_B A}=dda=0
\]
From $d=\{\omega,-\}^{\Omega_B A}$ we deduce $0=d\omega=\{\omega,\omega\}^{\Omega_B A}$.
Applying $\Sigma$ we obtain $\{P,P\}=0$ (using Lemma \ref{ref-A.5.3-144}
again), finishing the proof.
\end{proof}
 We obtain
the following consequence of \eqref{ref-A.17-147}. 
\begin{corollarys}
Assume that $A$ is has a bi-symplectic form $\omega$ with
corresponding differential Poisson bracket $P$. If we equip $D_B A$ with
the differential $-\{P,-\}$ then $\Omega_B A$ and $D_B A$ become
isomorphic double DG-Gerstenhaber algebras. 
\end{corollarys}
The \emph{De Rham cohomology} of $A$ is the cohomology of $\DR_B(A)$
for the differential~$d$. Likewise the \emph{Poisson cohomology} of
$(A,P)$ \cite{weyer} is the cohomology of the complex $D_B/[D_B A,D_B
A]$ for the differential $\{P,-\}$. We obtain
the following corollary. 
\begin{corollarys}
If $A$ is equipped with a bi-symplectic form then its Poisson cohomology 
coincides with its De Rham cohomology. 
\end{corollarys}
This is an analogue of a well-known commutative result.
\subsection{The moment map}
We keep the assumptions of the previous section. We assume in addition
that $B=ke_1+\cdots+ke_n$.\footnote{In \cite{CBEG} it is only assumed
  that $B$ is semi-simple. The role of the $(e_i)_i$ is played by a separability
idempotent.} In \cite{CBEG} the very beautiful
observation is made that if $A$ is equipped with a bi-symplectic form
$\omega$ then $A$ is automatically Hamiltonian in a suitable sense
(and hence this is true for all representation spaces).  For
completeness we give the construction of the moment map in our
present context. No originality is intended. 

According to Definition \ref{ref-2.6.4-40} we need to find $\mu_i$ such
that for all $a$ one has
\[
\ldb \mu_i,a\rdb_\omega=E_i(a)
\]
Or in other words
\[
H_{\mu_i}=E_i
\]
Applying $\imath(\omega)$ this is equivalent to
\[
d\mu_i=\imath_{E_i}(\omega)
\]
Since $\Omega_B A$ is acyclic in degrees $\ge 1$ (see \cite[\S
2.5]{CBEG}) the existence of $\mu_i$ follows from
$d\imath_{E_i}(\omega)=0$. By \eqref{ref-A.1-123} and the fact that
$\omega$ is closed we have
\[
d\imath_{E_i}(\omega)=\Lscr_{E_i}(\omega)=0
\]
where we have also used \eqref{ref-A.7-129}.


\def\cprime{$'$} \def\cprime{$'$}
\ifx\undefined\bysame
\newcommand{\bysame}{\leavevmode\hbox to3em{\hrulefill}\,}
\fi

\end{document}